\documentclass[11pt]{amsart}
\usepackage{verbatim, amscd, epsfig
}

\theoremstyle{plain}
\newtheorem{theorem}{Theorem}[section]
\newtheorem{lemma}[theorem]{Lemma}
\newtheorem{corollary}[theorem]{Corollary}

\newtheorem{axiom}{Axiom}

\theoremstyle{definition}
\newtheorem{definition}{Definition}[section]

\newtheorem*{conjecture*}{Conjecture}
\newtheorem*{definition*}{Definition}
\newtheorem{example}{Example}[section]
     
\theoremstyle{remark}

\newtheorem{remark}{Remark}

\DeclareMathOperator{\rank}{rank}

\DeclareMathOperator{\coker}{coker}
\DeclareMathOperator{\ord}{ord}
\DeclareMathOperator{\diag}{diag}
\DeclareMathOperator{\trace}{tr}
\renewcommand{\Re}{{\rm Re}\,}

\newcommand{\del}{\partial}
\newcommand{\delbar}{\bar{\del}}

\newcommand{\R}{\mathbb{ R}}
\newcommand{\C}{\mathbb{ C}}
\newcommand{\Z}{\mathbb{ Z}}
\renewcommand{\H}{\mathbb{ H}}

\newcommand{\N}{\mathbb{ N}}
\renewcommand{\P}{\mathbb{ P}}

\newcommand{\tL}{\tilde{L}}

\begin{document}
\title[Quaternionic holomorphic geometry]{Quaternionic holomorphic geometry:
Pl\smash{\"u}cker formula, Dirac eigenvalue estimates and energy estimates of
harmonic 2-tori}  
\author{D. Ferus, K. Leschke, F. Pedit and U. Pinkall}
%\author{The usual suspects}
\thanks{Authors supported by SFB288. Third author additionally
  supported by  NSF-grants DMS-9011083 and  DMS-9705479}
%\thanks{Most of the usual suspects are either supported by SFB288 or NSF-grants DMS-9011083, DMS-9705479}
\maketitle

\section{Introduction}\label{sec:introduction}
In the middle of the 19th century two important theories emerged: the first,
developed by Gauss, dealt with the fundamental equations of
surfaces in 3-space. The second, developed by Riemann, was concerned
with the properties of complex curves. The latter
evolved into what is now known as the
algebraic geometry of curves with  nontrivial examples,
powerful methods and beautiful global results.
In contrast to this, the former still looks rather immature even today.
After intensive work around the turn of the 19th century by Darboux, Bianchi and
later Blaschke -- who mainly studied the local aspects of
the theory -- the subject of surfaces in 3-space lost its prominence. 
The second half of the 20th century saw a renewed interest in
some of the global aspects of the theory.
Special surface classes 
characterized by curvature properties and variational equations 
were studied, including
minimal, constant curvature and Willmore surfaces. Despite this, 
even now the catalogue of explicit examples of compact surfaces in 3-space is
very small compared to the wealth of explicitly studied algebraic curves.

One of the reasons why algebraic curve theory apparently outperformed 
surface theory lies in the latter's analytic difficulty: the fundamental 
equation of algebraic curve theory is the linear, first order Cauchy-Riemann equation, 
whereas the fundamental equations of a surface in 3-space, the Gauss-Codazzi equations,
are a nonlinear, third order system. While 
it is possible to write down explicit formulas for meromorphic functions on a given 
Riemann surface, it is almost impossible
to write down an explicit conformal parameterization into 3-space of the same Riemann
surface, even though such a parameterization exists \cite{Garcia}.
The reader may want to list conformally parameterized surfaces of, say, genus two. 

The analytic differences of the two theories 
also influenced their respective methods: algebraic curve theory  
can be formulated in the language of holomorphic line bundles. The link to extrinsic
curve theory is given by the Kodaira embedding. The basis for much of the
theory consists of fundamental results such 
as the Riemann-Roch Theorem, the Clifford estimate, the Pl\"ucker relations and the Abel map.
Surface theory in 3-space has an entirely different flavor: there are no comparable
global theorems and the Gauss-Codazzi equations for various classes of surfaces give
rise to qualitatively different systems of differential equations, 
each demanding its own theory. These equations are usually formulated
in coordinates or moving frames, which often has the effect that the initial 
geometric properties get encoded in coordinate dependent quantities, resulting
in non-geometric equations. In the rare cases where one can find
first and second fundamental forms solving the Gauss-Codazzi equations globally,
it is still very hard to control the periods of the resulting surface
in 3-space. Other global properties, such as the homotopy type, are
almost impossible to detect from the infinitesimal data. 

A related issue is the description of surfaces whose defining properties are invariant
under the group of M\"obius transformations of 3-space, such as 
isothermic surfaces (admitting conformal curvature line
parameterizations) or Willmore surfaces 
(extremizing $\int H^2$ where $H$ is the mean curvature). Due to the larger
symmetry group -- compared to the euclidean case -- 
one expects fewer invariants and thus simpler equations. Unfortunately,
most descriptions of M\"obius invariant
surface geometries are given in terms of euclidean quantities, which obscure
the inherent symmetry. If M\"obius invariant quantities are used they tend to be
dependent, resulting in over determined systems of equations.

The theme of the present paper is the unification of 
algebraic curve theory and differential geometry of surfaces. 
This unified theory is immediately applicable to surface theory and 
addresses a number of its above listed shortcomings. The theory of algebraic curves
appears as a special, some might say singular, case and all of its basic constructions, 
methods and results acquire new meaning when viewed from the
unified perspective. For instance, the
classical Pl\"ucker formula of a holomorphic
curve becomes a more fundamental relation in this unified setup, with applications
to Dirac eigenvalue estimates and to estimates on the energy of harmonic 2-tori.

To guide the reader through our exposition, it 
will be helpful to explain the fundamental ideas underlying it.
The basic notion in the theory of holomorphic curves is that of a 
meromorphic function $f$ (or holomorphic map onto the Riemann sphere $\C\P^1$) 
on a Riemann surface $M$. In differential geometric language a meromorphic function
is a branched, conformal immersion and as such satisfies the conformality condition
%\begin{equation}\label{eq:conf_imm_plane}
\[
df(JX)=i\,df(X)\,,
\]
%\end{equation}
where $J\colon TM\to TM$ is the complex structure and $X$ is a tangent vector to $M$. Using
the Hodge-star operator to denote pre-composition by the complex 
structure $J$, the above relation takes the more familiar form
%\begin{equation}\label{eq:Cauchy_Riemann}
\[
*df=i\,df
\]
%\end{equation}
of the Cauchy-Riemann equation, 
expressing the fact that the $(0,1)$-part of $df$ vanishes. On the other hand,
differential geometry of surfaces studies conformal immersions $f$ of a Riemann
surface  $M$ into $\R^3$. In case the image of $f$ lies in a plane inside $\R^3$,
admittedly a rather singular case from the viewpoint of surface geometry, we are
just dealing with a holomorphic map. It is precisely in this sense that we regard
complex function theory as a  special case of conformal surface theory.

To translate this conceptual idea into a mathematical one, we rewrite
the conformality condition for an immersion into $\R^3$ using quaternions,
similar to expressing the Cauchy-Riemann
equation in terms of complex numbers. 
Recall that an immersion $f\colon M\to\R^3$ is conformal if the induced metric
is in the conformal class of $M$. In other words, images under $df$
of orthogonal tangent vectors 
of equal length (with respect to any metric in the conformal class) remain orthogonal
and are again of equal length. If $N\colon M\to S^2\subset \R^3$ denotes the unit normal
map which makes $f$ positively oriented then it is easy to check that
we can write the conformality condition as
%\begin{equation}\label{eq:conf_imm}
\[
*df=N\times df\,.
\]
%\end{equation}
If we regard $\R^3$ as the imaginary quaternions, $\text{Im}\,\H\,$, 
the vector cross product  between  perpendicular
vectors is given by quaternionic multiplication. Thus, our conformality equation for
an immersion $f$ with unit normal $N$ becomes the Cauchy-Riemann type equation
%\begin{equation}\label{eq:conf_imm_quat}
\[
*df=N\,df\,
\]
%\end{equation}
with ``varying $i$''.
If $f$ takes values in a 2-plane its unit normal is a constant map,
say $N=i$, in which case we recover the usual Cauchy-Riemann equation
%\eqref{eq:Cauchy_Riemann} 
for $f$. With this in mind we arrive at the following 
extended notion of holomorphicity \cite{icm}:
\begin{definition*}
%\label{def:quat_hol_function}
A map $f\colon M\to\H$ of a Riemann surface into $\H$ is {\em holomorphic}
if there exists a map $N\colon M\to\H$ with $N^2=-1$ so that
\[
*df=N\,df\,.
\]
\end{definition*}
The condition $N^2=-1$ simply means that $N$ is an imaginary quaternion of length one,
i.e., a map into the 2-sphere in $\R^3$.
Note that we do not require $f$ to be immersed, thus allowing for
isolated branch points whose structure is described -- in a more general 
setting -- in Section~\ref{sec:Pluecker}. At immersed points
our notion of holomorphicity is equivalent to conformality and  
the map $N$ is uniquely determined: the Gauss map of $f$ takes values in the
oriented $2$-plane Grassmannian $S^2\times S^2$ of $\R^4$, 
and $N$ is the projection into the
first $S^2$. In particular, if $f$ is $\R^3$ valued then $N$ is the 
unit normal map. A detailed development of conformal surface theory using
quaternionic valued functions, with applications to the theory of Willmore surfaces, 
can be found in \cite{bflpp}. For now it is sufficient to keep in mind
that a holomorphic function into
$\H$ is a conformal map into $\R^4$ with isolated branch points.

Even though our notion of holomorphicity for a function $f$ is formulated in euclidean 
terms, it is
invariant under M\"obius transformations of $\R^4$. This yields the concept of 
a holomorphic map into $\H\P^1=S^4$ and, more generally, into quaternionic Grassmannians,
found in Section~\ref{subsec:holo_curve}. For the M\"obius invariant formulation and
in order to connect to the holomorphic bundle theory, we have to 
interpret the holomorphic function $f$ as the ratio of holomorphic sections of a suitable
line bundle.  Consider the trivial quaternionic line bundle
\[
L=M\times \H
\]
and denote by $\nabla$ the trivial quaternionic connection given by differentiation of 
$\H$ valued functions on $M$. If $\phi$ is a fixed constant section then any
other section $\psi$ is of the form $\psi=\phi\lambda$ for some function $\lambda\colon M\to\H$.
The bundle $L$ carries the quaternionic linear endomorphism
\[
J\phi:=\phi N
\]
defined in the basis $\phi$ by the map $N$. Since $N^2=-1$ the endomorphism $J$ is
a complex structure making $L$ into a rank 2 complex bundle. As such 
\[
L=E\oplus E
\]
is the double of a complex line bundle $E$ whose degree, 
in case $M$ is compact with genus $g$,
is the mapping degree  of $N\colon M\to S^2$. Here $E$ is the $i$-eigenspace bundle of
$J$ consisting of vectors $\psi$ in $L$ for which $J\psi=\psi\,i$.
% We let the (quaternionic) degree of 
%the line bundle $L$ to be the degree of $E$, i.e., half the degree of the rank 2
%complex bundle $L$. 
Using the complex structure $J$ on $L$ the 
trivial connection $\nabla$ decomposes into types
\[
\nabla=\nabla'+\nabla''\,.
\]
The resulting first order linear operator
\[
\nabla''\colon \Gamma(L)\to \Gamma(\bar{K}L)
\]
satisfies the Leibniz rule
\[
\nabla''(\psi\lambda)=(\nabla''\psi)\lambda+(\psi\,d\lambda)''
\]
over quaternionic valued functions. Since the constant section $\phi$ is parallel,
it clearly satisfies $\nabla''\phi=0$. But there is a second section
in the kernel of $\nabla''$: due to the holomorphicity of $f$ also $\psi=\phi f$
satisfies $\nabla''\psi=0$,  as one sees from the Leibniz rule. In other words,
given a holomorphic map $f\colon M\to \H$ we have attached to it a quaternionic line bundle
with complex structure and a certain quaternionic linear first order operator. 
Its kernel contains a $2$-dimensional linear subspace of sections with basis
$\psi, \phi\in\Gamma(L)$ such that $f$ is the ratio 
\[
\psi=\phi\, f\,.
\]
Any other choice of basis is related by an element of $\text{\bf Gl}(2,\H)$ which 
acts by M\"obius transformations on $f$. 

This construction tells us what the notion of a holomorphic
structure on a quaternionic vector bundle with complex structure has to be: an operator
of the type $\nabla''$. Thus, reminiscent of the theory of complex $\delbar$-operators,
we describe in Section~\ref{subsec:quat_holo} a {\em holomorphic structure} 
by a first order quaternionic linear
operator
\[
D\colon \Gamma(L)\to \Gamma(\bar{K}L)
\]
satisfying the above quaternionic Leibniz rule. Sections are called {\em holomorphic}
if they are contained in the kernel of $D$, and the space of all
holomorphic sections is denoted by $H^0(L)$. Given a $2$-dimensional linear 
subspace of holomorphic sections of a quaternionic holomorphic line bundle, 
any choice of basis $\psi,\phi$ determines a map $f\colon M\to \H$ by
$\psi=\phi\, f$, provided that $\phi$ never vanishes. We use the complex structure $J$
on $L$ to define the map $N\colon M\to S^2\subset \text{Im}\,\H$ by
\[
J\phi=\phi N\,,
\]
and then the quaternionic Leibniz rule for $D$ implies the holomorphicity condition
$*df=N\,df$ for the function $f$. A different choice of basis would have resulted in 
a M\"obius transform of $f$. 

This correspondence between M\"obius equivalence classes of quaternionic valued
holomorphic functions, i.e., branched conformal immersions,
and $2$-dimensional linear subspaces of holomorphic sections of 
quaternionic holomorphic line bundles is a special instance of the 
{\em Kodaira embedding} described in Section~\ref{subsec:Kodaira}. In general, 
$n+1$-dimensional linear subspaces of holomorphic sections (without base points) of 
a quaternionic holomorphic line bundle correspond to projective 
equivalence classes of holomorphic maps into $n$-dimensional quaternionic projective
space. The components of such a map are maps into
$\H$ that are all holomorphic with respect to the same $N$. Thus, through
quaternionic linear combinations and projections, a
single holomorphic curve into $\H\P^n$ yields families of branched conformal immersions 
into $\R^4$. 

At this point we have a description of conformal surface theory in terms of
quaternionic holomorphic line bundles, in precisely the
same way complex holomorphic curves are described by complex holomorphic line
bundles. We now show how the complex theory fits into our extended setup.
In general, a holomorphic structure $D$ does not commute with the
complex structure $J$ on the bundle $L$. The $J$-commuting part of $D$ 
is a complex $\delbar$-operator on $L=E\oplus E$, in fact the double of
a $\delbar$-operator on $E$, and therefore describes a complex holomorphic
structure on $E$. The $J$-anticommuting part $Q$ of $D$ can be viewed as a
complex antilinear endomorphism of $E$, the {\em Hopf field} of $D$.
The resulting decomposition of a quaternionic holomorphic structure  
\[
D=\delbar+Q
\]
allows us to define in Section~\ref{subsec:Willmore_energy}
the {\em  Willmore energy} of the holomorphic structure $D$, 
\[
W=2\int <Q\wedge *Q>\,,
\]
the $L^2$-norm of the Hopf field. The  holomorphic
structure $D$ on $L=E\oplus E$ is a complex holomorphic structure on $E$ if and only
if the Willmore energy of $D$ vanishes. Put differently, the
moduli space of quaternionic holomorphic line bundles is a bundle over the Picard
group of $M$ whose fibers are sections of antilinear endomorphism of the base points $E$. 
This (infinite dimensional) bundle carries the natural energy function $W$ whose zero level
set is the Picard group of $M$.

We now interpret the above
quantities in terms of the geometry of a conformal immersion $f\colon M\to\R^3$ with
unit normal map $N\colon M\to S^2$. We have seen that the relevant quaternionic line bundle
is $L=M\times \H$ with some fixed constant trivialization $\phi$. 
The complex structure
$J$ on $L$ is defined via the unit normal by $J\phi=\phi N$. As we
have seen, in case $M$ is compact with genus $g$, the degree of $L$ is
the mapping degree $1-g$ of $N$. The holomorphic
structure $D=\nabla''$ is the $(0,1)$-part of the trivial connection $\nabla$ on $L$,
and we know that both $\phi$ and $\psi=\phi f$ are holomorphic sections of $L$. 
If we decompose
\[
\nabla=\nabla_{+}+\nabla_{-}
\]
into $J$-commuting and anticommuting parts then $\nabla_{+}$ is a complex connection
on $L=M\times \H$, i.e., $[\nabla_{+},J]=0$. As such, $\nabla_{+}$ can be viewed as a
complex connection on the underlying complex line bundle $E$ of $L=E\oplus E$.
As one expects, $\nabla_{+}$ is closely related to the Levi-Civita connection on $M$ with 
respect to the induced metric $|df|^2$. In fact, $E^2$ can be identified with the
tangent bundle and $\nabla_{+}$ corresponds to the Levi-Civita
connection under this isomorphism. On the other hand, the $J$-anticommuting part 
$\nabla_{-}$ is an endomorphism of $L$ which, in the trivialization $\phi$, is given
by the shape operator $\tfrac{1}{2}N\,dN$ of $f$. So we see that
the above decomposition of the
trivial connection $\nabla$ is a somewhat more refined version of the usual decomposition
of derivatives along a surface into Levi-Civita connection and second fundamental form.
The holomorphic structure 
\[
D=\nabla''=\nabla_{+}''+\nabla_{-}''=\delbar+Q
\]
can now be easily calculated: $\delbar$ is a square root of the holomorphic 
structure of $M$, i.e., its dual is a holomorphic spin structure on $M$. The Hopf field $Q$
is given in the trivialization $\phi$ by 
the trace free shape operator 
\[
Q=\tfrac{1}{2}N(dN-Hdf)\,.
\]
The Willmore energy of the quaternionic holomorphic line bundle $L$ thus becomes
\[
W=\int(H^2-K)|df|^{2}\,,
\]
$K$ denoting the Gauss curvature, 
which is the classical M\"obius invariant Willmore energy of the immersion $f$.
Notice that the holomorphic structure $D=\delbar+Q$ is indeed invariant under
M\"obius transformations of $f$: the complex holomorphic structure $\delbar$,
the dual of the induced spin structure on $M$, 
does not change under M\"obius transformations and neither does the 
trace free part of the shape operator. Also note that the Willmore energy
vanishes if and only if the surface $f$ is totally umbilic, in which case
$f$ takes values in $S^2$ and can be viewed as a complex holomorphic map.

The study of the critical points and values of the Willmore energy 
on conformal immersions of a given compact Riemann surface is an important,
yet unresolved, problem of classical surface theory. From our discussion
we see that this problem translates to the study of the critical
points and values of the Willmore energy on the moduli space of quaternionic
holomorphic line bundles with the constraint that the bundles admit at least two 
linearly independent holomorphic sections. As we have seen above, the condition
on the number of holomorphic sections insures that
such line bundles indeed give rise to a conformal map $f\colon M\to\H\P^1$ via the Kodaira
embedding. When $M=S^2$ this program has recently been carried out
for an arbitrary number of independent holomorphic sections \cite{Katrin}.
The critical points, called {\em Willmore bundles}, can all be constructed 
from rational data generalizing earlier results of Bryant \cite{Bryant} and others 
\cite{Ejiri, Montiel}. Related to these issues are the
following more general questions:
\begin{enumerate}
\item
What can one say about the number $h^{0}(L)$ of holomorphic sections of a quaternionic
holomorphic line bundle $L$? 
\item
Given a linear subspace $H\subset H^{0}(L)$ of holomorphic sections of a quaternionic
holomorphic line bundle $L$, is there a lower bound on the Willmore energy $W$
in terms of fundamental invariants such as the dimension of $H$, degree
of the bundle $L$, genus of $M$, etc.?
\end{enumerate}
The general answer to the first question is given by the Riemann-Roch Theorem 
in Section~\ref{subsec:Riemann-Roch} which 
holds in verbatim form for quaternionic holomorphic line bundles:
\[
h^{0}(L)-h^{0}(KL^{-1})=\deg L -g +1\,.
\]
Here the degree of $L=E\oplus E$ 
is the degree of the underlying complex bundle $E$, and all dimension counts are
quaternionic. In particular, we get the existence of holomorphic sections with a lower
bound on $h^{0}(L)$ when the degree is large. 

Whereas complex holomorphic bundles
of negative degree do not admit holomorphic sections, quaternionic holomorphic
bundles generally do. We have seen above that the bundle induced from a conformal
immersion into $\R^3$ has degree $1-g$ and has at least two independent 
holomorphic sections. Unless $g=0$, this would be impossible for complex holomorphic bundles.
Thus the Willmore energy cannot be zero in such cases. This leads in 
Section~\ref{sec:Pluecker} to 
a central result of this paper, the {\em quaternionic Pl\"ucker formula},
providing an answer to the second question above.
Let $H\subset H^{0}(L)$ be an
$n+1$-dimensional linear subspace of holomorphic sections of a quaternionic
holomorphic line bundle $L$ of degree $d$. Then the Willmore energy satisfies
\[
\tfrac{1}{4\pi}(W-W^{*})=(n+1)(n(1-g)-d)+\ord H\,,
\]
which generalizes the classical Pl\"ucker formula of a com\-plex holo\-morphic curve
\cite{bible}. Rough\-ly speaking, $\ord H$
counts the singularities of all the higher osculating curves 
of the Kodaira embedding of $L$, and $W^{*}$ is the Willmore energy of 
the dual curve, i.e., the highest osculating curve. For example,
a holomorphic curve $f$ in $\H\P^n$ has $\ord H=0$ if and only if $f$ is a Frenet curve 
which, for $n=1$, simply means that $f$ is immersed. Note that in the complex case,
where $W=0$, the only Frenet curve is the rational normal curve. This demonstrates again
the rather special flavor of the complex theory from our viewpoint.
Since $W^{*}$ is nonnegative we obtain the lower bound
\[
\tfrac{W}{4\pi}\geq (n+1)(n(1-g)-d)+\ord H
\]
for the Willmore energy. Equality holds if the
bundle $L$, or rather its Kodaira embedding, 
is the dual curve of the twistor projection into $\H\P^n$ 
of a complex holomorphic curve in $\C\P^{2n+1}$. In other words, 
a quaternionic holomorphic line bundle $L$ for which equality holds
can be calculated from a complex holomorphic line bundle by
taking derivatives and performing algebraic operations.
In the case of genus zero, the above estimate is sharp and 
geometric examples include Willmore spheres \cite{Bryant, Ejiri, Montiel},
Willmore bundles \cite{Katrin} 
and, more generally, soliton spheres \cite{Taimanov, Paul}. 

To obtain a useful lower bound in higher genus $g\geq 1$, we use $\ord H$ to compensate
for the negative quadratic term. This results in  general quadratic estimates
in terms of the number $n+1$ of holomorphic sections for any $g\geq 1$,
\[
W\geq
\begin{cases}
\tfrac{\pi}{g}((n+g-d)^{2}-g^{2}) &\text{if}\qquad n\geq 0\,,\, n\geq d\,,\, d\leq g-1\\
\tfrac{\pi}{g}((n+1)^2-g^2) &\text{if}\qquad n\geq d-g+1\,,\, n\geq g-1\,,\, d\geq g-1\,,
\end{cases}
\]
which are discussed in 
Section~\ref{subsec:Willmore_estimates}. To summarize,
if a quaternionic holomorphic line bundle has more sections
than allowed in the complex holomorphic theory, its Willmore energy 
has to grow quadratically in the number of sections. 

Classically, the Pl\"ucker formula is proven by counting the singularities,
i.e., the zeros of derivatives, of the various osculating curves to a holomorphic curve into
$\C\P^n$. In the quaternionic setting the osculating curves fail to extend smoothly
across the singularities, requiring a more intrinsic viewpoint. For this purpose we develop
in Section~\ref{sec:jets} the theory of holomorphic jets for holomorphic structures
$\delbar+Q$ on complex vector bundles where $Q$ is complex antilinear. The main 
advantage of our axiomatic development of jet theory is its global and geometric flavor
needed in our applications to holomorphic curve theory. 
The section is self contained, including
a discussion on the basic local analytical properties of the first order elliptic 
operators $\delbar +Q$, which historically appeared as Carlmen-Bers-Vekua operators
in the literature \cite{Vekua}. 
\begin{comment}
which makes it immediately applicable to the case at hand. 
To keep the exposition self contained, we 
Even though there is a vast literature on the local analysis of first order elliptic 
operators of the form $\delbar +Q$,
so called Bers-Vekua operators, we need very little . To keep the exposition self contained, 
we boroughed the fundamental identity result for these operators from the literature 
and used it to derive the few analytical results we needed,
i.e., behavior of zeros etc..
\end{comment}

The paper concludes with two rather different applications of the quaternionic 
Pl\"ucker formula, demonstrating the extended scope of 
quaternionic holomorphic geometry. The
first application in Section~\ref{sec:Dirac} concerns the spectrum of the Dirac
operator on surfaces.  We give quantitative lower bounds for the 
eigenvalues in terms of their multiplicities $m$, 
\begin{equation*}
\lambda^2\text{\rm area}_M\geq
\begin{cases}
4\pi m^2 &  \text{if}\qquad g=0\\
\tfrac{\pi}{g}(m^2-g^2)&\text{if}\qquad g\geq 1\,,
\end{cases}
\end{equation*} 
where the Dirac operator $\mathcal{D}$ is defined on a Riemannian spin bundle 
over a compact surface of genus $g$ with Riemannian metric. Taking a constant
curvature metric on $S^2$, we see that the bounds are sharp in genus zero. 
Note that a Riemannian spin bundle $L$ has a canonical quaternionic holomorphic structure
$D=\delbar+Q$:  up to Clifford multiplication, $\delbar$ is the Dirac operator $\mathcal{D}$
and the Hopf field $Q$ is the identity map. The eigenvalue equation
$\mathcal{D}\psi=\lambda\psi$ therefore becomes the holomorphicity condition 
$(\delbar+\lambda Q)\psi=0$, the multiplicity of $\lambda$ is given by $h^{0}(L)$, and
the Willmore energy is $W=\lambda^2\text{area}_{M}$.
Applying the above estimates for $W$ to this setting, and keeping in mind
that spin bundles have degree $g-1$, we obtain our lower bounds on the Dirac
eigenvalues.

In genus zero partial results have already been known: the case of multiplicity one
has been proven by other methods in \cite {Bar}. Methods from soliton
theory, for metrics admitting continuous symmetries, were used for general
multiplicity in \cite{Taimanov}. 

Our estimates for surfaces of higher genus are new.

The second, and final, application discussed in Section~\ref{sec:harmonic} 
deals with energy estimates
of harmonic 2-tori in $S^2$. It is well known that non-conformal harmonic maps in $S^2$ correspond
to constant mean curvature, CMC, surfaces in
$3$-dimensional space forms. 
Since the area of the resulting surface is, up
to a constant depending on the size of the mean curvature, 
equal to the energy of the harmonic map, we are simultaneously dealing
with area estimates of CMC tori in $3$-dimensional space forms. 
This relationship is most explicit
in the case of CMC surfaces in $\R^3$: the harmonic map into
$S^2$ is the Gauss unit normal map of the surface and its energy is twice the
area of the CMC surface if the mean curvature $H=1$. 
A fundamental property of harmonic $2$-tori
is given by the fact that they are solutions
to an algebraically completely integrable hierarchy.
This means that every
harmonic 2-torus and CMC torus is parametrized by theta functions on a 
hyperelliptic Riemann surface of genus $g$, the {\em spectral genus}.
In genus zero the harmonic map is a geodesic on $S^2$ and the associated
CMC surface in $\R^3$ is the cylinder. Spectral genus one gives
the normals to the rotationally symmetric Delauney surfaces described by elliptic functions.
The Wente torus is an example for 
spectral genus two, and numerical examples have been computed and visualized up to 
spectral genus five \cite{Mathias}. These
experiments suggest that the area, the first integral of motion in
the integrable hierarchy, grows with the spectral genus. Applying
our lower bounds for the Willmore energy of quaternionic holomorphic line bundles
we obtain lower bounds, quadratic in the spectral genus, for the energy of 
harmonic $2$-tori and the area of CMC tori.  For instance, if a CMC torus $f\colon M\to \R^3$  
has spectral genus $g$ then its area satisfies 
\[
\text{\rm area}(f)\geq
\begin{cases}
\tfrac{\pi}{4}(g+2)^2 & \,\,\text{if $g$ is even, and}\\
\tfrac{\pi}{4}((g+2)^2-1)&\,\,\text{if $g$ is odd}\,. 
\end{cases}
\]
These estimates arise from quaternionic holomorphic bundle theory by 
considering the line bundle $L=M\times \H$ induced by a surface $f\colon M\to \R^3$ 
that was previously described in some detail. Its complex structure is given in 
terms of the unit normal map $N\colon M\to S^2$, and the holomorphic structure $D$ is
the $(0,1)$-part $\nabla''$ of the trivial connection $\nabla$. Furthermore, the
Willmore energy becomes
\[
W=\int(H^2-K)|df|^2=\int H^2|df|^2=\text{\rm area}(f)
\]
when $f\colon M\to\R^3$ is a $2$-torus of constant mean curvature $H=1$. 
The constant mean curvature condition on $f$, or equivalently the harmonicity of $N$,
is expressed by the flatness of the family of connections
\[
\nabla_{\lambda}=\nabla+(\lambda-1)\nabla_{-}'
\]
on $L$. Since $\nabla_{\lambda}''=D$ for all $\lambda$, parallel sections of 
$\nabla_{\lambda}$ are holomorphic. For a parallel section to exist over the torus $M$
it has to be a common fixed point of the holonomy $H_{\lambda}$ 
of the flat connection $\nabla_{\lambda}$.
This is where the spectral curve comes into play: the values $\lambda\in\C\P^1$ for 
which the holonomy $H_{\lambda}$ has coinciding eigenvalues, together with
$\lambda=0$ and $\lambda=\infty$, define a hyperelliptic compact Riemann surface 
of genus $g$ branched at those values. Since $\nabla_{\lambda}$ is a family
of $\text{\bf Sl}(2,\C)$-connections, unitary for $\lambda$ on the unit circle, there
are genus $g$ many pairs of branch values (not including $\lambda=0,\infty$)
that are symmetric with respect to the unit circle. Each such pair gives rise to
a parallel quaternionic section of $L$ with $\Z_2$ holonomy. Therefore we obtain $g$ 
independent parallel, 
and thus holomorphic, sections of $L$ -- possibly on a $4$-fold covering of $M$.
At $\lambda=1$ the connection $\nabla_{\lambda}$ is trivial to start with, 
adding one more holomorphic section to give $h^{0}(L)=g+1$. Up to here we only
used the periodicity of the normal map $N$. The closing condition of the CMC torus
$f\colon M\to\R^3$ yields an additional holomorphic section, finally resulting in the above area
estimate.   
 
We close with a number of geometric applications of these estimates. 
Reversing the viewpoint,
we obtain a quantitative bound on the genus of the spectral curve for CMC and
harmonic tori. In particular, we get a different proof that
the spectral curve has finite genus. 

For minimal tori in $S^3$ our estimates imply that the area, or Willmore energy,
is greater than $9\pi>2\pi^2$ for spectral genus four and higher. Therefore, 
verifying the Willmore conjecture \cite{Willmore} in this case ``only'' requires an analysis
of the situation in genus two and three. A similar remark applies to the
Lawson conjecture which states that the only embedded minimal
torus in $S^3$ is the Clifford torus. Since it is known \cite{Chinese}
that an embedded minimal torus has area less then $16\pi$ it ``suffices" to
check minimal tori resulting from spectral genera at most five.

A non-constant harmonic $2$-torus in $S^2$ can have energy arbitrarily close to zero,
which can be seen by mapping a very thin torus to a geodesic on $S^2$. 
Nevertheless, one expects the harmonic map to take values in an
equator if its energy is ``small enough''. Indeed, if the energy of a harmonic torus is
below $4\pi$ then  our estimates imply that it must be an equator. This value 
misses the conjectured sharp value of $2\pi^2$ by the
scalar factor of $\tfrac{\pi}{2}$. The value $2\pi^2$ 
is attained by the family of normals to the Delauney surfaces, whose limit is
the cylinder.

\begin{comment}
the endomorphism bundle
$\text{End}(L)$ splits into $J$-commuting and anticommuting 
endomorphisms $\text{End}_{\pm}(L)$, both of which are invariant under $\nabla_{+}$.
Via the trivialization $\phi$ this decomposition corresponds to the
splitting 
\[
\H=df(TM)\oplus (\R N\oplus \R)
\]
of $\H$ into tangent and normal spaces to $f$ and $\nabla_{+}$ respects this 
decomposition. This implies that $\nabla_{+}$ on $\text{End}_{-}(L)$
is the Levi-Civita connection of $f$. But $\text{End}_{-}(L)$ 
 
with $\H$ via the trivialization $\phi$. Then $\nabla_{+}$, as
a connection on endomorphisms, leaves the decomposition of $\H$ into 
$N$-commuting and anticommuting
subspaces invariant. The latter is the image of $TM$ under $df$ and the former is
the trivial complex bundle $\R\oplus\R N$. Thus the tangent part of $\nabla_{+}$ 
on endomorphisms of $L$ is
the Levi-Civita connection on $M$ with  
\end{comment}

%%% Local Variables: 
%%% mode: latex
%%% TeX-master: "willmore"
%%% End: 

\section{Quaternionic holomorphic geometry}\label{sec:quat_holo_geom}
\subsection{Preliminaries}\label{subsec:prelim}
We set up some basic notation used throughout the paper.
A Riemann surface $M$ is a $2$-dimensional, real manifold with an
endomorphism field $J_{M}\in\Gamma(\text{End}(TM))$ satisfying $J_{M}^2=-1$.
If $V$ is a vector bundle over $M$, we denote the space of $V$ valued $k$-forms 
by $\Omega^{k}(V)$.
If $\omega\in\Omega^{1}(V)$, we set
\begin{equation}\label{eq:hodge*}
*\omega:=\omega\circ J_{M}\,.
\end{equation}
For example, sections of the canonical bundle $K$ of $M$ are those 
$\omega\in\Omega^{1}(\C)$ for which $*\omega=i\omega$. The form
$\omega$ is holomorphic if, in addition, $d\omega=0$.

In calculations it is often useful to identify $2$-forms on a Riemann surface $M$
with quadratic forms. If $\omega\in\Omega^{2}(V)$ then the associated quadratic form,
which we will again denote by $\omega$, is given by
\begin{equation}\label{eq:quadform}
\omega(X):=\omega(X,J_{M}X)\,,
\end{equation}
where $X\in TM$. We will frequently use this identification in the
following setup: let $V_k$, $k=1,2,3$,
be vector bundles over $M$ which have a pairing $V_1\times V_2\to V_3$. 
If $\omega\in\Omega^1(V_1)$ and $\eta\in\Omega^1(V_2)$ then the
$V_3$ valued $2$-form $\omega\wedge\eta$, where the wedge is over the given pairing,
corresponds via \eqref{eq:quadform} to
the $V_3$ valued quadratic form
\begin{equation}\label{eq:identify}
\omega\wedge\eta=\omega(*\eta)-(*\omega)\eta\,.
\end{equation}
In particular
\[
\omega\wedge *\eta=-*\omega\wedge\eta\,.
\]

Most of the 
vector bundles occurring will be {\em quaternionic} vector bundles,
i.e., the fibers are quaternionic vector spaces and the local
trivializations are quaternionic linear on each fiber.
Using a transversality
argument, one can show that every quaternionic vector bundle over a $2$ or 
$3$-dimensional manifold is trivializable. 

A {\em quaternionic connection} on
a quaternionic vector bundle satisfies the usual Leibniz rule over
quaternionic valued functions.

We adopt the convention 
that all quaternionic vector spaces are {\em right} vector spaces. Using
quaternionic conjugation, any occurring left vector space can be made
into a right vector space via $\lambda\psi=\psi\bar{\lambda}$. 
For example, the dual $V^{-1}$ of all quaternionic linear forms on $V$ is naturally a left
vector space, but we will always consider it as a right vector space.

If $V_1$ and $V_2$ are quaternionic vector bundles, we denote by  $\text{Hom}(V_1,V_2)$
the bundle of quaternionic linear homomorphisms. As usual, 
$\text{End}(V)=\text{Hom}(V,V)$ denotes the quaternionic linear endomorphisms.
Notice that $\text{Hom}(V_1,V_2)$ is {\em not} a quaternionic bundle.

Any $J\in\Gamma(\text{End}(V))$ with $J^2=-1$ gives a complex structure
on the quaternionic bundle $V$. Since $J$ is quaternionic linear, this 
complex structure is compatible with the quaternionic structure.
The $\pm i$-eigenbundles $V_{\pm}=\{\psi\in V\,;\, J\psi=\pm\psi\,i\}$
are complex vector bundles, via $J$, whose complex rank equals the
quaternionic rank of $V$. Multiplication by, say $j$, is a complex linear
isomorphism between $V_{+}$ and $V_{-}$. Thus, up to complex isomorphisms,
$V=W\oplus W$ is a double for a unique
complex vector bundle $W$. Conversely, every 
complex vector bundle $W$ can be made into a quaternionic vector bundle
$V=W\oplus W$ with a natural complex structure by doubling: the quaternionic
structure is then defined by $(\psi,\phi)i:=(i\psi,-i\phi)$ and $(\psi,\phi)j:=
(-\phi,\psi)$ and $ij=k$. The above two operations are inverse to each other on the
respective isomorphism classes of bundles.

Any complex structure $J$ on $V$ can be dualized to give a complex
structure $J$ on $V^{-1}$ via $<J\alpha,\psi>=<\alpha,J\psi>$, where
$\alpha\in V^{-1}$ and $\psi\in V$. If $V=W\oplus W$ then
$V^{-1}=W^{-1}\oplus W^{-1}$. 

If $V$ is a quaternionic bundle with complex structure $J$ and $E$ is a complex bundle,
we can tensor over $\C$ to obtain the quaternionic bundle $EV$ with
complex structure $J(e\psi):=eJ(\psi)$. If $W$ is the complex bundle underlying $V$, then
$EW$ is the complex bundle underlying $EV$. Over a Riemann surface typical
examples of such bundles are 
\begin{equation}\label{eq:KV}
KV=\{\omega\in\Omega^1(V)\,;\,*\omega=J\,\omega\}\;\;\text{and}\;\;
\bar{K}V=\{\omega\in\Omega^1(V)\,;\,*\omega=-J\,\omega\}\,.
\end{equation}
The quaternionic linear splitting 
\begin{equation}\label{eq:type_splitting}
T^{*}M\otimes V= KV\oplus\bar{K}V
\end{equation}
induces the {\em type decomposition} of $V$ valued $1$-forms $\omega\in\Omega^1(V)$ into
\begin{equation}\label{eq:types}
\omega=\omega'+\omega''\,,
\end{equation}
with $K$-part $\omega'=\frac{1}{2}(\omega-J*\omega)\in\Gamma(KV)$ and 
$\bar{K}$-part $\omega''=\frac{1}{2}(\omega+J*\omega)\in\Gamma(\bar{K}V)$.

Given two quaternionic bundles $V_k$ with complex structures, the real 
bundle $\text{Hom}(V_1,V_2)$ has two natural complex structures:
a left complex structure given by
composition with $J_2$ and a right complex structure given by
pre-composition with $J_1$. Unless noted, we will use the left complex
structure on Hom-bundles. $\text{Hom}(V_1,V_2)$ splits into the 
direct sum of the two 
complex subbundles $\text{Hom}_{\pm}(V_1,V_2)$ of complex linear, respectively antilinear,
homomorphisms. With respect to the left complex structure, 
\begin{equation}\label{eq:homs}
\text{Hom}_{+}(V_1,V_2)=\text{Hom}_{\C}(W_1,W_2)\;\;\text{and}\;\;
\text{Hom}_{-}(V_1,V_2)=\text{Hom}_{\C}(\overline{W_1},W_2)\,,
\end{equation}
where $W_k$ are the underlying complex bundles to $V_k$.

If we have a  $1$-form $\omega$ with values in $\text{Hom}(V_1,V_2)$,
it can be of type $(1,0)$ with respect to the left or
right complex structure. We use the notation
\begin{equation}\label{eq:lrK}
K\text{Hom}(V_1,V_2)=\{\omega\,;\,*\omega=J_2\,\omega\}\;\;
\text{and}\;\;
\text{Hom}(V_1,V_2)K=\{\omega\,;\,*\omega=\omega\,J_1\}
\end{equation}
for the $(1,0)$-forms with respect to the left and right
complex structures and call such $1$-forms {\em left $K$}, respectively
{\em right $K$}. Similar notations apply to $(0,1)$-forms.
For example, $K\text{Hom}_{-}(V_1,V_2)=\text{Hom}_{-}(V_1,V_2)\bar{K}$.

On a Riemann surface $K\wedge K=0$, so that for 
$1$-forms $\omega\in\Gamma(\text{Hom}(V_2,V_3)K)$ and $\eta\in\Gamma(K\text{Hom}(V_1,V_2))$
we have
\[
\omega\wedge\eta=0\,,
\]
where the wedge is over composition.
A typical example of such type considerations is the following: let
$\omega\in\Gamma(K\text{Hom}_{-}(V_2,V_3))$ and $\eta\in\Gamma(\bar{K}\text{Hom}(V_1,V_2))$,
then $\omega\wedge\eta=0$ since $\omega$ is also right $\bar{K}$.
 
If $V$ is a quaternionic bundle with complex structure over a compact,
$2$-dimensional, oriented manifold $M$ we define its degree to be the degree
of the underlying complex bundle $W$, 
\begin{equation}\label{eq:degree}
\deg V:=\deg W\,.
\end{equation}
For any complex bundle $E$ over $M$ we then have 
\begin{equation}\label{eq:degmult}
\deg EV=\rank V\deg E +\rank E\deg V\,.
\end{equation}
As in the complex case degrees can also be calculated by curvature integrals. 
If $B\in\text{End}(V)$ we denote by
\begin{equation}\label{eq:trace}
<B>:=\tfrac{1}{4}\trace_{\R}(B)
\end{equation}
the trace of $B$ viewed as a real endomorphism. Note that the trace of a quaternionic
endomorphism is not defined, but the real part of the trace is. Our normalization is such that 
$<B>=\Re\trace(B)$ and thus $<Id>=\rank V$. 
Now assume that $\nabla$ is a connection on $V$ so that $\nabla J=0$.
Then we have
\begin{equation}\label{eq:chern}
2\pi\deg V=\int<JR^{\nabla}>
\end{equation}
where $R^{\nabla}$ is the curvature $2$-form of $\nabla$.

Unless we need to emphasize it, the adjective {\em quaternionic} will be dropped.

\subsection{Quaternionic holomorphic structures}\label{subsec:quat_holo}
Let $V$ be a quaternionic vector bundle with complex structure $J$
over a Riemann surface $M$.
\begin{definition}\label{def:quat_holo} 
A {\em quaternionic holomorphic structure} is given by a quaternionic linear map
\begin{equation}
D\colon \Gamma(V)\to\Gamma(\bar{K}V)
\end{equation}
satisfying the Leibniz rule
\begin{equation}\label{eq:Leibniz}
D(\psi\lambda)=(D\psi)\lambda+(\psi d\lambda)''\,,
\end{equation}
for sections $\psi\in\Gamma(V)$ and
quaterionic valued functions $\lambda\colon M\to\H$. Here 
$(\psi d\lambda)''=\tfrac{1}{2}(\psi d\lambda+J\psi *d\lambda)$
denotes the $\bar{K}$-part \eqref{eq:type_splitting} of the $L$ valued 1-form $\psi d\lambda$.
Note that $d\lambda''$ would not make
sense  since $\H$ has no distinguished complex structure.

A section $\psi\in\Gamma(V)$ is called {\em holomorphic} if $D\psi=0$.
We denote the quaternionic subspace of holomorphic sections by
\[
H^{0}(V):=\ker D\subset \Gamma(V)\,.
\]
\end{definition}
We will see below that $D$ is a zero-order perturbation of a complex
$\delbar$-operator which is elliptic. This implies that over a compact
Riemann surface $H^{0}(V)$ is finite dimensional.

A bundle homomorphism $T\colon V\to\tilde{V}$ between two quaternionic holomorphic bundles
is called {\em holomorphic} if $T$ is complex
linear and preserves the holomorphic structures, i.e., $\tilde{J}T=TJ$ and $\tilde{D}T=TD$.  
A subbundle $\tilde{V}\subset V$ is a {\em holomorphic subbundle} if the inclusion
map $\tilde{V}\subset V$ is holomorphic. If $\tilde{V}\subset V$ is
a holomorphic subbundle then the quotient bundle $V/\tilde{V}$, with the
induced structure, is quaternionic holomorphic and the 
quotient projection $\pi\colon V\to V/\tilde{V}$ is holomorphic.  

To describe the space of all quaternionic holomorphic structures we
decompose $D$ into $J$-commuting and anticommuting parts
\[
D=\tfrac{1}{2}(D-JDJ)+\tfrac{1}{2}(D+JDJ)\,.
\]
It is easy to see that 
\begin{equation}\label{eq:delbar}
\delbar:=\tfrac{1}{2}(D-JDJ)
\end{equation} 
is again a holomorphic structure. 
Since $\delbar$ is $J$ commuting it is a {\em complex} holomorphic structure on 
$V=W\oplus W$ and
induces a complex holomorphic structure $\delbar_{W}$ on the underlying complex bundle $W$
(see section~\ref{subsec:prelim}). Thus
\begin{equation}\label{eq:delbar_double}  
\delbar=\delbar_{W}\oplus\delbar_{W}
\end{equation}
is the double of the complex holomorphic structure $\delbar_{W}$.
In the sequel we simply will use $\delbar$ and think of it as either the double
of a holomorphic structure on $W$ or the holomorphic structure on $W$ itself, depending on the
context. The $J$-anticommuting part
\begin{equation}\label{eq:Q}
Q:=\tfrac{1}{2}(D+JDJ)
\end{equation}
is a section of $\bar{K}\text{End}_{-}(V)$, i.e.,
\[
*Q=-JQ=QJ\,.
\]
For reasons explained in \cite{icm}, we call $Q$ the {\em Hopf field} of
the holomorphic structure $D$.
We have thus decomposed every holomorphic structure $D$ on $V$ as
\begin{equation}\label{eq:hol_decomp}
D=\delbar+Q\,,
\end{equation}
with $\delbar$ a complex holomorphic structure (on the underlying complex vector bundle $W$),
and $Q$ a section in  $\bar{K}\text{End}_{-}(V)$.
The complex holomorphic theory is characterized by the absence of the Hopf field $Q$.

A holomorphic bundle map $T\colon V\to \tilde{V}$ is complex linear and 
has $\tilde{D}T=TD$. The latter is equivalent to 
$\tilde{Q}T=TQ$ and $\delbar T=T\delbar$. Thus, a quaternionic holomorphic bundle map
is a complex holomorphic bundle map intertwining the respective Hopf fields.

An {\em antiholomorphic} structure on a quaternionic vector bundle $V$ with complex structure
$J$ is a holomorphic structure on the bundle $V$ with the opposite complex structure $-J$.
Thus, every antiholomorphic structure is of the form $\del+A$ with $\del$ 
(the double of) a complex antiholomorphic structure (on the underlying 
complex bundle $W$) and $A$ a section of $K\text{End}_{-}(V)$.
 
There is an obvious relation between holomorphic structures and quaternionic
connections. If $\nabla$ is a quaternionic connection on the bundle $V$
with complex structure $J$, we have the type decomposition \eqref{eq:types}
\[
\nabla=\nabla'+\nabla''\,,
\]
where
\begin{equation}\label{eq:nabla_type}
\nabla'= \tfrac{1}{2}(\nabla-J*\nabla)\qquad\text{and}\qquad\nabla''=\tfrac{1}{2}(\nabla+J*\nabla)\,.
\end{equation}
Decomposing $\nabla'$ and $\nabla''$ further into  
$J$-commuting and anticommuting parts, we obtain
\begin{equation}\label{eq:nabla_decomposition}
\nabla=\nabla'+\nabla''=(\del+A)+(\delbar+Q)\,.
\end{equation}
In particular, $\nabla'$ and $\nabla''$ are antiholomorphic, respectively 
holomorphic, structures on $V$. Of course,
every quaternionic holomorphic structure $D$ on $V$ can be augmented 
to a quaternionic connection in various ways by adding 
antiholomorphic structures.

Note that $\del+\delbar$ is a complex connection, i.e., $(\del+\delbar)J=0$, so that
\[
\nabla J=(\del+\delbar)J+[(A+Q),J]=2(*Q-*A)\,.
\]
Here we used $*A=JA=-AJ$ and $*Q=-JQ=QJ$. Differentiating once more, we obtain
\[
d^{\nabla}(\nabla J)=[R^{\nabla}, J]= 2\,d^{\nabla}(*Q-*A)\,,
\]
where $d^{\nabla}$ is the exterior derivative on $\text{End}(V)$ valued forms.
We summarize the basic two formulas,
\begin{equation}\label{eq:dJ}
\nabla J=2(*Q-*A)\qquad \text{and}\qquad [R^{\nabla}, J]=2\,d^{\nabla}(*Q-*A)
\end{equation}
for a quaternionic connection $\nabla$ on a complex quaternionic vector bundle.

\subsection{Pairings and Riemann-Roch}\label{subsec:Riemann-Roch}
A standard construction of the complex holomorphic theory
is to use the product rule to relate holomorphic structures on suitably paired
bundles. We will now discuss their quaternionic counter parts.
\begin{definition}\label{def:mixed}
Let $V$ be a complex quaternionic vector bundle.
A {\em mixed} structure is a quaternionic linear map
\[
\hat{D}\colon \Gamma(V)\to \Omega^{1}(V)
\]
satisfying the usual Leibniz rule \eqref{eq:Leibniz} and 
the condition $*\hat{D}=-\hat{D}J$.
\end{definition}
Taking the $J$-commuting and anticommuting parts of a mixed structure, we learn
that
\[
\hat{D}_{+}=\tfrac{1}{2}(\hat{D}-J\hat{D}J)=\tfrac{1}{2}(\hat{D}+*J\hat{D})=\hat{D}''\,,
\]
and similarly,
\[
\hat{D}_{-}=\hat{D}'\,.
\]
Since $J$-anticommuting parts are always tensorial we obtain
\[
\hat{D}=\delbar+A
\]
with $\delbar$ a complex holomorphic structure, and
$A$ a section of $K\text{End}_{-}(V)$. This explains the term {\em mixed}. 
The motivation for those structures comes from the fact that under
the evaluation pairing of $V^{-1}$ with $V$, holomorphic structures and
mixed structures correspond.

\begin{lemma}\label{lem:mixed_holomorphic}
Let $V$ be a quaternionic vector bundle with complex structure $J$ and
holomorphic structure $D=\delbar+Q$.
Then, for $\alpha\in\Gamma(V^{-1})$ and $\psi\in\Gamma(V)$, the product rule
\begin{equation}\label{eq:mixed_holomorphic}
\tfrac{1}{2}(d<\alpha,\psi>+*d<\alpha,J\psi>)=<\hat{D}\alpha,\psi>+<\alpha,D\psi>
\end{equation}
defines the mixed structure $\hat{D}=\delbar-Q^{*}$ on $V^{-1}$.

Conversely, given a mixed structure $\hat{D}=\delbar+A$ on $V^{-1}$, the above
product rule defines the holomorphic structure $D=\delbar-A^{*}$ on $V$.

In the complex setting, $Q=0$, the relation \eqref{eq:mixed_holomorphic} defines 
the dual holomorphic
structure on the dual bundle $V^{-1}$.
\end{lemma}
\begin{proof}
Given $D$, or $\hat{D}$, one first shows that \eqref{eq:mixed_holomorphic} is tensorial in $\psi$,
or $\alpha$. The product rule then uniquely defines the map $\hat{D}$, or $D$. It is
easy to check their required properties. Finally, if $Q$ is a section of
$\bar{K}\text{End}_{-}(V)$ then $Q^{*}$ is a section of $K\text{End}_{-}(V)$ since
\[
<*Q^{*}\alpha,\psi>=<\alpha,*Q\psi>=<\alpha,QJ\psi>=<JQ^{*}\alpha,\psi>\,.
\]
\end{proof}

A fundamental construction of holomorphic bundle theory is the 
Riemann-Roch pairing 
\[
KV^{-1}\times V\to T^{*}M\otimes\H:(\omega,\psi)\mapsto <\omega,\psi>\,.
\]
Unlike in the case discussed above, holomorphic structures do correspond via
this pairing \cite{icm}. If $V$ is a quaternionic holomorphic bundle, then
there is a unique
quaternionic holomorphic structure $\tilde{D}$ on $KV^{-1}$ such that 
\begin{equation}\label{eq:Riemann_Roch_pairing}
d<\omega,\psi>=<\tilde{D}\omega,\psi>-<\omega\wedge D\psi>
\end{equation}
for $\omega\in\Gamma(KV^{-1})$ and $\psi\in\Gamma(V)$. The respective Hopf fields
are related by
\begin{equation}\label{eq:Hopf_field_conjugation}
<\tilde{Q}\omega,\psi>-<\omega\wedge Q\psi>=0\,.
\end{equation}
This is precisely the setup needed for the Riemann-Roch Theorem.
\begin{theorem}\label{thm:Riemann-Roch}
Let $V$ be a quaternionic holomorphic vector bundle with holomorphic
structure $D$ over a compact Riemann surface
of genus $g$. Let  $\tilde{D}$ be the Riemann-Roch paired holomorphic structure
on $KV^{-1}$. Then
\begin{equation}\label{eq:Riemann-Roch}
\dim_{\H}H^{0}(V)- \dim_{\H}H^{0}(KV^{-1})=\deg V-(g-1)\rank V\,.
\end{equation}
\end{theorem}
\begin{proof}
Integrating \eqref{eq:Riemann_Roch_pairing} over the compact Riemann surface shows
that $\tilde{D}$ is the adjoint elliptic operator to $D$. Using the invariance
of the index of an elliptic operator under continuous deformations and the 
Riemann-Roch Theorem in the complex holomorphic setting, we obtain
\begin{align*}
\dim_{\H}H^{0}(V) &- \dim_{\H}H^{0}(KV^{-1})\\
&=\dim_{\H}\ker D-\dim_{\H}\ker\tilde{D}
=\dim_{\H}\ker D-\dim_{\H}\coker D\\
&=\tfrac{1}{4}\text{index}\,D=\tfrac{1}{4}\text{index}\,\delbar\\
&=\tfrac{1}{2}\dim_{\C}H^{0}(W\oplus W)-
\tfrac{1}{2}\dim_{\C}H^{0}(KW^{-1}\oplus KW^{-1})\\
&=\dim_{\C}H^{0}(W)-\dim_{\C}H^{0}(KW^{-1})
=\deg_{\C}W -(g-1)\rank_{\C} W\\
&=\deg V -(g-1)\rank V\,.
\end{align*}
\end{proof}

\subsection{The Willmore energy of quaternionic holomorphic structures}
\label{subsec:Willmore_energy}  
Over a compact Riemann surface we can assign to a quaternionic holomorphic structure an
energy. 
\begin{definition}\label{def:Will_energy}
Let $V$ be a holomorphic vector bundle over a compact Riemann surface
with holomorphic structure $D=\delbar+Q$. The {\em Willmore energy}
of $V$ (or, more accurately $D$) is defined by
\[
W(V)=2\int<Q\wedge *Q>\,,
\]
where $<\,>$ is the trace form introduced in \eqref{eq:trace}.
\end{definition}
Due to \eqref{eq:Hopf_field_conjugation}
the Riemann-Roch paired bundle $KV^{-1}$ has the same Willmore energy
$W(KV^{-1})= W(V)$.
When $V$ is a line bundle the trace form is
definite on $\text{End}_{-}(V)$ so that the Willmore energy of
a holomorphic structure $D$ vanishes if and only if $D=\delbar$ is
a complex holomorphic structure.

The relationship to the usual Willmore energy on immersed surfaces 
is discussed extensively in \cite{bflpp}, and we will remark on this
differential geometric application throughout the paper.  

\subsection{Holomorphic maps into Grassmannians}\label{subsec:holo_curve}
A geometric situation where the above holomorphic bundle
theory applies is that of maps into quaternionic Grassmannians $G_{k}(\H^n)$.
Let $\Sigma\to G_{k}(\H^n)$ denote the 
tautological $k$-plane bundle whose fiber over  $V\in  G_{k}(\H^n)$ is
$\Sigma_{V}=V\subset \H^n$. A map $f\colon M\to G_k(\H^n)$ is the same as a rank $k$ subbundle 
$V\subset\underline{\H^n}$ of the trivial $\H^n$-bundle over $M$
via $V=f^{*}\Sigma$, i.e., $V_p=\Sigma_{f(p)}=f(p)\subset \H^n$ for $p\in M$.
Usually it is clear from the context that one thinks of a vector space as a trivial
bundle over $M$, so we can avoid the cumbersome notation $\underline{\H^{n}}$. 

It will be notationally and conceptually useful 
to adopt a slightly more general setup which
includes maps with holonomy: we replace the
trivial $\H^n$-bundle by a rank $n$ quaternionic vector bundle $H$ over $M$
with a flat connection $\nabla$. A rank $k$ subbundle $V\subset H$ gives rise
to a map $f\colon \tilde{M}\to G_k(\H^n)$ from the universal cover of $M$
which is equivariant with respect to the holonomy representation of $\nabla$ in
$\text{\bf Gl}(n,\H)$. Conversely, any such equivariant map 
defines a flat rank $n$ bundle with a $k$-plane subbundle. 
In case the connection $\nabla$ is trivial,
i.e., has no holonomy, $H$ can be identified with the trivial ${\H^n}$ bundle
 over $M$, and $\nabla$ is just the directional derivative
of $\H^n$ valued functions. From now on, we will make no distinction
between an (equivariant) map $f$ into the Grassmannian $G_k(\H^n)$
and the corresponding subbundle $V\subset H$. 

The derivative of $V\subset H$ is given by the $\text{Hom}(V,H/V)$ valued $1$-form
\begin{equation}\label{eq:derivative}
\delta=\pi\,\nabla\,,
\end{equation}
where $\pi\colon H\to H/V$ is the canonical projection. Under the usual identification
of $TG_k(\H^n)=\text{Hom}(\Sigma,{\H^n}/\Sigma)$ the $1$-form
$\delta$ is the derivative $df$ of $f\colon M\to G_k(\H^n)$. We will
say that $V\subset H$ is {\em immersed} if its derivative $\delta$
has maximal rank, i.e., $\delta_p(T_p M)\subset \text{Hom}(V,H/V)_p$ is 
a real $2$-plane for $p\in M$.
\begin{definition}\label{def:holo_curve}
Let $M$ be a Riemann surface, $H$ a flat quaternionic $n$-plane bundle
and $V\subset H$ a rank $k$ subbundle. We say that $V$ is a {\em holomorphic curve}
(thought of as an equivariant map into the Grassmannian)
if there exists a complex structure $J\in\Gamma(\text{End}(V))$, $J^2=-1$,
so that 
\begin{equation}\label{eq:right_holo}
*\delta=\delta\,J\,.
\end{equation}
Put differently, $V\subset H$ is a holomorphic curve if there is a
complex structure $J$ on $V$ such that $\delta$ is a section
of $\text{Hom}(V,H/V)K$.

In this setting the {\em degree} of the holomorphic curve is the degree
of the dual bundle $V^{-1}$.  
\end{definition}
It is important not to confuse the notion of {\em holomorphic} just given with
that of a holomorphic subbundle in Definition~\ref{def:quat_holo}: 
the bundle $H$ does not have a holomorphic structure so
it makes no sense to demand $V\subset H$ to be a holomorphic subbundle. 
Equation \eqref{eq:right_holo} is a Cauchy-Riemann-type relation on the (equivariant) map 
$f\colon \tilde{M}\to G_k(\H^n)$. However, if all occurring objects were complex (instead
of quaternionic) then  \eqref{eq:right_holo} just gives the usual condition
for complex holomorphic curves into Grassmannians: $J$ would be multiplication by 
$i$, and the subbundle $V$ would then also be a holomorphic subbundle. 
In this case everything ties together.

Note that our definition of holomorphicity of a curve $f\colon M\to G_k(\H^n)$ is 
projectively invariant, i.e., $f$ is holomorphic if and only if 
$Af$ is holomorphic (with respect to the complex structure $AJA^{-1}$ on AV) for any
$A\in\text{\bf Gl}(n,\H)$.

A useful equivalent characterization of holomorphic curves $V\subset H$ is in terms of
mixed structures (see Definition~\ref{def:mixed}): the condition 
$*\delta=\delta J$ is, via \eqref{eq:derivative}, equivalent to the fact that
\begin{equation}\label{eq:holo_curve_mixed}
\hat{D}:=\tfrac{1}{2}(\nabla+ *\nabla J)\colon \Gamma(V)\to\Omega^{1}(V)
\end{equation}
leaves $V$ invariant. Moreover, $\hat{D}$ is a mixed structure on $V$
and as such induces a holomorphic structure $D$ on $V^{-1}$ by the
product rule \eqref{eq:mixed_holomorphic}: for $\alpha\in\Gamma(V^{-1})$ and
$\psi\in\Gamma(V)$ the holomorphic structure is given by
\begin{equation}\label{eq:holo_structure_induced}
<{D}\alpha,\psi>=\tfrac{1}{2}(d<\alpha,\psi>+*d<\alpha,J\psi>)-<\alpha,\hat{D}\psi>\,.
\end{equation}
We therefore obtain
\begin{theorem}\label{thm:curve_holo_structure}
Let $V\subset H$ be a holomorphic curve. Then there exists a unique quaternionic
holomorphic structure $D$ on $V^{-1}$ with the following property:
if $\alpha$ is a local parallel section of $H^{-1}$ then its restriction
to $V\subset H$ gives a local holomorphic section of $V^{-1}$.
\end{theorem}
\begin{definition}\label{eq:Will_energy_curve}
The {\em Willmore energy of a holomorphic curve} $V\subset H$ is given
by the Willmore energy $W(V^{-1})$ of the holomorphic bundle $V^{-1}$.  
\end{definition}
\begin{proof}
Equation \eqref{eq:holo_structure_induced} shows that local parallel sections 
$\alpha$ of $H^{-1}$ induce
local holomorphic sections $\alpha_{|V}$ of $V^{-1}$. For the uniqueness part 
take any other holomorphic
structure $D'$ on $V^{-1}$ which thus differs from $D$ by a section 
$B\in\Gamma(\bar{K}\text{End}(V^{-1}))$.
Since $D$ and $D'$ both annihilate parallel local sections of $H^{-1}$, also $B$ annihilates
such sections. But at each point $p\in M$ we may always choose a basis of $V^{-1}_p$
consisting of elements in $H^{-1}_p$, which then can be locally extended to
parallel sections and thus $B_p=0$.
\end{proof}
\begin{example}\label{ex:standard_1}
At this point it is perhaps helpful to address the most basic example, that
of holomorphic curves $f\colon M\to \H\P^{1}$. Recall that the projective
geometry of the quaternionic projective
line $\H\P^1$ is the conformal geometry of the 4-sphere $S^4$. As discussed in
\cite{bflpp}, immersed holomorphic curves in this setting are precisely the conformal
immersions into the 4-sphere. Moreover, the Willmore energy of the holomorphic
curve $f$ is the classical Willmore energy 
\[
W(f)=\int_{M}H^2-K-K^{\perp}\,,
\]
where $H$ is the mean curvature vector, $K$ the Gauss curvature, and $K^{\perp}$ the
curvature of the normal bundle of $f$, all calculated with respect to some conformally flat
metric on $\H\P^1$. Conceptually, it is helpful to think about holomorphic curves
into $\H\P^1$ as the ``meromorphic functions'' of the quaternionic holomorphic theory.
Curves with zero Willmore energy are meromorphic functions in the usual sense,
i.e., holomorphic maps $f\colon M\to\C\P^1$.
\end{example} 
For any rank $k$ subbundle $V\subset H$ we have the perpendicular rank $n-k$ subbundle
$V^{\perp}\subset H^{-1}$ of the dual bundle. $V^{\perp}$ consists of all linear
forms on $H$ which vanish on $V$. The pairing  between $\alpha\in V^{\perp}$ and  
$[\psi]\in H/V$ given by
\begin{equation}\label{eq:perp_identify}
<\alpha,[\psi]>:=<\alpha,\psi>
\end{equation}
identifies $H/V=(V^{\perp})^{-1}$. The derivative 
$\delta^{\perp}\in\Omega^{1}(\text{Hom}(V^{\perp},H^{-1}/V^{\perp}))$ then 
becomes the negative adjoint of the derivative of $V$, i.e.,
\begin{equation}\label{eq:perp_der}
\delta^{\perp}=-\delta^{*}\,.
\end{equation}
To see this, take $\alpha\in\Gamma(V^{\perp})$ and $\psi\in\Gamma(V)$ so that
\begin{align*}
0&=d<\alpha,\psi>=<\nabla\alpha,\psi>+<\alpha,\nabla\psi>\\
&=<\pi^{\perp}\nabla\alpha,\psi>+<\alpha,\pi \nabla\psi>\\
&=<\delta^{\perp}\alpha,\psi>+<\alpha,\delta\psi>\,,
\end{align*}
where we used the same notation $\nabla$ for the dual connection on $H^{-1}$.

If $V\subset H$ is a holomorphic curve, then generally $V^{\perp}\subset H^{-1}$ will not 
be holomorphic (in contrast to the complex holomorphic setup).
\begin{lemma}\label{lem:left_holo}
Let $V\subset H$ be a subbundle and $\delta\in\Omega^{1}(\text{Hom}(V,H/V))$ its derivative. 
Then $V^{\perp}\subset H^{-1}$ is a holomorphic curve 
if and only if there exists a complex structure $J$ on $H/V$ so that
\begin{equation}\label{eq:left_holo}
*\delta=J\delta\,.
\end{equation}
%In other words, $\delta\in\Gamma(K\text{Hom}(V,H/V))$ is left $K$.
\end{lemma}
\begin{proof}
By definition, $V^{\perp}\subset H^{-1}$ is a holomorphic curve if there exists
a complex structure $J^{\perp}$ on $V^{\perp}$ so that 
\[
*\delta^{\perp}= \delta^{\perp}J^{\perp}\,.
\]
Dualizing this relation and using \eqref{eq:perp_der} we obtain the lemma
with $J=(J^{\perp})^{*}$.
\end{proof}
In particular, both $V$ and $V^{\perp}$ are holomorphic curves if and only if
there exist complex structures on $V$ and $H/V$ such that 
\begin{equation}\label{eq:holo}
*\delta=J\delta=\delta J\,,
\end{equation}
in which case $\delta\in\Gamma(K\text{Hom}_{+}(V,H/V))$.

Now assume that the flat rank $n$ bundle $H$ has a complex structure
$J\in\Gamma(\text{End}(H))$, $J^2=-1$. According to
\eqref{eq:nabla_decomposition}, the flat connection $\nabla$ on $H$ decomposes into
\[
\nabla=\nabla'+\nabla''=(\del+A)+(\delbar+Q)\,,
\]
with $\nabla''=\delbar+Q$ a holomorphic structure on $H$. Then 
the dual connection $\nabla$ on $H^{-1}$ decomposes, with respect to the dual complex structure,
into
\begin{equation}\label{eq:nabla_dual}
\nabla=(\del-Q^{*})+(\delbar-A^{*})\,,
\end{equation}
so that $\delbar-A^{*}$ is a holomorphic structure on $H^{-1}$.
If $V\subset H$ is a $J$-stable subbundle then $V$, $H/V$ and hence  
$V^{\perp}=(H/V)^{-1}$ have induced complex structures, all of which we again denote by $J$.
By definition of these complex structures, the canonical projections
$\pi$, $\pi^{\perp}$ and the identification $V^{\perp}=(H/V)^{-1}$ are complex linear.
Now it makes sense to compare the
two notions of holomorphicity of $V$ as a subbundle or as a curve.
\begin{lemma}\label{lem:holo_compare}
Let $J$ be a complex structure on $H$ and $V\subset H$ a $J$ stable
subbundle. Then $V\subset H$ is a holomorphic subbundle with respect to the
holomorphic structure $\nabla''=\delbar+Q$ if and only if $V^{\perp}$ is a holomorphic curve
with respect to the complex structure $J$. The holomorphic structure 
\eqref{eq:holo_structure_induced} on 
$(V^{\perp})^{-1}=H/V$ is the one induced by $\nabla''$. 

The Willmore
energy of the holomorphic curve $V^{\perp}$ is given by
\[
W((V^{\perp})^{-1})=2\int<Q\wedge *Q\,_{H/V}>=2\int<Q^{*}\wedge *Q^{*}\,_{|V^{\perp}}>\,.
\] 
\end{lemma}
\begin{proof}
Recall \eqref{eq:holo_curve_mixed} that $V^{\perp}\subset H^{-1}$ is a holomorphic curve
if and only if the mixed structure 
\[
\hat{D}^{\perp}=\tfrac{1}{2}(\nabla+*\nabla J)=\delbar-Q^{*}
\]
leaves $V^{\perp}$ invariant. But this is equivalent \eqref{eq:perp_identify} to $\delbar+Q$
stabilizing $V\subset H$, i.e., to $V$ being a holomorphic subbundle of $H$.
Lemma~\ref{lem:mixed_holomorphic} then shows that the holomorphic structure $D$ on 
$(V^{\perp})^{-1}=H/V$
from Theorem~\ref{thm:curve_holo_structure} is given by $\delbar+Q$ on $H/V$.

If we denote by $Q_{H/V}$ the Hopf field of $D$ then \eqref{eq:Will_energy_curve} gives
\[
W((V^{\perp})^{-1})=2\int<Q\wedge *Q_{\,H/V}>=2\int<Q^{*}\wedge *Q^{*}_{\,|V^{\perp}}>\,.
\]
We used that
$<Q_{H/V}\wedge *Q_{H/V}>=<Q\wedge *Q_{\,H/V}>$ which holds since $Q$ stabilizes $V\subset H$.
The second equality follows from the fact that $Q^{*}$ stabilizes $V^{\perp}$ and
$(Q_{H/V})^{*}=Q^{*}_{|V^{\perp}}$.
\end{proof}
Putting the above two lemmas together we get
\begin{corollary}\label{cor:holo_compare}
Let $J$ be a complex structure on $H$ and $V\subset H$ a $J$ stable subbundle.
Then we have the following equivalent statements:
\begin{enumerate}
\item
$V$ and $V^{\perp}$ are holomorphic curves with respect to the induced $J$'s.
\item
$V\subset H$ and $V^{\perp}\subset H^{-1}$ are holomorphic subbundles
with respect to $\delbar +Q$ and $\delbar -A^{*}$.
\item
$*\delta=J\delta=\delta J$.
\item
$*\delta^{\perp}=J\delta^{\perp}=\delta^{\perp} J$.
\end{enumerate}
In particular, the Willmore energies of the holomorphic curves $V$ and $V^{\perp}$ are given by
\[
W(V^{-1})=2\int<A\wedge *A\,_{|V}> \qquad \text{and}\qquad W((V^{\perp})^{-1})=2\int<Q\wedge *Q\,_{H/V}>\,.
\]
\end{corollary}
We conclude this section by briefly sketching 
the relation between holomorphic curves and certain
{\em twistor lifts}. If $H$ is an $n$-dimensional quaternionic vector space we use the notation
$(H,i)$ do indicate that we consider $H$ as a $2n$-dimensional complex vector space via 
the complex structure given
by right multiplication by the quaternion $i$ on $H$. Let
\[
G_k^{*}(H,i)=\{W\subset (H,i)\,;\, \dim_{\C}W=k\;\text{and}\;W\cap Wj=0\}\subset G_k(H,i) \cong
G_k(\C^{2n})\,,
\]
which is an open subset of the Grassmannian of complex $k$-panes in $(H,i)$, and note that 
\[
G_1^{*}(H,i)=G_1(H,i)\cong\C\P^{n-1}\,.
\]
The familiar {\em twistor projection}
\begin{equation}\label{eq:twistor_proj}
\pi\colon  G_k^{*}(H,i)\to G_k(H)
\end{equation}
maps a complex $k$-plane $W\subset (H,i)$ to the quaterionic $k$-plane 
$V=W\oplus Wj\subset H$, i.e.,
\[
\pi(W)=W\oplus Wj\,.
\]
The fiber $\pi^{-1}(V)=G_k^{*}(V,i)$ consists of all complex $k$-planes $W\subset (V,i)$
such that $V=W\oplus Wj$. Any complex structure $J\in\text{End}(V)$, $J^2=-1$,
defines a splitting of $V$ into $\pm i$-eigenspaces of $J$, so that the 
fiber  $\pi^{-1}(V)=G_k^{*}(V,i)$ can also be viewed as the space of all
complex structures on $V$ compatible with the quaternionic structure. Thus we can
define the {\em twistor lift} of a quaternionic holomorphic curve $V$ in $G_k(H)$:
by Definition~\ref{def:holo_curve} $V$ has a complex structure $J$ so 
that we can
decompose $V=W\oplus Wj$
into the $\pm{i}$ eigenspaces of $J$ on $V$. Then $W$ is a smooth curve in $G_k^{*}(H,i)$
with $\pi(W)=V$, the {\em twistor lift} of the quaternionic holomorphic curve $V$. 
We denote by
\[
\delta_W=\pi_{W}\nabla
\]
the derivative of the curve $W$ where $\pi_W\colon H\to H/W$ is the 
complex linear quotient projection. Note that $\nabla$ is a flat complex connection
on $(H,i)=\C^{2n}$ since $\nabla$ is quaternionic linear. 
\begin{lemma}\label{lem:twistor_char}
Let $W$ be a smooth curve in $G_k^{*}(H,i)$ and let $V=\pi(W)$ be the twistor projection. 
Then $V$ is a quaternionic holomorphic
curve if and only if $W$ has vertical $\delbar$ derivative, i.e.,
$\delta_W''\in\Gamma(\bar{K}\text{Hom}_{\C}(W,V/W))$. 

In this situation 
\[
\delta_W''=-Q^{*}\,,
\]
where $Q\in\Gamma(\bar{K}\text{End}_{-}(V^{-1}))$
is the Hopf field \eqref{eq:Q} of the holomorphic structure induced by $V$ given in 
Theorem~\ref{thm:curve_holo_structure}. 

Moreover, if the complex structure $J$ on $V$ is the restriction of a complex structure 
$J$ on $H$
then 
\[
\delta_W''=A_{|V}\,,
\]
where $\nabla=\del+\delbar+A+Q$ is the usual decomposition \eqref{eq:nabla_decomposition} of
the flat connection on $H$.

In particular, the  quaternionic holomorphic
curves $V$ in $G_k(\H^{n})$ with zero Willmore energy are the twistor projections $V=\pi(W)$ of 
complex holomorphic curves $W$ in $G_k(\C^{2n})$ with $W\cap Wj=0$. Note that 
the last condition
is vacuous for $k=1$, so that quaternionic holomorphic
curves $L$ in $\H\P^n$ with zero Willmore energy correspond to complex holomorphic curves
$E$ in $\C\P^{2n+1}$ via the twistor projection.
\end{lemma}

\begin{proof}
Let $\hat{D}=\tfrac{1}{2}(\nabla+*\nabla J)$ and take $\psi\in\Gamma(W)$, i.e., $J\psi=\psi i$, 
then we get
\[
\delta_W''\psi=\pi_W \hat{D}\psi\,.
\]
Therefore $\delta_W''$ takes values in $\text{Hom}_{\C}(W,V/W)$ if and only if $\hat{D}$ leaves
$V$ invariant which, by \eqref{eq:holo_curve_mixed}, is equivalent to $V\subset H$ being
a holomorphic curve. But $G_k^{*}(V,i)$ is open inside $G_k(V,i)$ so that  
the tangent space to the fiber $\pi^{-1}(V)$ at $W$ is given by
$\text{Hom}_{\C}(W,V/W)$.

Since $V/W=Wj$ the quotient projection $\pi_W\colon V\to V/W$ is projection onto the
$-i$-eigenspace of $J$ and thus
\[
\delta_W''=\pi_W\hat{D}=\hat{D}_{-}=-Q^{*}\,.
\]
Note that $Q^{*}$ is a section of
$K\text{End}_{-}(V)$ which, by \eqref{eq:homs}, is the same as
$K\text{Hom}_{\C}(\bar{W},W)=\bar{K}\text{Hom}_{\C}({W},\bar{W})=
\bar{K}\text{Hom}_{\C}({W},Wj)=\bar{K}\text{Hom}_{\C}(W,V/W)$.

Let $V\subset H$ be $J$-stable for some complex structure on $H$ and
assume $V$ is holomorphic curve with respect to the induced complex structure $J$ on $V$.
Then, by Lemma~\ref{lem:holo_compare}, the holomorphic structure on $V^{-1}$ is given 
$\delbar-A^{*}$
and thus
\[
\delta_W''=A_{|V}\,.
\]    
\end{proof}
\subsection{Kodaira correspondence}\label{subsec:Kodaira}
The extrinsic geometry of holomorphic curves in Grassmannians
and the geometry of holomorphic vector bundles are related by the Kodaira correspondence.
To discuss this we need the notion of base points.
\begin{definition}\label{def:base_points}
Let $V$ be a quaternionic holomorphic vector bundle over a compact Riemann surface $M$.
A linear subspace $H\subset H^{0}(V)$ of holomorphic sections of $V$ is called
a {\em linear system}.
A point $p\in M$ is said to be a {\em base point} for $H$ if the evaluation map
\[
ev_p\colon H\to V_p\,,\qquad ev_p(\psi)=\psi_p
\]
is not surjective. Thus, $H$ is {\em base point free} if the bundle homomorphism
\[
ev\colon {H}\to V
\]
is surjective.
\end{definition}
We already have seen in Theorem~\ref{thm:curve_holo_structure}
that given a  holomorphic curve $f\colon M\to G_k(\H^n)$ the dual $V^{-1}$ of the induced 
bundle $V=f^{*}\Sigma\subset {\H^{n}}$ has a unique holomorphic structure such that linear
forms $\alpha\in (\H^{n})^{-1}$ restrict to holomorphic sections $\alpha_{|V}\in H^{0}(V^{-1})$.
If the curve $f$ is full, i.e., not contained in any lower dimensional Grassmannian, then
$(\H^{n})^{-1}$ injects into $H^{0}(V^{-1})$, and we obtain an $n$-dimensional linear system 
$H\subset H^{0}(V^{-1})$ without base points. Notice that
$H\subset H^{0}(V^{-1})$ is invariant under projective transformations of $f$. 

Thus we have
assigned to every holomorphic curve, in a projectively invariant way, a 
quaternionic holomorphic vector bundle and a base point free
linear system.
 
Conversely, we start with a rank $k$ holomorphic vector bundle $V^{-1}$ over $M$ and an 
$n$ dimensional subspace $H\subset H^{0}(V^{-1})$ of holomorphic sections without base points.
Since the evaluation map 
$
ev\colon {H} \to  V^{-1}
$ 
is a smooth surjective bundle homomorphism, the dual bundle map
\[
ev^{*}\colon V\to{H}^{-1}
\]
is injective and defines the smooth curve $f\colon M\to G_k(H^{-1})$,
where $f(p)=(ev_p)^{*}(V_p)$, which is full. We now want to show that this curve is
in fact holomorphic with respect to the dual complex structure $J$
coming from $V^{-1}$. Let us denote the evaluation pairing on $H$ by
$<\,,>_{H}$. Note that if $\alpha\in H$ and $\psi\in \Gamma(V)$ then
\[
<ev^{*}\psi,\alpha>_{H}=<\alpha,\psi>\,,
\]
where the latter is the evaluation pairing on $V$. Let $\hat{D}$ be the mixed structure
on $V$ coming from the holomorphic structure $D$ on $V^{-1}$ via 
Lemma~\ref{lem:mixed_holomorphic}. If $\alpha\in H\subset H^{0}(V^{-1})$ and
$\psi\in\Gamma(V)$, then the
product rule \eqref{eq:mixed_holomorphic} implies that 
\begin{align*}
<ev^{*}\hat{D}\psi,\alpha>_{H}&=
<\alpha,\hat{D}\psi>\\
&=\tfrac{1}{2}(d<\alpha,\psi>+*d<J\alpha,\psi>)\\
&=\tfrac{1}{2}(d<ev^{*}\psi,\alpha>_{H}+*d<ev^{*}J\psi,\alpha>_{H}\\
&=<\tfrac{1}{2}(\nabla ev^{*}\psi+*\nabla ev^{*}J\psi),\alpha>_{H}\,,
\end{align*}
since $\nabla\alpha=0$ ($\alpha$ is just a fixed vector in $H$). Thus,
\[
\hat{D}_{ev^{*}(V)}=ev^{*}\hat{D}
\]
which shows, using  \eqref{eq:holo_curve_mixed}, that $ev^{*}(V)\subset H^{-1}$ is a holomorphic
curve. The inclusion $(H^{-1})^{-1}$ into the
space of holomorphic sections of $ev^{*}(V)^{-1}=V^{-1}$ has image $H$, and
is given by the canonical identification of $H=(H^{-1})^{-1}$.
Let us summarize the above discussion:
\begin{theorem}\label{thm:Kodaira}
Let $M$ be a compact Riemann surface. Then there is a correspondence
between the following objects:
\begin{enumerate}
\item quaternionic holomorphic rank $k$ vector bundles $V$ with a base point 
free $n$-dimensional
linear system $H\subset H^{0}(V)$, and
\item projective equivalence classes of full quaternionic holomorphic curves 
$f\colon M\to G_k(\H^n)$ into 
$k$ plane Grassmannians of an $n$ dimensional quaternionic vector space.
\end{enumerate}
\end{theorem}
\begin{example}\label{ex:standard_2}
An obvious application of the Kodaira correspondence is the following: let $L^{-1}$ be
a quaternionic holomorphic line bundle over $M$ with a $2$-dimensional linear system
$H\subset H^{0}(L^{-1})$. If $H$ has no base points, we get
a holomorphic curve $f\colon M\to \H\P^1$. 
If this curve is immersed, then we have seen in Example~\ref{ex:standard_1}
that $f$ is a conformal immersion into the $4$-sphere. To check whether $f$ is immersed,
one has to study the vanishing of the differential 
$\delta\in \Gamma(\text{Hom}(L,{\H^2}/L)K)$.

More generally, if the linear system $H$ is $n+1$ dimensional and base point free, 
we get a holomorphic curve
$f\colon M\to \H\P^n$. A standard construction in the complex holomorphic setting is to assign to
such a curve its Frenet flag \cite{bible} consisting of the various
higher osculating curves. This construction is more subtle in 
the quaternionic case, since in general we cannot extend holomorphically
across zeros of the derivatives. A necessary step in the understanding of this
construction is to study the zeros of quaternionic 
holomorphic sections and
their higher order derivatives. To carry this program through,
we will work intrinsically with holomorphic jets.
\end{example} 

%%% Local Variables: 
%%% mode: latex
%%% TeX-master: "willmore"
%%% End: 

\section{Holomorphic jet bundles}\label{sec:jets}
In this section we  study the geometry of the jet complex  
of a first order linear differential equation
on  a complex vector bundle 
over a Riemann surface whose symbol is that of $\delbar$. The quaternionic symmetry is
auxiliary for these considerations and will be added at the end
as an important special case.

For our geometric point of view
it is advantageous to
give an axiomatic development of the jet complex. The more familiar
notion of jets as successive higher derivatives \cite{Palais} will become apparent
in our construction of an explicit model of the jet complex.

\subsection{The holomorphic jet complex}\label{subsec:holo_jet_complex}
Throughout $V$ will be a vector bundle with complex structure $J$ over a Riemann surface $M$.
\begin{definition}\label{def:holo_structure}
A {\em holomorphic structure} on $V$ is a real linear map 
\[
D\colon \Gamma(V)\to\Gamma(\bar{K}V)
\]
satisfying the Leibniz rule
\[
D(f\psi)=fD\psi+\tfrac{1}{2}(df\psi+*df J\psi)
\]
over {\em real} valued functions $f\colon M\to \R$. We denote by
\[
H^{0}(V)=\ker D\subset \Gamma(V)
\]
the real linear subspace of {\em holomorphic sections} of $V$.

A bundle map $T\colon V\to \tilde{V}$ between two bundles with holomorphic structures
is {\em holomorphic} if $T$ is complex linear and satisfies $\tilde{D}T=TD$. 
\end{definition}
\begin{remark}
The standard examples of holomorphic structures  are {\em quaternionic}
holomorphic structures, i.e., $D$ is quaternionic linear, and {\em complex} 
holomorphic structures, i.e.,  $D$ is complex linear. 
As we already have seen, the latter are a special
case of the former. 

Most constructions in section~\ref{sec:quat_holo_geom} have analogues for the more general
case of holomorphic structures. 
\end{remark}
Decomposing $D$ into $J$-commuting and 
anticommuting parts, we get
\[
D=\delbar+Q\,,
\]
where $\delbar$ is a complex holomorphic structure on $V$ and $Q$ is a section
of $\bar{K}\text{End}_{-}(V)$. A holomorphic bundle map $T\colon V\to \tilde{V}$
then is complex holomorphic and intertwines the respective $Q$'s.

We now come to the axiomatic development
of the {\em holomorphic} jet complex. As usual, $V$ denotes a  rank $r$ complex
vector bundle with complex structure $J$ and holomorphic structure $D$. 
\begin{axiom}\label{ja1}
There exist a sequence of real vector bundles $V_n$ of rank $2r(n+1)$ and linear
surjective bundle homomorphisms $\pi_n\colon V_n\to V_{n-1}$, $n\in\N$, 
with $V_0=V$ and $\pi_0=0$.
\end{axiom}
We denote by $N_n:=\ker\pi_n$ the rank $2r$ subbundle given by the kernel of
$\pi_n$, so that $V_n/N_n=V_{n-1}$.
If the context is clear, we will omit the subscripts and just write 
$\pi\colon V_n\to V_{n-1}$.
For a given $n\in\N$ and $k\leq n$ the subbundles
\[
F_k:=F_{n,k}:=\ker\pi^{k+1}\subset V_n
\]
have rank $2r(k+1)$ and give the flag 
\[
V_n=F_{n}\supset F_{n-1}\supset\cdots\supset F_{1}\supset F_{0}
\]
in $V_n$ with $V_n/F_{n-1}=V$ and $N_n=F_0$.
\begin{axiom}\label{ja2}
For each $n\in \N$  there exists a linear map
\[
d_n\colon \Gamma(V_n)\to\Omega^{1}(V_{n-1})\,,\qquad d_0=0
\]
satisfying the Leibniz rule
\[
d_n(f\psi)=fd_n\psi+df\pi_{n}(\psi)
\]
over real valued functions. The $d$'s commute with the $\pi$'s, so that
\[
\pi_{n-1}d_n=d_{n-1}\pi_n\,.
\]
Moreover, the $d$'s are ``flat", i.e., $d_{n-1}d_{n}=0$.
\end{axiom}
To understand the last condition, observe that $d$ can be uniquely extended to
an ``exterior derivative" on 
$1$-forms (in fact, $k$-forms) by the usual formalism: for $\omega\in\Omega^{1}(V_n)$
and vector fields $X$,$Y$ on $M$,
\[
d\omega_{X,Y}=d_{X}\omega(Y)-d_{Y}\omega(X)-\pi\omega([X,Y])\,.
\]
\begin{comment}
Conceptually, one should think of the $d$'s as a ``flat" connection on the
inverse limit of the $V_n$'s with $d\Gamma(V_n)\subset \Omega^{1}(V_{n-1})$.
\end{comment}
We shall see below that jets in the kernel of $d_n$  
are $n$-th derivatives of holomorphic sections of $V$.

It is important to notice that, due to the Leibniz rule, 
the maps $d_n$ become tensorial when restricted to the kernels $N_n$. 
Since $d$ and $\pi$ commute,
$d$ maps kernels into kernels and we obtain bundle maps
\[
\delta_{n}:={d_{n}}_{|N_n}\colon N_n\to T^{*}M\otimes N_{n-1}\,,
\]
which is to say that $\delta_n\in\Omega^{1}(\text{Hom}(N_n,N_{n-1}))$.
\begin{axiom}\label{ja3}
The kernels $N_n=\ker\pi_n$ have complex structures $J_n$ (and thus become
complex bundles of rank $r$), and  
\begin{equation}\label{eq:jet_delta}
\delta_n\colon N_n\to KN_{n-1}
\end{equation}
are complex linear bundle isomorphisms for $n\geq 1$. On $N_0=V$ we have $J_0=J$ and $\delta_0=0$.
\end{axiom}
To formulate the last axiom, which will tie the given holomorphic structure $D$ to the 
jet complex, note that the $\bar{K}$-part $d_1''=\tfrac{1}{2}(d_1+*J d_1)$ is  defined since 
it only involves the complex structure $J$ on $V_0=V$. 
\begin{axiom}\label{ja4}
\[
D\pi_1=d_1''\,.
\]
\end{axiom}
This last axiom assures that we only consider jets of holomorphic sections of $V$.

\begin{definition}\label{def:holo_jet_complex}
Let $V$ be a complex vector bundle with holomorphic structure $D$.
A sequence of real vector bundles $V_n$ with maps $d_n$ and $\pi_n$
satisfying the above axioms is called a {\em holomorphic jet complex} of $V$, 
and will be denoted by $\mathcal{V}$.

Let $V$ and $\tilde{V}$ be two complex vector bundles with holomorphic structures
$D$ and $\tilde{D}$ and
holomorphic jet complexes $\mathcal{V}$ and $\tilde{\mathcal{V}}$. 
A {\em homomorphism} of jet complexes
$T\colon \mathcal{V}\to\tilde{\mathcal{V}}$ is given by a sequence of real linear bundle homomorphisms
$T_n\colon V_n\to \tilde{V}_n$ commuting with $\pi$ and $d$. In addition, $T_0$ has
to be holomorphic.
\end{definition}

An immediate consequence of the axioms for a jet complex is the unique prolongation property
of holomorphic sections:
\begin{lemma}\label{lem:prolongation}
Let $V$ be a complex vector bundle with holomorphic structure $D$ and 
$\mathcal{V}$ a holomorphic jet complex of $V$.
Let $\psi\in\Gamma(V_n)$ with $D\psi=0$ for $n=0$ or with  $d\psi=0$ 
for $n\geq 1$. 
Then there exists a
unique section $\tilde{\psi}\in\Gamma(V_{n+1})$ with $\pi\tilde{\psi}=\psi$ and
$d\tilde{\psi}=0$.
\end{lemma}
\begin{proof}
Choose any lift $\hat{\psi}\in\Gamma(V_{n+1})$ of $\psi$, i.e., $\pi\hat{\psi}=\psi$.
Then 
\[
0=d\psi=d\pi\hat{\psi}=\pi d \hat{\psi}\,,
\]
so that $d \hat{\psi}$ is an $N_n$ valued $1$-form. Since $d^2=0$ and $d$ restricts
to $\delta$ on $N$'s, we obtain
\[
0=d^2 \hat{\psi}=\delta\wedge d \hat{\psi}\,,
\]
which implies that $d''\hat{\psi}=0$ in case $n\geq 1$. 
For $n=0$ we use axiom~\ref{ja4}, instead of $d^2=0$, to get 
\[
d''\hat{\psi}=D\pi\hat{\psi}=D\psi=0\,.
\]
If $\phi\in\Gamma(N_{n+1})$, then
$\tilde{\psi}:= \hat{\psi}-\phi$ is also a lift of $\psi$. Moreover, the equation
\[
0=d\tilde{\psi}=d\hat{\psi}-d\phi=d'\hat{\psi}-\delta\phi
\]
has the unique solution $\phi=\delta^{-1}d\hat{\psi}$.  
\end{proof}
\begin{corollary}\label{cor:prolongation}
For each $n\in\N$ there is an injective real linear map,
the {\em $n$-th prolongation} of holomorphic sections,
\[
P_n\colon H^{0}(V)\to\Gamma(V_n)\,,\qquad P_0=id\,, 
\]
satisfying
\[
\pi^{k}P_n=P_{n-k}\qquad \text{and}\qquad d P_n=0\,.
\]
\end{corollary}
We now turn to the uniqueness and existence of the holomorphic jet complex. For the
uniqueness part, we will prolong bundle homomorphisms using the same ideas as in the 
previous lemma.
\begin{theorem}\label{thm:uniqueness}
Let $V$, $\tilde{V}$ be complex vector bundles with holomorphic structures $D$, $\tilde{D}$ 
and holomorphic jet complexes $\mathcal{V}$,$\,\tilde{\mathcal{V}}$. 
Let $T_0\colon V\to\tilde{V}$ be holomorphic. Then there is a unique homomorphism 
$T\colon \mathcal{V}\to\tilde{\mathcal{V}}$ of
jet complexes extending $T_0$.

In particular, the holomorphic jet complex is unique.
\end{theorem}
\begin{proof}
We proceed by induction.  
Assume that we have constructed extensions $T_k$ for $k\leq n$ with the 
desired properties. Let $\hat{T}\colon V_{n+1}\to\tilde{V}_{n+1}$ be an extension of
$T_n$ such that $\tilde{\pi}\hat{T}=T_n\pi$. Any other extension is of the form
$T_{n+1}=\hat{T}-B$, where $B\colon V_{n+1}\to\tilde{N}_{n+1}$. 
The condition $\tilde{d}T_{n+1}=T_n d$ then becomes
\begin{equation}\label{eq:plus_B}
\tilde{\delta} B=\tilde{d}\hat{T}-T_n d\,.
\end{equation}
Applying $\tilde{\pi}$ we get
\[
\tilde{\pi}(\tilde{d}\hat{T}-T_n d)=(\tilde{d}T_n-T_{n-1} d)\pi=0
\]
by the induction assumption. Thus the right hand side of \eqref{eq:plus_B}
is a $\text{Hom}(V_{n+1},\tilde{N}_{n})$ valued $1$-form. Moreover, this $1$-form has
no $\bar{K}$-part since, for $n\geq 1$,
\[
\tilde{\delta}\wedge(\tilde{d}\hat{T}-T_n d)=\tilde{d}^2\hat{T}-\tilde{d}T_n d=-T_n d^2=0\,.
\]
On the lowest level, $n=0$, we obtain the same conclusion
\[
\tilde{d}''\hat{T}-T_0 d''=\tilde{D}\tilde{\pi}\hat{T}-T_{0}D\pi=(\tilde{D}\hat{T}-\tilde{D}T_{0})\pi=0
\]
due to axiom~\ref{ja4} and the assumptions on $T_0$.
Since $\tilde{\delta}\colon \tilde{N}_{n+1}\to K\tilde{N}_n$ is invertible, we now can
solve \eqref{eq:plus_B} uniquely for $B$.

Let $T$ and $\hat{T}$ be two extensions of $T_0$ to jet complex homomorphisms. Then their difference
is $\tilde{N}$ valued and in the kernel of $\tilde{\delta}$, and thus it is zero. 
This shows the uniqueness of 
extensions.
\end{proof}
\begin{remark}
We only used the invertibility of $\tilde{\delta}$ on the target jet complex.
\end{remark}

To construct an explicit model of the holomorphic jet complex and for further
computations it is useful
to work with ``adapted" connections.  
\begin{lemma}\label{lem:jet_nabla}
Let $\mathcal{V}$ be a holomorphic jet complex of a complex vector bundle $V$ with
holomorphic structure $D$ and fix $n\in\N$.  
\begin{enumerate}
\item
If $T\colon V_{n}\to V_{n+1}$ is a section of $\pi$, i.e., $\pi T= id$, then
\begin{equation}\label{eq:dP}
\nabla=d T
\end{equation}
is a connection on $V_{n}$ satisfying
\begin{equation}\label{eq:jet_nabla}
\pi\nabla=d\;\;\text{and}\;\; \pi R^{\nabla}=0\,.
\end{equation}
\item
Conversely, any connection $\nabla$ on $V_n$ satisfying \eqref{eq:jet_nabla} is of the form
\eqref{eq:dP} for a suitable section $T$ of $\pi$.
\item
Let $\nabla$ be an adapted connection. Then $\nabla+\omega$ is adapted, if and only if
$\omega$ is a section of $K\text{Hom}(V_n,N_n)$.  
\end{enumerate}
\end{lemma}
A connection on a jet bundle $V_n$ with the above properties will be called an
{\em adapted} connection. Since a section $T$ of $\pi$ is given by any choice of
splitting
\[
V_{n+1}=N_{n+1}\oplus H_n
\]
we have adapted connections on all the jet bundles $V_n$.  
\begin{proof}
Given $T$ we clearly have 
\[
\pi\nabla=\pi d T= d\pi T=d\,
\]
by the properties of $d$ and $\pi$. This immediately 
implies $\pi d^{\nabla}=d$ on $1$-forms and thus
\[
\pi R^{\nabla}=\pi d^{\nabla}dT=d^{2}T=0\,.
\]
For the converse we note that $\omega:=d-\nabla\pi$ is linear over real valued
functions and thus defines a $1$-form with values in 
$\text{Hom}(V_{n+1},V_n)$. But
\[
\pi\omega=\pi d-\pi\nabla\pi=\pi d-d\pi=0
\]
so that $\omega$ is $\text{Hom}(V_{n+1},N_n)$ valued. Moreover,
\[
\delta\wedge\omega=d^2-\pi d^{\nabla}\nabla\pi=-\pi R^{\nabla}\pi=0\,
\]
which shows that $\omega$ is a section of $K\text{Hom}(V_{n+1},N_n)$.
Notice also that $\omega=\delta$ when restricted to $N_{n+1}$. Thus we
have the splitting
\[
V_{n+1}=N_{n+1}\oplus \ker\omega
\] 
defining $T$ with $\pi T=id$, which gives 
\[
dT-\nabla=(d-\nabla\pi)T=\omega T=0\,.
\]
Assume now that $\omega=\tilde{\nabla}-\nabla$ is the difference of two adapted
connections. Then $\pi\omega=0$ so that $\omega\in \Omega^{1}(\text{Hom}(V_n,N_n))$. 
Using again  that $d=\pi d^{\nabla}$ for an adapted connection, we get
\[
0=\pi R^{\tilde{\nabla}}=\pi( R^{\nabla}+d^{\nabla}\omega +\omega\wedge\omega)=\pi d^{\nabla}\omega=
\delta\wedge\omega\,, 
\]
which shows that $\omega$ has no $\bar{K}$-part. 
\end{proof}
\begin{remark}
As the proof indicates, we never had to use the holomorphic structure $D$,
i.e., axiom~\ref{ja4}, in the construction
of adapted connections. This reflects the fact that those connections exist
for general jet complexes of vector bundles over any manifold. But we choose
not to present this subject in that generality, since we have no need for it
in the present context.
\end{remark}
To give an explicit model of the holomorphic jet complex, we  follow the idea that the 
successive jet bundles are comprised of successive higher derivatives. Let
$V$ be a complex vector bundle with holomorphic structure $D$.
Put
\[
V_1:=V\oplus KV
\]
with $\pi_1$ projection onto the first factor. Let $\nabla$ be any connection
on $V$ whose $\bar{K}$-part $\nabla''= D$. Define
\[
d_1\colon \Gamma(V_1)\to \Omega^{1}(V)\,,\qquad d_1(\psi,\omega)=\nabla\psi-\omega\,,
\]
so that $1$-jets in the kernel of $d_1$ are in fact first derivatives
of sections in $V$. We clearly have the required Leibniz rule for $d_1$.
Moreover, $N_1=KV$ and thus has a complex structure, 
and the restriction $\delta_1={d_1}_{|N_1}$ is the identity
map. To verify axiom~\ref{ja4}, we note that  $\omega''=0$ and hence  
\[
d_1''(\psi,\omega)=\nabla''\psi= D\pi_1(\psi,\omega)\,.
\]
Now let us assume that we have constructed the jet complex up to level $n$.
We put
\[
V_{n+1}=V_n\oplus KN_n\,,
\]
with $\pi_{n+1}$ projection onto the first factor. Since, by induction, 
$N_n$ has a complex structure
also $N_{n+1}=\ker\pi_{n+1}=KN_{n}$ has a complex structure. 
Let $\nabla$ be an adapted connection on $V_n$, and define
\[
d_{n+1}\colon \Gamma(V_{n+1})\to\Omega^{1}(V_n)
\]
by
\[
d_{n+1}(\psi,\omega)=\nabla\psi-\omega\,.
\]
Clearly, the required Leibniz rule holds.
Since $\pi\nabla=d$ it follows that 
$\pi_{n}d_{n+1}=d_{n}\pi_{n+1}$. 
The restriction of $d_{n+1}$ to $N_{n+1}=\ker\pi_{n+1}$ is again given by the
identity map $\delta_{n+1}\colon$ $ N_{n+1}\to KN_n$. To check $d^2=0$ we 
first note that $\delta_n\wedge \omega=0$ by type. Thus
\[
d_n d_{n+1}(\psi,\omega)=d_n(\nabla\psi+\omega)=d_n\nabla\psi=
\pi_n d^{\nabla}\nabla\psi=\pi_n R\psi=0\,.
\]
This concludes our construction of the holomorphic jet complex of a complex vector bundle
$V$ over a Riemann surface with holomorphic structure $D$.

\subsection{Complex structures on the holomorphic jet complex}
We are now going to discuss 
complex structures on the holomorphic jet complex.
Since $D$ is not complex linear, we cannot use 
Theorem~\ref{thm:uniqueness} to extend $J$ on $V$ canonically to the jet complex.
In fact, there are many ways to choose complex structures $S_n$ on
$V_n$ which, in a certain sense, are adapted to the jet complex:
\begin{definition} 
Let $S$ be a complex structure on the jet complex $\mathcal{V}$ of a complex vector bundle
$V$ with holomorphic structure $D$. We say that  $S$ is {\em adapted}
if the $\pi$'s are complex linear and $S$ induces the given complex structures
$J$'s on $\ker\pi$'s. In other words, $S$ satisfies
\begin{equation}\label{eq:compatible_S}
\pi S=S\pi\qquad \text{and}\qquad  S_{|N}=J\,.
\end{equation}
\end{definition}
Given an adapted complex structure $S$ the holomorphic jet bundles $V_n$ 
become complex vector bundles which, over a compact Riemann surface, have a degree.
It turns out that the degree of $V_n$ is independent
of the choice of adapted complex structure: by the axioms of the holomorphic jet complex
$V_n/N_n=V_{n-1}$
and $N_n=KN_{n-1}$ as complex bundles. Thus \eqref{eq:degmult} implies
\begin{align}\label{eq:jet_degree}
\deg V_n &=\sum_{k=0}^{n}\deg K^{k}V =\sum_{k=0}^{n}k(2g-2)\rank V+\deg V\\
&=(n+1)(n(g-1)\rank V+\deg V)
\end{align}
where $g$ denotes the genus of $M$.

Any adapted $S$ can be used
to construct holomorphic structures on every $V_n$. Since $d_{|N}=\delta$
we see that ${d''}_{|N}=0$, and thus $d''$ induces a linear map 
\begin{equation}\label{eq:D_n}
D_n\colon \Gamma(V_n)\to\Gamma(\bar{K}V_n)\,, \qquad D_n\pi_{n+1}=d_{n+1}''\,.
\end{equation}
Clearly, $D_n$ has the Leibniz rule over
real valued functions. Thus, we have holomorphic structures on every $V_n$
and $D_0=D$. Decomposing into $S$-commuting and anticommuting parts we get
\begin{equation}\label{eq:D_n_decomposed}
D_n=\delbar_n+Q_n\,,
\end{equation}
with complex holomorphic structures $\delbar_n$ on each of the $V_n$ and
sections $Q_n$ of $\bar{K}\text{End}_{-}(V_n)$. From \eqref{eq:D_n} we obtain
\begin{equation}\label{eq:pi_D_n}
\pi^{k} D_n=D_{n-k}\pi^{k}\,,
\end{equation}
and, since $\pi S=S\pi$, also
\begin{equation}\label{eq:Q_n}
\pi^{k} \delbar_n=\delbar_{n-k}\pi^{k}\qquad\text{and}\qquad \pi^{k}Q_n=Q_{n-k}\pi^{k}\,.
\end{equation}
In particular $D_n$, $\delbar_n$ and $Q_n$ stabilize the flag
\[
V_n\supset F_{n-1}\supset\cdots\supset F_1\supset F_0\,.
\]
Also note that \eqref{eq:D_n} and \eqref{eq:pi_D_n} imply
\begin{equation}\label{eq:ker_d}
\ker d_n =\ker D_n\cap\ker d_{n}'\,.
\end{equation}
The prolongation maps $P_n$ of Corollary~\ref{cor:prolongation} satisfy
$d_n P_n=0$ so that
\begin{equation}\label{eq:D_to_D_n}
P_n\colon H^{0}(V)\to H^{0}(V_n)\,.
\end{equation}
In other words, the $n$-th prolongation takes holomorphic sections of $V$
into holomorphic sections of $V_n$. 

The next lemma will show that the $\bar{K}$-part of an adapted connection on $V_n$
gives, as expected,  the holomorphic structure $D_n$:
\begin{lemma}\label{lem:nabla_holo}
Let $S$ be an adapted complex structure on $\mathcal{V}$ and $\nabla=(\del+\delbar)+({A}+{Q})$ 
an adapted
connection on $V_n$ for fixed $n\in \N$.  Then
\[
\nabla''=D_n\,,
\]
and thus also 
\[
\delbar=\delbar_n\qquad \text{and}\qquad Q=Q_n\,.
\]
\end{lemma}
\begin{proof}
Lemma~\ref{lem:jet_nabla} gives us $\nabla=d T$ for some bundle map 
$T\colon V_{n}\to V_{n+1}$ with $\pi T=id$. Using \eqref{eq:D_n} this immediately yields 
\[
\nabla''=d''T=D_n\pi T=D_n\,,
\]
and hence also the corresponding statements for the $S$-commuting and anticommuting parts.
\end{proof}

A more subtle feature of the holomorphic jet complex is the existence
of a canonical adapted complex structure. So far we have said nothing about the
compatibility of $S$ with the maps $d$ on the jet complex. Since $D_n\pi =d''_{n+1}$
we have
\begin{equation}\label{eq:[d'',S]}
2Q_n \pi=(D_n+SD_n S)\pi=d_{n+1}''+Sd_{n+1}''S=S[d_{n+1}'',S]
\end{equation}
so that $S$ and $d''$ commute only in the complex holomorphic setting $Q=0$. On the other hand,
as we will see in the following theorem, the condition $[d',S]=0$ can always
be fulfilled and singles out a unique adapted complex structure $S$ on the holomorphic jet complex. 
\begin{theorem}\label{thm:mcS}
Let $\mathcal{V}$ be the holomorphic jet complex of the holomorphic structure $D$
on the complex vector bundle $V$. Then there exists a unique adapted complex structure
$S$ on $\mathcal{V}$ satisfying
\[
[d',S]=0\,.
\]
This condition is
equivalent to the vanishing of the bundle map $[d'',S]$ on $F_{n-1}=\ker\pi^{n}\subset V_n$.

In particular \eqref{eq:[d'',S]}, the holomorphic structures $D_n=\delbar_n+Q_n$ satisfy 
${Q_n}_{|F_{n-1}}=0$ so that 
the restriction of $D_n$ to the flag $F_k$ induces the complex holomorphic structure
$\delbar_n$.
\end{theorem}
\begin{proof}
Again we proceed by induction. We already have a complex structure on $V=V_0$.
Assume now that we have $S$ up to $V_n$ with the desired properties.
To construct $S_{n+1}$, first choose any complex structure $\hat{S}$ on $V_{n+1}$ 
extending $S_n$, i.e., $\pi_{n+1}\hat{S}=S_n\pi_n$, and inducing $J_{n+1}$ on
$N_{n+1}$.  That $\hat{S}+R$, with $R\colon V_{n+1}\to V_{n+1}$, is another
such extension is equivalent to
\[
\text{im} R\subset N_{n+1} \subset \ker R
\]
and the fact that $R$ and $\hat{S}$ anticommute. 
Now put $S_{n+1}:=\hat{S}+R$ and determine $R$ so that
\begin{equation}\label{eq:S_{n+1}}
0=d'S_{n+1}-S_n d'=d'\hat S-S_n d'-\delta R
\end{equation}
holds. In this relation $d'$ is determined by $S_n$, which 
we already know. Furthermore, 
\[
\pi(d'\hat S-S_n d')=d'\pi\hat{S}-S_{n-1}d'\pi=(d'S_n-S_{n-1}d')\pi=0\,,
\]
which implies that $d'\hat S-S_n d'$ is a section of $\text{Hom}(V_{n+1},KN_n)$.
But $\delta_{n+1}\colon N_{n+1}\to KN_{n}$ is invertible, so we can
solve \eqref{eq:S_{n+1}} uniquely for $R$. The necessary properties for $R$
are now easily checked: restricted to $N_{n+1}$ 
\[
d'\hat{S}-S_n d'=\delta_{n+1}J_{n+1}-J_n \delta_{n+1}=0\,,
\]
since $\delta_{n+1}$ is complex linear. Thus $\text{im}R\subset N_{n+1}\subset \ker R$. Finally,
$R$ and $\hat{S}$ anticommute because
\[
\delta\hat{S}R=
\delta J_{n+1}R=J_n\delta R=S_n d'\hat{S}+d'=
-\delta R\hat{S}\,.
\]
It remains to show that $d'S=Sd'$ is equivalent to the vanishing of $Sd''-d''S$ on $F_{n-1}$.
Since $\pi F_{n-1}=F_{n-2}\subset V_{n-1}$ we conclude inductively that
\[
\pi(Sd''-d''S)_{|F_{n-1}}=(Sd''-d''S)_{|F_{n-2}}=0\,.
\]
Thus, $Sd''-d''S$ takes values in $N_n$ and, to show its vanishing, it suffices to 
calculate $\delta\wedge(Sd''-d''S)$ restricted to $F_{n-1}$: 
\[
\delta\wedge(Sd''-d''S)=dSd''-dd''S=d(Sd''+d'S)=Sd^{2}=0\,.
\]
We used $d^2=d(d'+d'')=0$ and the fact that $d''\colon \Gamma(F_{n-1})\to \Gamma(\bar{K}F_{n-2})$, 
which implies $Sd''=d''S$ on $F_{n-1}$ by induction. The converse is proven similarly. 
\end{proof}
\begin{remark}
The arguments used in the construction of the canonical complex structure $S$
are the same as those in the construction of the mean curvature sphere \cite{bflpp}
of a conformal immersion $f\colon M\to\H\P^1$. We revisit this aspect in more detail
later on, when we will discuss holomorphic curves in $\H\P^n$ and interpret $S$ as 
a congruence of osculating $\C\P^n$'s. 
\end{remark}
\begin{corollary}\label{cor:adapted_nabla_relations}
Let  $\nabla$ be an adapted connection \eqref{eq:jet_nabla} on the 
$n$-th jet bundle of the holomorphic jet complex
of $V$ and $S$ the canonical complex structure. If $\nabla=(\del+\delbar)+({A}+{Q})$
is the usual decomposition \eqref{eq:nabla_decomposition} then
\[
Q_{|F_{n-1}}=0\qquad\text{and}\qquad \pi_n A=0\,.
\]
\end{corollary}
\begin{proof}
From Lemma~\ref{lem:nabla_holo} we have $Q=Q_n$ which implies $Q_{|F_{n-1}}=0$
due to Theorem~\ref{thm:mcS}. On the other hand  \eqref{eq:dJ} gives 
$\nabla S=2*(Q-A)$ which, together with  $\pi\nabla=d$, implies
\[
0=d'S-Sd'=\pi\nabla'S=-2*\pi A\,.
\]
\end{proof}
\begin{remark}\label{rem:A=0}
If $S$ is the canonical complex structure on the holomorphic jet complex of $V$ then
there always exists an adapted connection on $V_n$ with $A=0$: let $\nabla$ be
any adapted connection, then Lemma~\ref{lem:jet_nabla} implies
that $\tilde{\nabla}=\nabla-A$ is also adapted, since $\pi A=0$.
\end{remark} 
\subsection{Zeros of solutions}
Since solutions of $D=\delbar+Q$ cannot vanish to infinite order
they vanish, up to higher order terms, like complex holomorphic sections. Thus, zeros
are isolated and we can assign a vanishing order. 
\begin{definition}\label{def:vanishing}
Let $V$ be a vector bundle over a Riemann surface $M$. A section $\psi\in\Gamma(V)$ 
{\em vanishes to order at least $k$} at 
$p\in M$ if 
\[
|\psi|\leq C|z|^k\,.
\]
Here $z$ is a centered holomorphic coordinate near $p$, $C>0$ some constant, and
$|\psi|$ denotes any norm on the vector bundle $V$. 

We will use the notation $\psi=O(k)$ to indicate that $\psi$ vanishes to
at least order $k$ at $p$. 
\end{definition}
\begin{lemma}\label{lem:orders}
Let $\mathcal{V}$ be the holomorphic jet complex of the holomorphic structure $D$
on the complex vector bundle $V$, and let $\psi\in\Gamma(V_{n})$ with $d\psi=0$. Then
\[
\psi=O(k)\;\;\text{if and only if}\;\; \pi\psi=O(k+1)\,.
\]
\end{lemma}
\begin{proof}
We work with the explicit holomorphic jet complex described above. This means that
$\psi=(\pi\psi,\omega)$, with $\omega$ a section of $KN_{n-1}$. The condition
$d\psi=0$ translates into
\[
\nabla\pi\psi=-\omega\,, 
\]
where $\nabla$ is any adapted connection \eqref{eq:jet_nabla} on the jet complex.
Using trivializations near $p$, we may assume that all occurring sections are vector valued
functions. Moreover, up to zero order terms, $\nabla$ is the directional derivative.
The lemma then follows.
\end{proof}
In the next lemma we analyze the zeros of holomorphic sections.
\begin{lemma}\label{lem:zeros}
Let $V$ be a complex vector bundle over a Riemann surface $M$ with holomorphic
structure $D=\delbar+Q$, and let $\psi\in H^{0}(V)$ be a non-trivial holomorphic section   
Then, at each $p\in M$, there exists a centered holomorphic coordinate $z$ and  
a local nowhere vanishing section $\phi$ of $V$, such that 
\[
\psi=z^{n}\phi+O({n+1})\,.
\]
The integer $n\in\N$  depends only on the section $\psi$.

In particular, non-trivial holomorphic sections have isolated zeros.
\end{lemma}
\begin{definition}\label{def:order}
In the situation of the lemma we denote by
\[
\ord_p\psi:=n
\]
the {\em order} of the zero of $\psi\in H^{0}(V)$ at $p\in M$.
\end{definition}
\begin{proof}
Let $\underline{\psi}$ be a local trivializing frame of $V$ near $p\in M$. 
Then $\psi=\underline{\psi}f$
with $f$ a local $\C^{r}$ valued map. Furthermore,
\[
D\psi=\delbar\psi+Q\psi=(\delbar\underline{\psi})f+
\underline{\psi}\delbar f+(Q\underline{\psi})\bar{f}\,,
\]
where we used that $Q$ is complex antilinear. Putting
\[
\delbar\underline{\psi}=\underline{\psi}\alpha\,,\qquad 
Q\underline{\psi}=\underline{\psi}\beta\,,
\]
for local sections $\alpha$ and $\beta$ of $\bar{K}\mathbf{gl}(r,\C)$, the
equation $D\psi=0$ is locally given by
\begin{equation}\label{eq:analysis}
\delbar f+ \alpha f +\beta \bar{f}=0\,.
\end{equation}
By standard results in analysis (e.g., \cite{HW},\cite{A}),
non-trivial solutions to the above equation cannot vanish to infinite order
(on a connected set). Let $n$ be the order of the smallest non-vanishing derivative of $f$
at $p$, i.e.,
\[
f=\sum_{k=0}^{n} z^{n-k}\bar{z}^{k}a_k+O({n+1})\,,
\]
with $a_k\in \C^{r}$ not all zero. Inserting back into \eqref{eq:analysis} we obtain
\[
\sum_{k=1}^{n} kz^{n-k}\bar{z}^{k-1} a_k d\bar{z}=O({n})\,,
\]
which is only possible if $a_k=0$ for $k=1,\dots,n$. Thus, $f=z^n a_0 + O({n+1})$
and $\psi=z^n\phi+ O({n+1})$ where $\phi=\underline{\psi}a_0$. Clearly,
$n$ depends only on $\psi$ and not on the choice of trivialization
or coordinates. 
\end{proof}
For later applications it will be useful to explicitly know 
the leading order term of prolongations of holomorphic sections.
\begin{lemma}\label{lem:prolongation_order}
Let $P_k\colon H^{0}(V)\to H^{0}(V_k)$ be the prolongation map. If $\psi\in H^{0}(V)$ and 
$\psi=z^n\phi+ O({n+1})$ then 
\[
P_k\psi=z^{n-k}(n(n-1)\cdots(n-k+1)\phi_k+ O(1))\,,
\]
where $\phi_k$ are local sections of $N_k=\ker\pi_k$ such that
\[
\delta^k \phi_k= (-1)^{k}dz^k\phi\,.
\]
In particular, 
$\ord_p\psi$ is characterized as the smallest $n$ such that the
$n$-th prolongation of $\psi$ does not vanish at $p$. 
\end{lemma}
\begin{proof}
Using the explicit representation of the jet complex
$V_{k+1}=V_k\oplus KN_k$
we have
\[
P_{k+1}\psi=(P_k\psi, -\nabla P_k\psi)\,,
\]
where $\nabla$ is an adapted connection \eqref{eq:jet_nabla} on the jet complex. 
To calculate $\nabla P_k\psi$
we decompose $\nabla=\hat{\nabla}+\omega$ so that $\tilde{\nabla}S=0$ for the canonical 
complex structure $S$ of Theorem~\ref{thm:mcS}. Then
\[
\nabla P_k\psi=n(n-1)\cdots (n-k+1)(n-k)z^{n-k-1}dz\,\phi_k+O(n-k)
\]
and thus
\[
P_{k+1}\psi=n(n-1)\cdots (n-k+1)(n-k)z^{n-k-1}(0,-dz\,{\phi}_k)+O({n-k})\,.
\]
The section $\phi_{k+1}:=(0,-dz\,{\phi}_k)$ is $N_{k+1}$ valued and  satisfies inductively
\[
\delta^{k+1}\phi_{k+1}=\delta^{k}(-dz\,{\phi}_k)=(-1)^{k+1}dz^{k+1}\phi\,.
\]
Thus we have shown 
\[
P_k\psi=n(n-1)\cdots (n-k+1)z^{n-k}\phi_k+O(n-k)
\]
or, by taking out $z^{n-k}$,
\[
P_k\psi=z^{n-k}(n(n-1)\cdots (n-k+1)\phi_k+O(1))
\]
in case when $k\leq n$. Note that the $O(1)$ in the last expression generally is not a 
smooth function. For $k\geq n+1$ the last formula reduces to $P_k\psi=z^{n-k} O(1)$ which
is trivially true, since the left hand side is smooth. In fact, the $O(1)$ then has to be a
$O(k-n)$.
\end{proof}
\subsection{The quaternionic holomorphic jet complex}
We now specialize the above considerations to quaternionic holomorphic structures: 
let $V$ be a quaternionic holomorphic vector bundle 
of rank $r$ over a Riemann surface $M$ with complex structure $J$ and quaternionic
holomorphic structure $D=\delbar+Q$. Let $\mathcal{V}$ be the holomorphic jet complex of $D$ as 
discussed above. Scalar multiplication by a quaternion 
$\lambda\colon V\to V$, which also is a complex linear isomorphism, 
preserves $D$ and hence lifts via Theorem~\ref{thm:uniqueness} to a unique isomorphism of
$V_n$ commuting with $\pi$ and $d$. Thus, $V_n$ becomes a quaternionic vector bundle
of rank $r(k+1)$ and
$\pi$ and $d$ become quaternionic linear. In particular, the kernels $\ker\pi_n=N_n\subset V_n$
and the flag 
\begin{equation}\label{eq:F_flag}
V_n\supset F_{n-1}\supset F_{n-2}\cdots\supset F_1\supset F_0\,,
\end{equation}
will be quaternionic with $\rank F_k=r(k+1)$. 

The canonical complex structure $S$ on the holomorphic jet complex introduced in 
Theorem~\ref{thm:mcS}
will be quaternionic linear: since  $\lambda S\lambda^{-1}$ also satisfies 
the properties in Theorem~\ref{thm:mcS},
by uniqueness we get $S=\lambda S\lambda^{-1}$. From \eqref{eq:D_n}, \eqref{eq:D_n_decomposed} we thus
get the quaternionic
holomorphic structures  
\begin{equation}\label{eq:quat_D_n}
D_n=\delbar_n+Q_n\colon \Gamma(V_n)\to \Gamma(\bar{K}V_n)
\end{equation}
on $V_n$ which induce the complex holomorphic structures $\delbar_n$ on the flag 
$F_k\subset V_n$ for $k< n$.
Since $\pi D_n= D_{n-1}\pi$, the maps $\pi$ are quaternionic holomorphic bundle maps.
The Hopf fields $Q_n$ are determined by $Q$ via \eqref{eq:Q_n}, i.e.,
\[
\pi^{n}Q_n=Q\pi^{n}\,.
\]
The prolongation maps $P_n\colon H^{0}(V)\to H^{0}(V_n)$
of Corollary~\ref{cor:prolongation} also become quaternionic linear. 
\begin{theorem}\label{thm:quat_holo_jet_complex}
Let $V$ be a quaternionic holomorphic vector bundle over a Riemann surface $M$. 
Then there exists a unique quaternionic
holomorphic jet complex $\mathcal{V}$. The jet bundles $V_n$ have canonical
complex structures $S_n$ given by Theorem~\ref{thm:mcS} and canonical
quaternionic holomorphic structures given by \eqref{eq:quat_D_n}. The  
projections $\pi$ are quaternionic holomorphic bundle maps and the prolongation
maps are quaternionic linear.
\end{theorem}

%%% Local Variables: 
%%% mode: latex
%%% TeX-master: "willmore"
%%% End: 

\section{Quaternionic Pl\"ucker formula}\label{sec:Pluecker}
We now come to the heart of this paper where we prove the quaternionic analogue
of the classical Pl\"ucker relations for a holomorphic curve, or more generally,
for a linear system. The classical Pl\"ucker formula
relates the basic integer invariants -- genus, degree and vanishing orders
of the higher derivatives -- of a holomorphic curve in
$\C\P^n$. In the quaternionic setting we already have seen that there is 
an additional non-trivial invariant, the  Willmore energy, which  also will enter into
the quaternionic Pl\"ucker relation. Thus, we will be able to use these relations to estimate
the Willmore energy of a holomorphic curve from below. 
Since the Willmore energy is zero for complex holomorphic
curves, we also recover the classical Pl\"ucker formula.    
\subsection{ Weierstrass points}
Let $V$ be a quaternionic holomorphic vector bundle over a Riemann surface $M$. We 
discuss the possible vanishing orders of a linear system
$H\subset H^{0}(V)$ at a given point $p\in M$.
Denote by
\[
n_0(p):=\min\{\ord_{p}\psi\,;\,\psi\in H\,\}
\]
the smallest vanishing order at $p$ among all holomorphic sections in $H$. 
Unless all sections in $H$ vanish to  precisely order $n_0(p)$, we define
\[
n_1(p):=\min\{\ord_{p}\psi>n_0(p)\,;\,\psi\in H\,\}
\]
and the linear subspace 
\[
H_1(p):=\{\psi\in H\,;\, \ord_p\psi\geq n_1(p)\,\}\subset H
\]
of sections of $H$ vanishing to at least order $n_1(p)$.
Continuing this procedure we arrive at a flag of linear subspaces
\[
H\supset H_1(p) \supset H_2(p) \supset \cdots \supset H_l(p)\supset 0
\]
and strictly increasing vanishing orders
\[
n_0(p)<n_1(p)<\cdots <n_l(p)
\]
satisfying
\[
H_k(p)=\{\psi\in H\,;\, \ord_p\psi\geq n_k(p)\,\}\,.
\]
If it is clear from the context, we will suppress the label $p$ and simply
write $n_k$ and $H_k$.
 
Note that the successive quotients $ H_k/H_{k+1}$ can have dimensions 
at least $\rank V$: fix $k$ and take the  $n_k$-th prolongation of the $H_j$
to $V_{n_k}$. Since the prolongation maps 
\eqref{eq:D_to_D_n} are injective we have $H_k/H_{k+1}=P_{n_k}H_k/P_{n_k}H_{k+1}$. 
By Lemma~\ref{lem:prolongation_order} the evaluation map
\[
ev_p\colon P_{n_k}H_k/P_{n_k}H_{k+1}\to \ker\pi_{n_k}
\]
is well defined and injective, and thus
\[
\dim H_k(p)/H_{k+1}(p)\leq \dim\ker\pi_{n_k}=\rank V\,.
\]
\begin{definition}\label{def:Weierstrass}
Let $V$ be a quaterionic holomorphic vector bundle and $H\subset H^{0}(V)$ a linear system.
We call
\[
H\supset H_1(p) \supset H_2(p) \supset \cdots \supset H_l(p)\supset 0
\]
the {\em Weierstrass flag} and 
\[
n_0(p)<n_1(p)< \cdots <n_l(p)
\]
the {\em Weierstrass gap sequence} of $H$ at $p\in M$. The successive quotients of the
Weierstrass flag satisfy
\[
\dim H_k(p)/H_{k+1}(p)\leq\rank V\,.
\]
\end{definition}

For a quaternionic holomorphic line bundle $L$ the situation simplifies considerably:
if $H\subset H^{0}(L)$ is an $n+1$-dimensional linear system the Weierstrass
flag at each $p\in M$ becomes the full flag 
\begin{equation}\label{eq:Weierstrass_flag}
H=H_n\supset H_{n-1}\supset\cdots\supset H_1\supset H_{0} 
\end{equation}
with $\dim H_k=k+1$. Note that we relabeled the elements of the flag 
according to their (projective) dimension. Generically one expects the Weierstrass gap sequence $n_k=k$.
The excess to the generic value gives an integer invariant of $H$:
\begin{definition}\label{def:order_H}
Let $L$ be a quaternionic holomorphic line bundle over a Riemann surface $M$
and $H\subset H^{0}(L)$ an $n+1$-dimensional linear system with
Weierstrass gap sequence $n_0<n_1< \cdots< n_{n}$. 
The {\em order at $p\in M$\,} of the linear system $H$  is defined by
\[
\ord_p H:=\sum_{k=0}^{n}(n_k(p)-k)=\sum_{k=0}^{n}n_k(p)-\tfrac{1}{2}n(n+1)\,.
\]
A point $p\in M$ is called a {\em Weierstrass point} (of the linear system $H$)
if $\ord_p H\neq 0$.

We will see below that Weierstrass points are isolated and thus there are only
finitely many of them on a compact Riemann surface $M$.
In this case we define the {\em order} of $H$ by 
\[
\ord H:=\sum_{p\in M}\ord_p H\,.
\]
The Weierstrass points of a holomorphic curve $f\colon M\to \H\P^n$ are the
Weierstrass points of the induced linear system $H\subset H^{0}(L^{-1})$
via the Kodaira correspondence (Theorem~\ref{thm:Kodaira}),
where $L=f^{*}\Sigma$ is the pull back of the tautological bundle. 
\end{definition}
For further applications it is important to understand the relationship
between Weierstrass points and prolongations of holomorphic sections.
Given the $n+1$-dimensional linear system $H\subset H^{0}(L)$ we define the
bundle map 
\begin{equation}\label{eq:ev_prolongation}
P\colon {H}\to L_{n}\qquad\text{by}\qquad P_p:=ev_p\circ P_{n}\,,
\end{equation}
where $P_{n}\colon H^{0}(L)\to H^{0}(L_{n})$ is the $n$-th prolongation map 
\eqref{eq:D_to_D_n}
of the holomorphic jet complex $\mathcal{L}$ of $L$.
This map will play a crucial role in our discussion of the Pl\"ucker relations
later on. 
\begin{lemma}\label{lem:weierstrass_char}
A point $p\in M$ is a Weierstrass point precisely when $P_p$ fails to be an isomorphism.
Away from Weierstrass points, $P$
maps the Weierstrass flag $H_k$ in $H$ isomorphically to the flag $F_k=\ker \pi^{k+1}$
in $L_{n}$ given by \eqref{eq:F_flag}. Therefore, away from Weierstrass points
the Weierstrass flag $H_k\subset H$ is a smooth flag.

In particular, if the linear system $H$ has no Weierstrass points, $P\colon {H}\to L_n$
is a bundle isomorphism mapping the Weierstrass flag $H_k$ to the flag $F_k$. 
\end{lemma}
\begin{proof}
If $p\in M$ is a Weierstrass point then there is a non-trivial section $\psi\in H$ with
$\ord_p\psi\geq n+1$.  Lemma~\ref{lem:prolongation_order} then implies that
$P_p(\psi)=0$, i.e., $P_p$ is not an isomorphism. On the other hand, if $p\in M$
is not a Weierstrass point then $P_p\colon H\to (L_n)_p$ is clearly injective: $P_p(\psi)=0$
implies that $\psi$ vanishes to at least order $n+1$ at $p$, and thus has to be identically
zero. From Corollary~\ref{cor:prolongation} we get $\pi^{k}P_n=P_{n-k}$, so that
$\pi^{k+1} P_p H_k=0$ and thus $P_p H_k\subset (F_k)_p$. Since $\dim H_k=\rank F_k$
the claim follows. 
\end{proof}
\subsection{Frenet curves}
The Kodaira correspondence assigns to a base point free linear system a holomorphic
curve. We will now discuss curves arising from linear systems without Weierstrass
points in more detail. They demonstrate a number of
conceptual features without the technical difficulties arising
from the existence of Weierstrass points. Whereas Weierstrass points are
generic in the complex case -- the only curve without Weierstrass points is the
rational normal curve -- they are in a certain sense ``non-generic'' in the quaternionic
setting.

Let $V$ be a flat quaternionic vector bundle of rank $n+1$ and $L\subset V$ a
line subbundle. Recall that by Definition~\ref{def:holo_curve} 
$L$ is a holomorphic curve if 
there exists a complex structure $J$ on $L$ such that 
$*\delta=\delta\, J$. Here
$\delta=\pi\nabla\in\Gamma(\text{Hom}(L,V/L)K)$ is the derivative of $L\subset V$,
where $\pi\colon V\to V/L$ is the canonical projection and $\nabla$ the flat connection on $V$.
In this situation we already have seen in Theorem~\ref{thm:curve_holo_structure}
that the dual bundle
$L^{-1}$ has a canonical holomorphic structure.  
\begin{definition}\label{def:Frenet_curve}
Let $L\subset V$ be a holomorphic curve. A {\em Frenet flag} for $L$ is given by the following
data:
\begin{enumerate}
\item
A full flag
\[
V_0=L\subset V_1\subset\cdots \subset V_{n-1}\subset V_{n}=V
\]
of quaternionic subbundles of $\rank V_k=k+1$ starting with $L$, such that
\[
\nabla\,\Gamma(V_k)\subset \Omega^{1}(V_{k+1})\,.
\]
In this case the derivatives of $V_k\subset V_{k+1}$ are given by
\[
\delta_k:=\pi_k\nabla\colon V_{k}/V_{k-1}\to T^{*}M\otimes V_{k+1}/V_{k}
\]
where $\pi_k\colon V\to V/V_k$ is the $k$-th quotient projection.
\item
Complex structures $J_k$ on the successive quotient line bundles  $V_{k}/V_{k-1}$, $J_0=J$,
such that the derivatives $\delta_k$ satisfy 
\[
*\delta_k=J_{k+1}\delta_k=\delta_k J_k\,,
\]
which is to say that $\delta_k\in \Gamma(K\text{Hom}_{+}(V_{k}/V_{k-1},V_{k+1}/V_{k}))$.
Additionally we demand each  
\[
\delta_k\colon V_{k}/V_{k-1}\to K\,V_{k+1}/V_{k}
\]
to be an isomorphism of line bundles. 
\end{enumerate}
A holomorphic curve $L\subset V$ is called a {\em Frenet curve} if $L$ has a 
Frenet flag (which then is necessarily unique). 
\end{definition}
The elements $V_k\subset V$ of the Frenet flag are the $k$-th osculating curves
to the curve $L$. We will see below that they too are holomorphic curves.
For now we content ourselves to observe that the derivative of the highest osculating curve 
$V_{n-1}\subset V$ satisfies
\[
*\delta_{n-1}=J_{n}\delta_{n-1}\,,
\]
where we interpret $\delta_{n-1}$ as a section of $K\text{Hom}(V_{n-1},V/V_{n-1})$ by
composing with the projection $V_{n-1}\to V_{n-1}/V_{n-2}$. 
Lemma~\ref{lem:left_holo} then implies that 
$V_{n-1}^{\perp}\subset V^{-1}$ is also a holomorphic curve. 
\begin{definition}\label{def:dual_curve}
Let $L\subset V$ be a Frenet curve with Frenet flag $V_k$. Analogous to the complex case
we call 
$V_{n-1}^{\perp}\subset V^{-1}$ the {\em dual curve} to $L$ and denote it by $L^{*}$.
\end{definition}
Of course, $L^{*}\subset V^{-1}$ is again a Frenet curve with Frenet flag $V_k^{\perp}$ and
$(L^{*})^{*}=L$. 
\begin{example}\label{ex:frenet_curve}
The simplest examples of Frenet curves are immersed holomorphic curves 
in $\H\P^1$, i.e., conformal immersions into the 4-sphere \cite{bflpp}.
If $L\subset V$ is immersed and $V$ has rank $2$, we use the derivative
$\delta$ to transport the complex structure $J$ to $V/L$: take any non-zero
tangent vector $X\in T_p M$, then $\delta_{X}$ is non-zero, and define $J_1$ on $V/L$ at $p$ by
\[
J_1\delta_{X}:=\delta_{X}J\,.
\]
Due to \eqref{eq:right_holo} this gives $J_1$ independently of the choice of non-zero $X$.  
By construction we then have
\[
*\delta=J_1\delta=\delta J
\]
and $L\subset V$ is a Frenet flag. In this case the dual curve $L^{*}$ is just $L^{\perp}$.  
\end{example}
We now show that the Weierstrass flag of a linear system without Weierstrass points
is a Frenet flag.
\begin{theorem}\label{thm:Kodaira_Frenet}
Let $L^{-1}$ be a holomorphic line bundle with an $n+1$-dimensional linear
system $H\subset H^{0}(L^{-1})$ without Weierstrass points and let $L\subset H^{-1}$
be the Kodaira embedding of $L$. Then the 
Weierstrass flag $H_k\subset H$ is the Frenet flag of the dual curve $L^{*}=H_0$.
In particular, the Kodaira embedding $L\subset H^{-1}$ is a Frenet curve with Frenet flag
$H_k^{\perp}$.
\end{theorem} 
\begin{proof}
Let $\tilde{\mathcal{L}}$ denote the jet complex of the holomorphic bundle $L^{-1}$.
Since $H$ has no Weierstrass points the bundle map
\[
P\colon {H}\to \tL_n\,,\qquad P_p=ev_p\circ P_n
\]
is an isomorphism by Lemma~\ref{lem:weierstrass_char} and maps the Weierstrass flag
$H_k$ to the flag $F_k=\ker\pi^{k+1}$ in $\tL_n$. It thus suffices
to show that $F_k$ is a Frenet flag. Put $V^{-1}:=\tL_n$ which 
has the trivial connection $\nabla$ coming from the directional 
derivative in $H$ via $P$. Since $dP_n=0$, but also $\nabla P_n=0$, we obtain
$\pi\nabla=d$, so that $\nabla$ is an adapted connection. But then 
\[
\pi^{k+2}\nabla=\pi^{k+1}d\alpha=d\pi^{k+1}
\]
which implies 
\[
\nabla \Gamma(F_k)\subset \Omega^{1}(F_{k+1})\,.
\]
The successive quotients $ F_k/F_{k-1}$ get mapped by 
$\pi^{k}$ isomorphically to the kernels $N_{n-k}\subset \tL_{n-k}$
and therefore have a complex structure. Now let 
$\delta_k\in\Omega^{1}(\text{Hom}(F_k/F_{k-1},F_{k+1}/F_k))$ be the derivatives 
of the flag $F_k$ and let $\tilde{\delta}$ be
the restrictions of $d$ to the kernels $\ker\pi$ in the jet complex $\tilde{\mathcal{L}}$.
Then, on $F_k/F_{k-1}$, we have 
\[
\pi^{k+1}\delta_k=\pi^{k+1}\nabla=\pi^{k}d=d\pi^{k}=\tilde{\delta}_{n-k}\pi^{k}\,
\]
which implies that $\delta_k\colon F_k/F_{k-1}\to K\,F_{k+1}/F_k$ is an isomorphism
of line bundles. Thus $F_k$ is a Frenet flag by Definition~\ref{def:Frenet_curve}. 

Finally, $H_0\subset H$ is the dual curve to the Kodaira embedding $L\subset H^{-1}$
if and only if $L=H_{n-1}^{\perp}\subset H^{-1}$. But this follows 
immediately since the kernel of the
evaluation map
\[
ev_p\colon H\to L^{-1}_p
\]
is $H_{n-1}(p)$ and the Kodaira correspondence is given by $L_p= (\ker ev_p)^{\perp}$. 
\end{proof}
Note that in the situation of the above theorem the dual isomorphism
\[
P^{*}\colon \tL_n^{-1}\to H^{-1}
\]
maps the dual flag $F_k^{\perp}$ to the Frenet flag $H_k^{\perp}$ of the Kodaira embedding
$L\subset H^{-1}$. Since $F_k^{\perp}=\tL_{n-k-1}^{-1}$ via the inclusion 
$(\pi^{k+1})^{*}\colon \tL_{n-k-1}^{-1}\to \tL_n^{-1}$  we see that the Frenet flag of
the curve $L\subset H^{-1}$ consists of the first $n$ jet bundles of the holomorphic
jet complex $\tilde{\mathcal{L}}$ of $L^{-1}$. The next theorem shows that 
the converse also holds.
\begin{theorem}\label{thm:Frenet_jet_complex}
Let $V_k$ be the Frenet flag of the Frenet curve $L\subset V$. Then $\tL_k:=V_k^{-1}$ are the first
$n$ jet bundles of the holomorphic jet complex $\tilde{\mathcal{L}}$ of the holomorphic 
line bundle $L^{-1}$. 

The linear system $H\subset H^{0}(L^{-1})$, induced by the Kodaira correspondence,
is given by $H=V^{-1}$ and the bundle map $P\colon H\to\tL_n$ introduced in
\eqref{eq:ev_prolongation} is the identity map. 
In particular, Frenet curves have no Weierstrass points.
\end{theorem}
\begin{proof}
We first verify the axioms of the holomorphic jet complex:
dualizing the Frenet flag we obtain the sequence of surjective bundle homomorphisms
\[
\tL_n\overset{\pi_n}{\to} \tL_{n-1}\overset{\pi_{n-1}}{\to}\cdots 
\overset{\pi_2}{\to}\tL_1\overset{\pi_1}{\to}\tL_0=L^{-1}
\]
whose kernels are
\[
N_k=\ker\pi_k=(V_k/V_{k-1})^{-1}\,.
\]
Thus we have complex structures on each $N_k$ by dualizing the complex structures on $V_k/V_{k-1}$.
Next we define 
$
d_k\colon \Gamma(\tL_k)\to \Omega^{1}(\tL_{k-1})
$
by the obvious product rule
\begin{equation}\label{eq:def_d}
<d_k\alpha,\psi>:=d<\pi_k\alpha,\psi>-<\alpha,\nabla\psi>\,,
\end{equation}
where $\alpha\in \Gamma(\tL_k)$, $\psi\in\Gamma(V_{k-1})$ and $\nabla$ 
is the flat connection on $V=V_n$. Since $\nabla\colon \Gamma(V_k)\to \Omega^{1}(V_{k+1})$
this is well defined. It is routine to check the required Leibniz rule for
$d_k$. Furthermore, $d_{k-1}d_k=0$ since $\nabla$ is flat. To calculate the restriction
$\tilde{\delta}_k$
of $d_k$ to $N_k$ we take $\alpha\in \Gamma(N_k)$, i.e., $\alpha$ is a section of
$V_k^{-1}$ vanishing on $V_{k-1}$, and $\psi\in \Gamma(V_{k-1})$. Then
\[
<\tilde{\delta}_k\alpha,\psi>=<d_k\alpha,\psi>=-<\alpha,\nabla\psi>=-<\alpha,\delta_{k-1}\psi>
\]
so that 
\[
\tilde{\delta}_k=-\delta_{k-1}^{*}\,.
\]
The properties of $\delta_k$ now imply that
$
\tilde{\delta}_k\colon N_k\to K N_{k-1}
$
are complex linear isomorphisms. To verify the last axiom $D\pi_1=d_1''$, where $D$ denotes
the holomorphic structure on $L^{-1}$, we take $\alpha\in\Gamma(\tL_1)$ and $\psi\in\Gamma(L)$
to calculate
\begin{align*}
<d_1''\alpha,\psi>&=\tfrac{1}{2}<d_1\alpha+*Jd_1\alpha,\psi>\\
&=\tfrac{1}{2}(d<\pi_1\alpha,\psi>+*d<\pi_1\alpha,J\psi>)-
\tfrac{1}{2}(<\alpha,\nabla\psi+*\nabla J\psi>)\\
&= \tfrac{1}{2}(d<\pi_1\alpha,\psi>+*d<\pi_1\alpha,J\psi>)-
\tfrac{1}{2}(<\pi_1\alpha,\nabla\psi+*\nabla J\psi>)\\
&=<D\pi_1\alpha,\psi>\,.
\end{align*}
We used \eqref{eq:holo_curve_mixed} to see that 
$\nabla\psi+*\nabla J\psi$ is again a section of $L$
and the definition \eqref{eq:holo_structure_induced} of the holomorphic structure on $L^{-1}$.
But the jet complex up to any level is determined uniquely by the axioms 
in section~\ref{subsec:holo_jet_complex}.

By Definition~\ref{def:Frenet_curve} any Frenet curve $L\subset V$ is full.
The linear system $H\subset H^{0}(L^{-1})$ induced by the Kodaira correspondence
(Theorem~\ref{thm:Kodaira}) is therefore equal to $H=V^{-1}$. If $\alpha\in H$ then 
\[
\pi^{n}\alpha=\alpha_{|L}\qquad\text{and}\qquad d_n\alpha=0\,,
\]
where the latter follows from \eqref{eq:def_d}. Therefore
the $n$-th prolongation of $\alpha\in H$ is given by $P_n(\alpha)=\alpha$, i.e.,
$P\colon H\to \tL_n$ is the identity map. Lemma~\ref{lem:weierstrass_char} then implies that
the curve $L$, or equivalently, the linear system $H$, has no Weierstrass points.
\end{proof}
Putting the last two theorems together we get bijective correspondences between
the following objects:
\begin{enumerate}
\item
Frenet curves in $\H\P^n$ up to projective equivalence,
\item
$n+1$-dimensional linear systems without Weierstrass points, and
\item
Holomorphic jet complexes of holomorphic line bundles (up to isomorphisms) whose $n$-th jet bundles
admit a trivial, adapted connection.
\end{enumerate}
To see that the osculating curves $V_k$ of a Frenet curve $L\subset V$ also are holomorphic,
we need to put a complex structure on $V$ which stabilizes the Frenet flag and induces the
given complex structures on the quotients $V_k/V_{k-1}$. By the above correspondence such
a complex structure always comes from an adapted complex structure on the 
holomorphic jet complex of $L^{-1}$. As we have seen in Theorem~\ref{thm:mcS}
there is a canonical complex structure on the holomorphic jet complex. 
\begin{theorem}\label{thm:Pluecker_unramified}
Let $L\subset V$ be a Frenet curve with Frenet flag $V_k$.
Then there exists a unique complex structure $S$ on $V$
such that
\begin{enumerate}
\item
$S$ stabilizes the flag $V_k$ and induces the given complex structures 
$J_k$ on the quotients $V_k/V_{k-1}$.
\item
If $\nabla=\del+A+\delbar+Q$ is the type decomposition \eqref{eq:nabla_decomposition}
with respect to $S$ then
\[
Q_{|V_{n-1}}=0\qquad\text{or, equivalently}\qquad A(V)\subset L\,.
\]
\end{enumerate}
The elements $V_k$ of the Frenet flag  are holomorphic, both as 
curves and as subbundles, and so is the Frenet flag for the dual curve
$L^{*}$. Moreover, if $M$ is compact of genus $g$
we have the {\em unramified} Pl\"ucker relation
\begin{equation}\label{eq:Pluecker_unramified}
W(L^{-1})-W((L^{*})^{-1})=4\pi (n+1)(n(1-g)-d)
\end{equation}
where $d=\deg L^{-1}$ denotes the degree of the holomorphic curve $L$.
\end{theorem}
In the case when $L$ is a complex holomorphic curve  also $L^{*}$ is complex holomorphic. 
Therefore their Willmore energies vanish 
and we recover the unramified Pl\"ucker relation
$n(1-g)=d$. These imply $g=0$ so that $L$ is a degree $n$ rational curve in $\C\P^n$, 
the {\em rational normal curve}.
\begin{proof}
We have seen in Theorem~\ref{thm:Frenet_jet_complex} that the dualized Frenet flag
$V_k^{-1}$ consists of the first $n$ jet bundles of the holomorphic jet complex 
$\tilde{\mathcal{L}}$ of $L^{-1}$. Theorem~\ref{thm:mcS} provides us with a canonical
adapted complex structure $S$ on $\tilde{\mathcal{L}}$. Therefore the dual complex
structure $S_n^{*}$ on $V$, which we again call $S$, leaves the Frenet flag invariant and
induces $J_k$ on the quotients $V_k/V_{k-1}$. 
On $V$ respectively $V^{-1}$ we can now 
decompose \eqref{eq:nabla_decomposition} the connections 
\[
\nabla=\del+A+\delbar+Q\qquad\text{and}\qquad \nabla^{*}=\del- Q^{*}+\delbar-A^{*}
\]
with respect to $S$. Due to \eqref{eq:def_d} the dual connection $\nabla$ on $\tL_n=V^{-1}$ 
satisfies $\pi_n\nabla=d_n$ and hence is adapted \eqref{eq:jet_nabla}. 
Corollary~\ref{cor:adapted_nabla_relations} then implies
\[
\pi Q^{*}=0\qquad\text{and}\qquad A^{*}_{|L^{\perp}}=0
\]
which is equivalent to 
\[
Q_{|V_{n-1}}=0 \qquad\text{and}\qquad A(V)\subset L\,.
\]
Since
\[
*\delta_k=S\delta_k=\delta_k S
\]
Corollary~\ref{cor:holo_compare} implies that $V_k\subset V$ and
$V_k^{\perp}\subset V^{-1}$ are holomorphic,
both as subbundles and curves. 

To see that $S$ is unique we note that the conditions of the theorem imply that 
the dual complex structures $S_k$ on $\tL_k=V_k^{-1}$ satisfy the requirements
of Theorem~\ref{thm:mcS} for the first $n$ jet bundles. Thus, $S_k$ are the unique
adapted complex structures on the holomorphic jet complex of $\tilde{\mathcal L}$.

To verify the unramified Pl\"ucker relations we integrate the trace of the
curvature 
\[
R^{\nabla}=R^{\del+\delbar}+d^{\del+\delbar}(A+Q)+(A+Q)\wedge (A+Q)\,
\]
of the flat connection $\nabla$ on $V$.
Except $d^{\del+\delbar}(A+Q)$, which anticommutes with $S$, all other terms are
$S$ commuting so that 
\[
0=R^{\nabla}_{+}=R^{\del+\delbar}+A\wedge A +Q\wedge Q\,,
\]
where we used that $A\wedge Q=Q\wedge A=0$ by type. Multiplying this last relation by $S$ 
we obtain 
\[
SR^{\del+\delbar}=A\wedge *A-Q\wedge *Q\,
\]
and, taking traces, also 
\[
<SR^{\del+\delbar}>=<A\wedge *A>-<Q\wedge *Q>\,.
\]
Considering that $Q_{|V_{n-1}}=0$ and $A(V)\subset L$, we obtain
\[
<SR^{\del+\delbar}>=<A\wedge *A\,_{|L}>-<Q\wedge *Q\,_{V/V_{n-1}}>
\]
which, together with  
Corollary~\ref{cor:holo_compare} and \eqref{eq:chern},
implies
\[
4\pi\deg V=W(L^{-1})-W((L^{*})^{-1})\,.
\]
Since $V^{-1}=\tL_n$ we obtain from \eqref{eq:jet_degree}
\[
\deg V=-\deg \tL_n=(n+1)(n(1-g)-d)
\]
which finishes the proof. 
\end{proof}
Before we continue we give a more geometric interpretation
of the above result by explaining its connection to the 
{\em  mean curvature sphere} of a conformally immersed surface. 
For this assume $V={\H}^{n+1}$ to be the trivial bundle so that
$L\subset V$ is a Frenet curve in $\H\P^n$. The holomorphic curves $V_k$ then are the 
$k$-th osculating curves in $G_{k+1}(\H^{n+1})$ of $L$. The canonical complex structure
$S$ on $V$, given by the above theorem, describes a ``congruence'' of osculating $\C\P^{n}$'s
along the curve $L$.  

To see this we regard $\H^{n+1}=\C^{2n+2}$ via right multiplication by the
quaternion $i$.
Let $W\subset \H^{n+1}$ be a complex linear subspaces satisfying $W\oplus Wj=\H^{n+1}$.
Then the map
\[
W\supset w\C\mapsto w\H\subset \H^{n+1}
\]
is injective and defines a $\C\P^n\subset \H\P^n$. Note that $W$ and $Wj$ give rise to the same
$\C\P^n$. On the other hand, splittings $\H^{n+1}=W\oplus Wj$ are the same as
linear endomorphisms $S$ on $\H^{n+1}$ with $S^2=-1$ by declaring $W$ to be the 
$i$-eigenspace of $S$. Therefore, any complex structure $S$ on $\H^{n+1}$ gives rise
to the $\C\P^n\subset \H\P^n$ consisting of all quaternionic lines fixed by $S$.

To see that the congruence of $\C\P^n$'s given in Theorem~\ref{thm:Pluecker_unramified}
is indeed osculating the curve $L$ we note that, since $S(L)\subset L$, the 
$\C\P^n$'s pass through the curve $L$ at each point. But we also
know that the derivative $\delta$ of $L$ is a section of 
$K\text{Hom}_{+}(L,{\H}^{n+1}/L)$. Since the tangent space to $\C\P^n\subset \H\P^n$
at the point $L_p$ is given by $\text{Hom}_{+}(L,{\H}^{n+1}/L)_p$, we see that our
curve $L$ is tangent to the congruence of $\C\P^n$'s. Of course, all of this would also hold
for any adapted complex structure $S$. The property $Q_{|V_{n-1}}=0$, picking out a unique
osculating congruence of $\C\P^n$'s, is an $n+1$-st order tangency condition on the curve $L$, which 
is easiest to interpret in the case of curves into $\H\P^1$. 
\begin{example}\label{ex:4.2}
In our standard example of a conformal immersion $f\colon M\to\H\P^1$, which is a Frenet curve,
the congruence of osculating $\C\P^1$'s is the classical mean curvature sphere congruence
\cite{bflpp}, or {\em conformal Gauss map}, of the immersion $f$. The mean curvature 
sphere at a given point $p\in M$ is the unique $2$-sphere in $S^4=\H\P^1$ touching $f$
at $f(p)$ and having the same mean curvature vector as $f$ at this point.
Notice that to formulate this last -- second order -- condition we had to break the M\"obius
symmetry and introduce euclidean quantities. The formulation given in 
Theorem~\ref{thm:Pluecker_unramified} is intrinsically M\"obius invariant.
\end{example}
Given a Frenet curve $L$ in $\H\P^n$ the unramified Pl\"ucker relation gives an estimate
of its Willmore energy from below
\begin{equation}\label{eq:Willmore_estimate}
W(L^{-1})\geq 4\pi (n+1)(n(1-g)-d)
\end{equation}
in terms of dimension, genus, and degree.
{F}rom Theorem~\ref{thm:Pluecker_unramified} we see that
equality holds if and only if the dual curve $L^{*}$ has zero Willmore energy. 
To characterize Frenet curves for which \eqref{eq:Willmore_estimate} becomes an
equality it thus suffices to describe Frenet curves with zero Willmore energy.
\begin{lemma}\label{lem:Frenet_curve_0}
Let $L$ be a Frenet curve in $\H\P^n$ with $W(L^{-1})=0$. Then $L$ is the 
twistor projection \eqref{eq:twistor_proj} of a complex holomorphic curve
$E$ in $\C\P^{2n+1}$
whose Weierstrass gap sequence $n_k=k$ for $k=0,\dots, n$ and whose
$n$-th osculating curve $W_n\subset \H^{n+1}=\C^{2n+2}$ satisfies $W_n\oplus W_{n}j=\H^{n+1}$.   
\end{lemma}
\begin{proof}
Let $V_k\subset \H^{n+1}$ be the Frenet flag of $L$ and let $S$ be the complex structure
on $\H^{n+1}$ given in Theorem~\ref{thm:Pluecker_unramified}. As usual we consider
$\H^{n+1}=\C^{2n+2}$ via right multiplication by the quaternion $i$. Since $S$ stabilizes
the flag $V_k$ we have
\[
V_k=W_k\oplus W_k j
\]
with $W_k$ the $i$-eigenspace of $S$ on $V_k$. From Lemma~\ref{lem:twistor_char} we know that
$L$ has zero Willmore energy if and only if the twistor lift $E=W_0$ is a complex holomorphic
curve in $\C\P^{2n+1}$. The derivatives of the Frenet flag
\[
\delta_k\colon V_k/V_{k-1}\to K V_{k+1}/V_k
\]
commute with $S$ and thus induce 
\[
\delta_k\colon W_k/W_{k-1}\to K W_{k+1}/W_k\,,
\]
which are easily seen to be the derivatives of the first $n$ osculating curves $W_k$ of the
complex holomorphic curve $E$. Since $\delta_k$ are isomorphisms the curve $E$ has the
Weierstrass gap sequence $n_k$ with $n_k=k$  for $k\leq n$.    
\end{proof}
Combining this lemma with the remarks above we obtain the following characterization
of Frenet curves for which equality holds in \eqref{eq:Willmore_estimate}. 
\begin{corollary}\label{cor:Willmore_equality}
Let $L$ be a Frenet curve in $\H\P^n$ of degree $d$ and genus $g$.
Then 
\[
W(L^{-1})= 4\pi (n+1)(n(1-g)-d)
\]
if and only if the dual
curve $L^{*}$ is the twistor projection of a holomorphic curve $E$ in $\C\P^{2n+1}$
whose Weierstrass gap sequence $n_k=k$ for $k=0,\dots, n$ and whose
$n$-th osculating curve $W_n\subset \H^{n+1}=\C^{2n+2}$ satisfies $W_n\oplus W_n j=\H^{n+1}$.   
\end{corollary}
\begin{proof}
From Theorem~\ref{thm:Pluecker_unramified} we see that equality holds in 
\eqref{eq:Willmore_estimate} if and only if the dual curve $L^{*}$ has zero Willmore energy.
By Lemma~\ref{lem:twistor_char} this is equivalent to $L^{*}=\pi(E)$ being the twistor
projection of a holomorphic curve in $\C\P^{2n+1}$.
\end{proof} 

\subsection{The general Pl\"ucker relation}
In Theorem~\ref{thm:Pluecker_unramified} we gave the quaternionic version of the unramified
Pl\"ucker relation, that is, the Pl\"ucker formula for a linear system $H\subset H^{0}(L)$
without
Weierstrass points. In this situation we conclude from Lemma~\ref{lem:weierstrass_char}
that the bundle map 
\[
P\colon H\to L_n\,,\qquad P_p=ev_p\circ P_n\,,
\]
where $L_n$ is the $n$-th jet bundle of the holomorphic jet complex $\mathcal{L}$
of $L$, is an isomorphism. Moreover, by Theorem~\ref{thm:Kodaira_Frenet} the flag 
$F_k=\ker\pi^{k+1}\subset L_n$ 
is the Frenet flag of the dual curve $(L^{*})^{-1}=F_0$ of the Kodaira embedding 
$L^{-1}\subset L_n^{-1}\cong H^{-1}$.
Here we view $L_n$ as a trivial bundle via the trivial connection $\nabla$ induced
by pushing forward the directional derivative in $H$ by $P$. 
By construction $\nabla P=0$ and, since $d P=0$, we conclude
$\pi\nabla=d$, so that $\nabla$ is a flat adapted connection on ${L_n}$. Let
\[
\nabla=\hat{\nabla}+A+Q
\]
be the decomposition into $S$-commuting and anticommuting parts \eqref{eq:nabla_decomposition}
with respect to the canonical complex structure $S$ on $\mathcal{L}$ given
in Theorem~\ref{thm:mcS}. Then 
\[
W(L)=2\int_{M}<Q\wedge *Q>\,,\qquad W(L^{*})=2\int_{M}<A\wedge *A>\,,
\]
and the unramified Pl\"ucker relation for the linear system $H\subset H^{0}(L)$ is given by
\begin{equation}\label{eq:Pluecker_unramified_1}
W(L)-W(L^{*})=4\pi (n+1)(n(1-g)-d)\,,
\end{equation}
where $g$ is the genus of $M$ and $d=\deg L$ is the degree of the holomorphic line bundle $L$.
Note that this formula is identical to the one given in Theorem~\ref{thm:Pluecker_unramified}
after exchanging the roles of $L$ and $L^{-1}$ .

In the general case, when the linear system $H$ has Weierstrass points, the bundle map
$
P\colon H\to L_n
$
fails to be an isomorphism precisely at the Weierstrass points by Lemma~\ref{lem:weierstrass_char}.
As we shall see below in 
Lemma~\ref{lem:P-matrix}, there are only finitely many such points 
$p_{\alpha}\in M$, $\alpha=1,\dots,N$, for the given linear system $H$.
Denote by $M_0=M\setminus\{p_1,\dots,p_N\}$ the open
Riemann surface punctured at the Weierstrass points. Over $M_0$ the bundle map 
$P$ is an isomorphism and, as in the unramified case above, we push forward the 
trivial connection on $H$
to obtain the trivial, adapted connection $\nabla$ on ${L_n}_{|M_0}$. We can decompose
\begin{equation}\label{eq:nabla_m0}
\nabla=\hat{\nabla}+A+Q
\end{equation}
where $\hat{\nabla}=\del+\delbar$ is a complex connection, i.e., 
$\hat{\nabla}S=0$, and $\delbar+Q$ extends to all of $M$ since, by Lemma~\ref{lem:nabla_holo},
it is the holomorphic structure $D_n$ on $L_n$. 
The flag $F_k\subset L_n$ is the Frenet flag of the dual curve $F_0=(L^*)^{-1}$ over $M_0$,
and $2<A\wedge *A>$ is the Willmore integrand for the holomorphic bundle $L^{*}\to M_0$
with holomorphic structure $\delbar-A^{*}$. Note that this holomorphic structure does not
extend to $M$, since $A$ is singular at the Weierstrass points but, as we will
see below, its Willmore
energy $2\int_{M}<A\wedge *A>$ remains finite.  
\begin{theorem}\label{thm:Pluecker}
Let $L$ be a quaternionic holomorphic line bundle of degree $d$ 
over a compact Riemann surface $M$ of genus $g$ and $H\subset H^{0}(L)$ an 
$n+1$-dimensional linear system. Then the Willmore energy of the dual curve  
\begin{equation}\label{eq:W_dual_curve}
W(L^{*}):=2\int_{M}<A\wedge *A>\;\,<\,\infty\,,
\end{equation}
is finite and we have the Pl\"ucker formula
\begin{equation}\label{eq:Pluecker}
\tfrac{1}{4\pi}(W(L)-W(L^{*}))=(n+1)(n(1-g)-d)+\ord H\,.
\end{equation} 
Here the non-negative 
integer $\ord H$ denotes the order of the linear system $H$ given in Definition~\ref{def:order_H}.
\end{theorem}
This theorem immediately implies the unramified Pl\"ucker formula
\eqref{eq:Pluecker_unramified_1} when $\ord H=0$. Since $W(L^{*})\geq 0$ we obtain a generalization
of the estimate \eqref{eq:Willmore_estimate}:
\begin{corollary}\label{cor:Pluecker_estimate}
Let $L$ be a quaternionic holomorphic line bundle of degree $d$ 
over a compact Riemann surface $M$ of genus $g$ and $H\subset H^{0}(L)$ an 
$n+1$-dimensional linear system. Then the Willmore energy is bounded from below by
\begin{equation}\label{eq:Pluecker_estimate}
\tfrac{1}{4\pi}W(L)\geq(n+1)(n(1-g)-d)+\ord H\,.
\end{equation}
\end{corollary}
Note that the presence of Weierstrass points increases the Willmore energy. This 
important feature will be used in the applications to eigenvalue estimates in the 
next section. 

Of course, if $L$ is a {\em complex} holomorphic line bundle then so is $L^{*}$ and, since
$W(L)=W(L^{*})=0$, we recover
the classical Pl\"ucker relation 
\[
\ord H= (n+1)(n(g-1)+d)
\]
for an $n+1$-dimensional linear system $H\subset H^{0}(L)$.

We now will work towards a proof of the above theorem. Since some of the
calculations and constructions have a technical flavor it will help to begin with
an outline of the basic ideas: the $S$-commuting part of the curvature 
of the flat connection $\nabla=\hat{\nabla}+A+Q$ in \eqref{eq:nabla_m0} is given, 
as in the proof of Theorem~\ref{thm:Pluecker_unramified}, by
\[
0=\hat{R}+A\wedge A+Q\wedge Q\,.
\]
Multiplying 
by $S$ and taking the trace we obtain 
\begin{equation}\label{eq:unintegrated_m0}
<S\hat{R}>+<Q\wedge *Q>-<A\wedge *A>=0
\end{equation}
over $M_0$. From Lemma~\ref{lem:P-matrix} below it will follow that
the complex connection $\hat{\nabla}$ has only logarithmic singularities
at the Weierstrass points $p_{\alpha}$ contributing an additional $-2\pi\ord_{p_{\alpha}} H$ to  
its curvature integral, so that
\begin{equation}\label{eq:R_hat_integral}
\int_{M}<S\hat{R}>=2\pi(\deg L_n-\ord H)\,.
\end{equation}
Since $W(L)=2\int_{M}<Q\wedge *Q>$ we conclude from \eqref{eq:unintegrated_m0}
that $\int_{M}<A\wedge *A>$ is finite. Finally we
integrate \eqref{eq:unintegrated_m0} over $M$,
insert \eqref{eq:R_hat_integral}, and recall \eqref{eq:jet_degree} to obtain
the general Pl\"ucker relation
\[
\tfrac{1}{4\pi}(W(L)-W(L^{*}))=(n+1)(n(1-g)-d)+\ord H\,.
\]

After this prelude we are going to fill in some more details. The crucial
technical lemma concerns the behavior of the smooth bundle map 
$
P\colon H\to L_n
$
across the Weierstrass points. We analyze this behavior by
calculating a local matrix expression of $P$ with respect to a suitable local frame in
$L_n$: let $\psi_k$ be a basis of $H$ whose vanishing
orders at some fixed point $p_0\in M$ are given by the Weierstrass gap sequence $n_k$ of 
$H$ at $p_0$. Then, by Lemma~\ref{lem:zeros},
\[
\psi_k=z^{n_k}\phi_k+O(n_k+1)
\]
with $\phi_k(p_0)\neq 0$ and $z$ a centered, holomorphic coordinate near $p_0$. 
Scaling each $\psi_k$ by a constant we may assume that
$\phi_k(p_0)=\phi_{p_0}\in L_{p_0}$ and that $S\phi_{p_0}=\phi_{p_0}i$.  
Extending $\phi_{p_0}$ to a local nowhere vanishing section $\phi$ of $L$ with $S\phi=\phi i$,
we obtain for all $k=0,1,\dots,n$
\[
\psi_k=z^{n_k}\phi+z^{n_k}(\phi_k-\phi)+O(n_k+1)=z^{n_k}\phi+O(n_k+1)\,.
\]
We apply Lemma~\ref{lem:prolongation_order} to calculate
the $l$-th prolongation
\begin{equation}
\label{eq:pre_basis}
P_{l}\psi_k=z^{n_k-l}(n_k(n_k-1)\cdots(n_k-l+1)\phi_l+O(1))\,,
\end{equation}
where $\phi_l$ are local nowhere vanishing sections of $N_l=\ker\pi_l\subset L_l$ satisfying
$\delta^l\phi_l=(-1)^{l}dz^l\phi$. Since $S\phi=\phi i$ and $\delta$ commutes with $S$ we also
have $S\phi_l=\phi_l i$. We now define a local frame $\tilde{\psi}_k$ of $L_n$ by 
requiring
\begin{equation}\label{eq:frame}
\pi^{n-k}\tilde{\psi}_k=\phi_k \qquad\text{and}\qquad S \tilde{\psi}_k=\tilde{\psi}_k i\,.
\end{equation}
Notice that the frame $\tilde{\psi}_k$ is adapted to the flag $F_k=\ker\pi^{k+1}$ which 
is invariant under $S$. 
\begin{lemma}\label{lem:P-matrix}
Let $P\colon H\to L_n$ be the bundle map $P_p=ev_p\circ P_n$ where $P_n$ is the $n$-th
prolongation map. Let $B$ be the local matrix expression of $P$  near $p_0\in M$ with respect to
the bases $\psi_k$ and $\tilde{\psi}_k$ defined above, i.e., 
$P(\underline{\psi})= \underline{\tilde{\psi}}B$. Then
\[
B=Z(B_0+O(1))W\,
\]
and
\[
dB=Z((-\diag(0,1,\dots,n)B_0+B_0\diag(n_0,n_1,\dots,n_n))\tfrac{dz}{z}+O(0))W\,,
\]
where $Z =\diag(1,z^{-1},\dots,z^{-n})$, $W=\diag(z^{n_0},z^{n_1},\dots,z^{n_n})$
and $B_0$ is an invertible matrix with integer coefficients,
in fact the Wronskian of the independent functions $z^{n_k}$ evaluated at $z=1$.

In particular, the locus where $P$ is not invertible is isolated, so that on a compact 
Riemann surface there are only finitely many Weierstrass points for a given linear system.
\end{lemma}
\begin{proof}
We apply $\pi^{n-l}$ to the equation
\[
P_n(\psi_k)=\sum_{j=0}^{n}\tilde{\psi}_j B_{jk}
\]
and obtain by Corollary~\ref{cor:prolongation}
\[
P_l(\psi_k)=\sum_{j=0}^{l}\pi^{n-l}\tilde{\psi}_j B_{jk}\,,
\]
where we also used that  $\pi^{n-l}\tilde{\psi}_j=0$
for $j>l$ by \eqref{eq:frame}. Inserting \eqref{eq:pre_basis} we thus get
\[
z^{n_k-l}(n_k(n_k-1)\cdots(n_k-l+1)\phi_l+O(1))=
\phi_l  B_{lk} +\sum_{j=0}^{l-1}\pi^{n-l}\tilde{\psi}_j B_{jk}\,,
\]
and, since $\pi^{n-l}\tilde{\psi}_j$ for $j\leq l$ are linearly independent, also
\[
z^{l}B_{lk}z^{-n_k}= n_k(n_k-1)\cdots(n_k-l+1)+O(1)
\]
which proves the claim regarding $B$. 

To obtain the statement for $dB$ we first note that $ZB_0W$ is smooth since
$(B_0)_{lk}=0$ when $l>n_k$. Thus we can write 
\[
B=ZB_0W+ZO(1)W
\]
with a smooth map $ZO(1)W$ which implies
\[
d(ZO(1)W)=ZO(0)W\,.
\]
Taking into account the special form of $Z$ and $W$, one gets the
stated expression of $dB$ by a direct calculation.
\end{proof}
We now are ready to finish the proof of the Pl\"ucker formula: recall that we only need to
verify \eqref{eq:R_hat_integral}
\[
\int_{M}<S\hat{R}>=2\pi(\deg L_n-\ord H)\,
\]
where $\nabla=\hat{\nabla}+A+Q$ is the flat adapted connection \eqref{eq:nabla_m0} over $M_0$ 
induced by $P$. To calculate this curvature integral of the singular complex connection
$\hat{\nabla}$ we compare $\hat{\nabla}$ to a suitably chosen non-singular 
complex connection. Let $z_{\alpha}$ be centered holomorphic coordinates
around the punctures $p_{\alpha}\in M$ and let $B_{\alpha}\subset M$ be the images of 
$\epsilon$-disks under $z_{\alpha}$
so that each $B_{\alpha}$ contains only the Weierstrass point $p_{\alpha}$.
We denote by $M_{\epsilon}$ the Riemann surface with boundary obtained by removing the
disks $B_{\alpha}$ from $M$. Let $\underline{\tilde{\psi}_{\alpha}}$ be the local frame
near the Weierstrass point $p_{\alpha}$ given in \eqref{eq:frame}. Using partition of unity we
define a connection $\tilde{\nabla}$ on $L_n$ over $M$ satisfying 
\begin{equation}\label{eq:tilde_nabla}
\tilde{\nabla}\underline{\tilde{\psi}_{\alpha}}=0 \qquad\text{and}\qquad \tilde{\nabla}S=0\,.
\end{equation}
By the last requirement $\tilde{\nabla}$ is a complex connection so that
\begin{equation}\label{eq:hat_tilde}
\hat{\nabla}=\tilde{\nabla}+\omega
\end{equation}
with $\omega\in\Omega^{1}(M_0,\text{End}_{+}(L_n))$. Therefore, the curvatures are related
by
\[
\hat{R}=\tilde{R}+d^{\tilde{\nabla}}\omega +\omega\wedge\omega\,
\]
over $M_0$.
Multiplying by $S$ and taking the trace \eqref{eq:trace} we get
\[
<S\hat{R}>=<S\tilde{R}>+d<S\omega>\,,
\]
where we used that $S\omega\wedge\omega$ is trace free due to $[S,\omega]=0$. 
We now integrate over $M_{\epsilon}$ and
apply Stokes to obtain
\begin{equation}\label{eq:M_epsilon_integration}
\int_{M_{\epsilon}}<S\hat{R}>=\int_{M_{\epsilon}}<S\tilde{R}>+\int_{\del M_{\epsilon}}<S\omega>\,.
\end{equation}
Since $\tilde{\nabla}$ is a complex connection over $M$, by \eqref{eq:chern} 
the first term on the right hand side 
approaches $2\pi\deg L_n$ as $\epsilon$ tends to zero.
 
This leaves us with calculating the limit of 
\begin{equation}\label{eq:B-integration}
\int_{\del M_{\epsilon}}<S\omega>=-\sum_{\alpha=1}^{N} \int_{\del B_{\alpha}}<S\omega>
\end{equation}
as $\epsilon$ goes to zero.
Let $\underline{\tilde{\psi}_{\alpha}}=\underline{\tilde{\psi}}$ 
be the local frame around the Weierstrass point $p_{\alpha}$
constructed
in \eqref{eq:frame}, and $\underline{\psi}$ the
basis of the linear system $H$ giving rise to that local frame. By \eqref{eq:hat_tilde}
\[
\nabla=\hat{\nabla}+A+Q=\tilde{\nabla}+\omega+A+Q\
\]
and, since $P(\underline{\psi})=\underline{\tilde{\psi}}B$  by Lemma~\ref{lem:P-matrix}, 
we obtain on the punctured disk 
\begin{equation}\label{eq:dB=}
0=\nabla P(\underline{\psi})=\underline{\tilde{\psi}}dB+\omega(\underline{\tilde\psi})B+
(A+Q)(\underline{\tilde\psi})B\,.
\end{equation}
Here we used $\nabla P=0$ and, by construction \eqref{eq:tilde_nabla} of the
connection $\tilde{\nabla}$, also
$\tilde{\nabla}\underline{\tilde{\psi}}=0$. 
Taking the $S$-commuting part of \eqref{eq:dB=} gives 
\begin{equation}\label{eq:omega}
\omega(\underline{\tilde{\psi}})=-\underline{\tilde{\psi}}(dB B^{-1})_{+}\,,
\end{equation}
where $(dB B^{-1})_{+}$ is the $i$-commuting part, i.e., the complex part
in the decomposition $\H=\C\oplus \C j$, since we have chosen 
$S(\underline{\tilde{\psi}})=\underline{\tilde{\psi}}i$. 
From Lemma~\ref{lem:P-matrix} we get
\[
B^{-1}=W^{-1}(B_0^{-1}+O(1))Z^{-1}\,,
\]
therefore
\[
-dB B^{-1}=Z((\diag(0,1,\dots,n)-B_0\diag(n_0,n_1,\dots,n_n){B_0}^{-1})\tfrac{dz}{z}+O(0))Z^{-1}\,
\]
and thus also
\[
-(dB B^{-1})_{+}=
Z((-\diag(0,1,\dots,n)+B_0\diag(n_0,n_1,\dots,n_n){B_0}^{-1})\tfrac{dz}{z}+O(0))Z^{-1}\,. 
\]
This last equation, together with \eqref{eq:omega} and Definition~\ref{def:order_H}, implies 
\[
<S\omega>=i\tfrac{dz}{z}\sum_{k=0}^{n}(k-n_k)+O(0)=-i\tfrac{dz}{z}\ord_{p_{\alpha}}H +O(0)\,
\]
on the disk $B_{\alpha}$. 
Integrating this expression over the boundary of $B_{\alpha}$ yields
\[
\int_{\del B_{\alpha}}<S\omega>=2\pi\ord_{p_{\alpha}}H+O(\epsilon)\,,
\]
and therefore \eqref{eq:B-integration}
\[
\int_{\del M_{\epsilon}}<S\omega>=-2\pi\ord H+O(\epsilon)\,.
\]
Inserting this last into \eqref{eq:M_epsilon_integration} results in
\[
\int_{M}<S\hat{R}>=\lim_{\epsilon\to 0}\int_{M_{\epsilon}}<S\hat{R}>=2\pi(\deg L_n-\ord H)\,,
\]
which we set out to verify. We thus have completed our discussion of the proof of
the Pl\"ucker relation given in Theorem~\ref{thm:Pluecker}.

\subsection{Equality in the Pl\"ucker estimate}
As in the unramified case, we want to characterize those holomorphic line
bundles $L$ for which equality  
\begin{equation}\label{eq:Pluecker_equality}
\tfrac{1}{4\pi}W(L) = (n+1)(n(1-g)-d)+\ord H\,
\end{equation}
holds in \eqref{eq:Pluecker_estimate}. The main technical ingredient in this discussion
will concern the extendibility of the Weierstrass flag $H_k$ into the Weierstrass points.
\begin{lemma}\label{lem:flag_extension}
Let $L$ be a quaternionic holomorphic line bundle, $H\subset H^{0}(L)$ an $(n+1)$-dimensional
linear system, and $M_0$ the open Riemann surface punctured at the Weierstrass points $p_{\alpha}$
of $H$. Let $P\colon H\to L_n$ be the bundle map  given in \eqref{eq:ev_prolongation} and $S$ be 
the canonical complex structure on $H$ over $M_0$ induced by $P$. Then the Weierstrass flag 
$H_k$ and the complex structure $S$
extend continuously into the Weierstrass points $p_{\alpha}$, where $H_k$ has limit
$H_k(p_{\alpha})$. 
\end{lemma}
\begin{proof}
Let $p_{\alpha}\in M$ be a Weierstrass point and $\psi_k$, $k=0,\dots,n$, a basis of $H$
whose vanishing orders at $p_{\alpha}$ are the Weierstrass gap sequence $n_k$ for $p_{\alpha}$.
Let $\tilde{\psi}_k$ be the corresponding local framing of $L_n$ constructed in 
\eqref{eq:frame}. Then both frames are adapted, meaning that
$\psi_{n-k},\dots,\psi_n$ are a basis of $H_{k}(p_{\alpha})$ and 
$\tilde{\psi}_{n-k},\dots,\tilde{\psi}_n$ frame $F_k$ in a neighborhood around $p_{\alpha}$. In 
Lemma~\ref{lem:weierstrass_char} we have seen that, away from Weierstrass points, 
$P$ maps the Weierstrass flag $H_k$ isomorphically to $F_k$. Therefore, expressing
$P(\underline{\psi})=\underline{\tilde{\psi}}B$ as in Lemma~\ref{lem:P-matrix},
we see that $\underline{\psi}B^{-1}$ is an adapted local frame for the flag $H_k$
away from $p_{\alpha}$. Our aim is to modify this frame so that it stays adapted
to the flag $H_k$ but approaches $\underline{\psi}$ as $p$ approaches the Weierstrass point 
$p_{\alpha}$.
Using the same notation as in Lemma~\ref{lem:P-matrix} we have $B=Z\tilde{B}W$
with $\tilde{B}=B_0+O(1)$. We decompose 
\[
\tilde{B}=\tilde{U}\tilde{L}
\]
into upper, with ones along the diagonal, and lower diagonal matrices. Then the upper
diagonal matrix
\[
U:=Z\tilde{U}Z^{-1}
\]
approaches the identity as $p$ goes to $p_{\alpha}$. With the lower diagonal matrix
\[
L:=Z\tilde{L}W\,,
\]
which fails to be invertible at $p_{\alpha}$, we obtain the upper-lower decomposition 
\[
B=UL
\]
of $B$. Since $L$ is lower diagonal, $\underline{\hat{\psi}}:=\underline{\psi}L$
is also an adapted basis for the flag $H_k(p_{\alpha})$ in a punctured neighborhood
of $p_{\alpha}$. But then 
\[
\underline{\hat{\psi}}B^{-1}=\underline{\psi}U^{-1}
\]
is an adapted frame for the Weierstrass flag $H_k$ near $p_{\alpha}$ approaching 
$\underline{\psi}$ as $p$ goes to $p_{\alpha}$. Thus the Weierstrass flag
$H_k$ approaches $H_k(p_{\alpha})$.

To calculate the limit of $S$ on $H$ as $p$ tends to $p_{\alpha}$ we again work with the frames
$\underline{\psi}$ in $H$ and $\underline{\tilde{\psi}}$ in $L_n$. From \eqref{eq:frame} we 
know that the matrix
of $S$ on $L_n$ in the frame  $\underline{\tilde{\psi}}$ is given by
\[
S(\underline{\tilde{\psi}})=\underline{\tilde{\psi}}I_{n+1}i\,.
\]
Thus the matrix of $S$ on $H$ in the frame $\underline{\psi}$ becomes
\[
P^{-1}SP(\underline{\psi})=\underline{\psi}B^{-1}I_{n+1}iB\,.
\]
From $B=Z(B_0+O(1))W$ and the fact that $Z$, $W$ and $B_0$ commute with $i$, we obtain
\[
B^{-1}I_{n+1}iB=W^{-1}(B_0^{-1}+O(1))Z^{-1}I_{n+1}iZ(B_0+O(1))W=I_{n+1}i+O(1)\,,
\]
which shows that $S$ on $H$ has a limit as $p$ tends to $p_{\alpha}$.
\end{proof}
Now let $L$ be a holomorphic line bundle with an $n+1$-dimensional linear system
$H\subset H^{0}(L)$ for which the equality \eqref{eq:Pluecker_equality} holds. 
By Theorem~\ref{thm:Pluecker} we see that this is the case precisely when
$W(L^{*})=0$, i.e., when $A=0$ on $M_0$. But
$P\colon H\to L_n$ 
is an isomorphism  over $M_0$  mapping the 
Weierstrass flag $H_k\subset H$ into the Frenet flag $F_k=\ker\pi^{k+1}\subset L_n$.  
Since $A=0$ we can apply Lemma~\ref{lem:twistor_char} and Lemma~\ref{lem:Frenet_curve_0} to conclude   
that $H_k=\pi(W_k)$ is the twistor projection \eqref{eq:twistor_proj} of
the complex holomorphic curve $W_k\colon M_0\to G_k^{*}(H,i)$, i.e., $H_k=W_k\oplus W_k j$,
where  $W_k\subset H_k$ denotes the $i$-eigenspace of $S$ on $H$. Moreover, $W_k$ are the first $n$ 
osculating curves of the complex holomorphic curve $W_0\colon M_0\to \C\P(H,i)$.
From Lemma~\ref{lem:flag_extension} we know that $H_k$ and $S$ extend continuously into the 
Weierstrass points, and hence  
the $i$-eigenspaces $W_k$ also extend continuously. But then 
$W_k$ extends complex holomorphically and we obtain a complex holomorphic curve
$W_0$, together with its first $n$ osculating curves $W_k$, on all of $M$. This implies that
both $H_k=W_k\oplus W_k j$ and $S$ are in fact smooth across the Weierstrass points.
In particular, the derivatives of the higher osculating curves
\[
\hat{\delta}_k\colon H_k/H_{k-1}\to K H_{k+1}/H_k
\]
are obtained from the first $n$ higher derivatives of the complex holomorphic curve $W_0$ by
quaternionic extension so that
\[
*\hat{\delta}_k=S\hat{\delta}_k=\hat{\delta}_{k}S\,.
\]
In other words, the flag $H_k$ is holomorphic, both as curves and as subbundles, 
by Corollary~\eqref{cor:holo_compare}.
 
We now find ourselves in the following set up: the trivial bundle $H$ over $M$ has the smooth
complex structure $S$ such that $PS=SP$. The $\bar{K}$-part $\nabla''=\delbar +Q$ 
of the trivial connection $\nabla$ defines a holomorphic structure on $H$. Over $M_0$
we have $\nabla P=P\nabla$, and taking $\bar{K}$-parts we also have
\[
(\delbar +Q_n)P=P(\delbar+Q)\,.
\]
Here $\delbar +Q_n$ denotes the
canonical holomorphic structure on the $n$-th jet bundle $L_n$ coming from 
the holomorphic structure $\delbar+Q_{L}$ on $L$ as described in 
Theorem~\ref{thm:quat_holo_jet_complex}.
Since all occurring  objects are defined over $M$, this last relation extends
to $M$ and the bundle map $P\colon H\to L_n$
is holomorphic, i.e., $P$ intertwines $\delbar$'s and $Q$'s. 
In particular, $P$ maps the holomorphic flag
$H_k$ into the holomorphic flag $F_k$. From Corollary~\ref{cor:adapted_nabla_relations}
we know that
$Q_{n}$ vanishes on $F_{n-1}$ and thus $Q$ also vanishes on $H_{n-1}$. This implies
that $P\colon H\to L_n$ induces the holomorphic bundle maps 
\begin{equation}\label{eq:P_k}
P_k\colon H_k/H_{k-1}\to F_k/F_{k-1}=N_{n-k}\,,
\end{equation}
where the latter identification is done by $\pi^{k}$.
In particular, 
\begin{equation}\label{eq:image_P}
P_n\colon H/H_{n-1}\to L_n/F_{n-1}=L
\end{equation}
is a holomorphic bundle map from the inverse $H/H_{n-1}=(H_0^{*})^{-1}$ 
of the dual curve $H_0^{*}=H_{n-1}^{\perp}\subset H^{-1}$ of 
the holomorphic curve $H_0$. Since the holomorphic curve $H_0^{*}\subset H^{-1}$ is full, its
induced linear system
\[
\hat{H}=\{\psi+H_{n-1}\,;\,\psi\in H\}\subset H^{0}(H/H_{n-1})
\]
is $n+1$-dimensional and base point free. It is easy to check that the holomorphic
bundle map \eqref{eq:image_P} induces the linear isomorphism
\[
P_n(\psi+H_{n-1})=\psi
\]
between the linear systems $\hat{H}\subset H^{0}(H/H_{n-1})$ and $H\subset H^{0}(L)$.
 
Let us summarize what we have shown so far:
a holomorphic line bundle $L$ with an $n+1$-dimensional linear system $H$ for which  
equality holds in the estimate \eqref{eq:Pluecker_estimate} is, up to the holomorphic bundle map
\eqref{eq:image_P}, given by the induced linear system  of the dual curve $H_0^{*}$ 
of the twistor projection of a complex holomorphic curve $W_0$ in $\C\P(H,i)$. The 
$n$-th osculating curve $W_n$ of $W_0$ has to satisfy  $W_n\oplus W_n j=H$. 

It remains
to calculate the degree and the first $n$ Weierstrass gaps of the complex holomorphic curve 
$W_0$ in terms of the corresponding data of $L$ and $H$.     
Consider the diagram
\begin{equation}\label{eq:diagram}
\begin{CD}
H_{k}/H_{k-1}@>{\hat{\delta}_k}>>K H_{k+1}/H_{k}\\
@V{P_k}VV @VV{P_{k+1}}V\\
F_k/F_{k-1}@>\delta_{n-k}>> K F_{k+1}/F_{k}
\end{CD}
\end{equation}
where the bottom row is given by the
bundle isomorphisms $\delta_k\colon N_k\to K N_{k-1}$ of the holomorphic jet complex of $L$,
as in \eqref{eq:jet_delta}, and the vertical maps are given in \eqref{eq:P_k}.
Since the diagram commutes away from  
Weierstrass points, it commutes on all of $M$. The vanishing orders of $P_k$ and
$\hat{\delta}_k$ at any given point are thus related by
\begin{equation}\label{eq:vanishing_orders}
\ord\hat{\delta}_k=\ord P_k-\ord P_{k+1}\,.
\end{equation}
To calculate $\ord P_k$ at a given point $q\in M$, we again take a basis $\underline{\psi}$ of $H$
whose vanishing orders are given by the Weierstrass gaps at this point, i.e., 
$\ord\psi_k=n_k$. Then $\underline{\psi}$ is an adapted basis for the Weierstrass
flag at $q$. Let $\underline{\hat{\psi}}=\underline{\psi}C$ be a smooth local
extension to an adapted frame  $\underline{\hat{\psi}}$ of $H_k$ near $q$. If 
$\underline{\tilde{\psi}}$ is the frame \eqref{eq:frame} in $L_n$ adapted to the
flag $F_k$, then
\[
P(\underline{\hat{\psi}})=\underline{\tilde{\psi}}BC
\]
with $B=Z(B_0+O(1))W$ from Lemma~\ref{lem:P-matrix}. Since $P(H_k)\subset F_k$ the matrix $BC$ is
lower diagonal and 
\[
(BC)_{k,k}=z^{n_k-k}(B_0)_{k,k}+O(1)\,,
\]
where we used the fact $C=I+O(1)$. Because $\underline{\hat{\psi}}$ is adapted to
the flag $H_k$, the local sections $\hat{\psi}_{n-k},\cdots,\hat{\psi}_{n}$ frame 
$H_k$ so that 
\begin{equation}\label{eq:order_P_k}
\ord P_k=n_{n-k}-(n-k)
\end{equation}
at $q$. Together with \eqref{eq:vanishing_orders} we thus obtain
\begin{equation}\label{eq:delta_orders}
\ord\hat{\delta}_k=n_{n-k}-n_{n-k-1}-1
\end{equation}
for the vanishing orders of the higher order derivatives $\hat{\delta}_k$ of the holomorphic
curve $H_0$. Let us denote the Weierstrass gap sequence of the curve $H_0\subset H$ by
$n^{*}_k$. To calculate $\ord\hat{\delta}_k$ in terms of $n^{*}_k$ we recall that the
derivatives of the dual flag $H_k^{\perp}$ are given by
\[
\delta^{\perp}_k=-\hat{\delta}^{*}_{k}\,.
\] 
We denote the osculating flag of the dual curve $H_0^{*}=H_{n-1}^{\perp}$ by 
$H_k^{*}=H_{n-k-1}^{\perp}$ so that the derivatives of this flag are $-\hat{\delta}^{*}_{n-k-1}$.
Applying \eqref{eq:delta_orders} to the smooth holomorphic flag $H_k^{*}$ we obtain
\[
\ord\hat{\delta}_k=\ord\hat{\delta}^{*}_k=n^{*}_{k+1}-n^{*}_k-1\,,
\]
and thus
\[
n^{*}_{k+1}-n^{*}_k=n_{n-k}-n_{n-k-1}\,.
\]
Telescoping these identities gives the desired relationship 
\begin{equation}\label{eq:gap_mirror}
n^{*}_{k}-n^{*}_0=n_n-n_{n-k}
\end{equation}
between the Weierstrass gaps of the holomorphic curve $H_0$ and its dual curve $H_0^{*}$.
In our case $n_0^{*}=0$ since $H_0$ is a smooth curve in projective space, in fact the
twistor projection of the complex holomorphic curve $W_0$ in $\C\P(H,i)$. Since the 
derivatives of the quaternionic flag $H_k$ are just the quaternionic extensions of the
first $n$ derivatives of the complex holomorphic curve $W_0$, we deduce that 
$W_0$ must have the first $n$ Weierstrass gaps $n_n-n_{n-k}$.

Finally we calculate the degree $d^{*}$ of the curve $W_0$. First note that
the curve $H_0^{*}=H_{n-1}^{\perp}$ has degree $d-\sum_{p}n_0(p)$ since the holomorphic bundle map 
\eqref{eq:image_P}
\[
P_n\colon (H_{n-1}^{\perp})^{-1}\to L
\]
vanishes to order $n_0(p)$ at $p\in M$ by \eqref{eq:order_P_k}.
Specifically, if $L$ had no base points then $(H_{n-1}^{\perp})^{-1}$
would be isomorphic to $L$ via $P_n$. The degrees of the bundle maps 
\[
\hat{\delta}_k\colon H_k/H_{k-1}\to K H_{k+1}/H_k
\]
are given by
\[
\deg\hat{\delta}_k=\deg K +d_{k+1}-d_k
\]
where $d_k$ denotes the degree of the line bundle $H_k/H_{k-1}$, in particular
$d^{*}=-d_0$ and $d_{n}=d-\sum_{p}n_0(p)$. Telescoping the latter identity, and recalling
\eqref{eq:delta_orders}, we finally get
\begin{equation}\label{eq:dual_degree}
\sum_{p}(n_n(p)-n)=d+d^{*}+n\deg K\,.
\end{equation}
We are now ready to formulate the analog of Corollary~\ref{cor:Willmore_equality}
in the presence of Weierstrass points, namely the characterization of
those line bundles for which we have equality in \eqref{eq:Pluecker_estimate}. 
\begin{theorem}\label{thm:equality}
Let $L$ be a quaternionic holomorphic line bundle of degree $d$ over a
compact Riemann surface of genus $g$ and $H\subset H^{0}(L)$ an $n+1$-dimensional
linear system with Weierstrass gap sequence $n_k$. 
Then
\[
\tfrac{1}{4\pi}W(L)=(n+1)(n(1-g)-d)+\ord H
\]
if and only if there is a holomorphic line bundle $\hat{L}$ and a holomorphic bundle map 
\[
P\colon \hat{L}\to L\
\]
vanishing to order $n_0(p)$ at $p\in M$, where
$\hat{L}^{-1}$ is the dual curve of a twistor projection of a 
complex holomorphic curve $E$ in $\C\P^{2n+1}$ of degree $\sum_{p}(n_n(p)-n)-d-n\deg K$
satisfying
\begin{enumerate}
\item there exists a quaternionic structure $j$, i.e, a complex antilinear map with
$j^{2}=-1$, on $\C^{2n+2}$ so that the $n$-th osculating bundle $E_n$ of $E$ satisfies
$E_n\oplus j(E_n)=\C^{2n+2}$, and
\item
the first $n$ elements of the Weierstrass gap sequence $n^{*}_k$ of $E$ are 
$n^{*}_k=n_{n}-n_{n-k}$\,. 
\end{enumerate} 
In this setting $\H^{n+1}=\C^{2n+2}$ via $j$, the twistor projection of $E$
is the holomorphic curve $E\oplus j(E)$ in $\H\P^{n}$, and $\hat{L}^{-1}=(E\oplus j(E))^{*}$.

In particular, if $L$ had no base points to start with, then $L$ would simply be the 
dual line bundle of the dual curve of a twistor projection of a complex holomorphic
curve in $\C\P^{2n+1}$ with the above properties. In this sense holomorphic line bundles
for which equality \eqref{eq:Pluecker_equality} holds in the estimate \eqref{eq:Pluecker_estimate}
are obtained from complex holomorphic data.
\end{theorem}
To finish the proof of the theorem we have to check that, starting with a complex holomorphic
curve with the right properties, we obtain equality in the Pl\"ucker estimate 
\eqref{eq:Pluecker_estimate}. Let $\H^{n+1}=\C^{2n+2}$, as usual, and consider a complex
holomorphic curve $E\subset \C^{2n+2}$ whose $n$-th osculating bundle $E_n$ satisfies
\[
E_n\oplus E_n j=\H^{n+1}\,.
\]
Define the complex structure $S$ on $\H^{n+1}$ by requiring $S$ to be multiplication by 
$i$ on $E_n$. Then the twistor projections $H_k:=E_k\oplus E_k j$, where $E_k$ are the 
first $n$ osculating bundles of $E$, have derivatives $\delta_k$ satisfying 
\[
*\delta_k=S\delta_k=\delta_k S\,.
\]
The dual curve $H_0^{*}=H_{n-1}^{\perp}$ of $H_0$ is holomorphic and we denote by
$\hat{L}$ the
holomorphic line bundle $(H_0^{*})^{-1}$ with induced 
$(n+1)$-dimensional linear system $\hat{H}\subset H^{0}(\hat{L})$. 
From our discussion above, and the assumptions in the theorem,
we deduce that the Weierstrass gap sequence for $\hat{H}$ is given by
\[
\hat{n}_k=n^{*}_n-n^{*}_{n-k}=n_{n}-n_{0}-n_{n}+n_{k}=n_{k}-n_{0}\,.
\]
We now apply the Pl\"ucker formula \eqref{eq:Pluecker} to the line bundle $\hat{L}$
to get
\[
\tfrac{1}{4\pi}W(\hat{L})=(n+1)(n(1-g)-\hat{d})+\ord\hat{H}\,,
\]
where we recall that $W(\hat{L}^{*})=0$,  because $\hat{L}^{*}=H_0^{-1}$ was obtained
from a twistor projection of a complex holomorphic curve. Since the bundle map
$P$ is assumed to vanish to order $n_0(p)$ at $p\in M$, it maps the
linear system $\hat{H}$ to $H$ and the line bundle $\hat{L}$ has degree
$\hat{d}=d-\sum_{p} n_0(p)\,$. The holomorphic bundle map $P\colon \hat{L}\to L$ intertwines
the Hopf fields
\[
P\hat{Q}=QP
\]
so that the respective Willmore functionals 
\[
W(L)=W(\hat{L})
\]
agree. But then
\begin{align*}
\tfrac{1}{4\pi}W(L)&=\tfrac{1}{4\pi}W(\hat{L})\\
&=(n+1)(n(1-g)-\hat{d})+\ord\hat{H}\\
&=(n+1)(n(1-g)-d)+\ord{H}\,,
\end{align*}
which finishes the proof.
\subsection{Estimates on the Willmore energy}\label{subsec:Willmore_estimates}
In this section we want to derive a general estimate for the 
Willmore energy of a holomorphic line bundle $L$ of degree $d$ over a compact Riemann
surface of genus $g$ in terms of the number $n+1$
of holomorphic sections in a given linear system $H\subset H^{0}(L)$.
As it stands, the estimate \eqref{eq:Pluecker_estimate} is not useful for higher genus 
due to the dominating term $(1-g)n^2$. 

Let us begin by discussing the case of genus zero:
since $\ord H$ is non-negative we get
\begin{equation}\label{eq:Willmore_estimate_0}
W(L)\geq 4\pi (n+1)(n-d)\,.
\end{equation}
If the degree $d$ of $L$ is non-negative then the unique complex holomorphic line bundle of degree
$d$ has exactly $d+1$ many holomorphic sections by Riemann-Roch. Hence $n=d$ in this case,
and we have the trivial estimate. On the other hand, if the degree of $L$ is negative
then, to obtain holomorphic sections, the Willmore energy $W(L)$ has to be at least
$4\pi (n+1)(n-d)$, i.e., the Willmore energy grows quadratically in the number of sections.
Using Theorem~\ref{thm:equality} it can be shown that the above estimate 
\eqref{eq:Willmore_estimate_0} is sharp. A geometrically
interesting special case occurs when $L$ has spin bundle degree $d=-1$. Then the estimate
becomes
\begin{equation}\label{eq:Willmore_spin_0}
W(L)\geq 4\pi (n+1)^2\,,
\end{equation}
and equality holds for spin bundles induced from special conformally immersed spheres in 3-space,
the so called {\em Dirac spheres} and {\em soliton spheres} \cite{Richter}, \cite{Taimanov}.

We now consider the case of non-zero genus. To still obtain an estimate which is
quadratic in $n$, we utilize the term $\ord H$ in \eqref{eq:Pluecker_estimate}: let
$H_{k} \subset H$ be the $(k+1)$-dimensional linear subspace  defined by the $n-k$
independent linear relations on $H$ guaranteeing the existence of a holomorphic
section $\psi\in H_k$ with a zero of order $n-k$ at a given point. Then the Weierstrass
gap sequence of the linear system $H_k\subset H$ at this point is at least
\[
n-k<n-k+1<\dots <n
\]
so that 
\[
\ord H_k \geq (k+1)(n-k)\,.
\]
Applying \eqref{eq:Pluecker_estimate} to the linear system $H_k$ we thus obtain
\[
\tfrac{1}{4\pi}W(L)\geq (k+1)(k(1-g)-d)+ (k+1)(n-k)=(k+1)(n-d-kg)
\]
for all $k=0,\dots, n$. In other words, for each $k$ the linear function
\begin{equation}\label{eq:f_k}
f_k(n):=(k+1)(n-d-kg)
\end{equation}
is bounding the Willmore energy from below. It is easy to check that the parabola
\[
f(n):=\tfrac{1}{4g}((n+g-d)^2-g^2)
\]
contains the points where $f_k=f_{k+1}$ through which the $f_k$'s cut as secants. 
Therefore
we get the lower bound
\begin{equation}\label{eq:Willmore_estimate_1}
W(L)\geq \tfrac{\pi}{g}((n+g-d)^2-g^2)
\end{equation}
for any  holomorphic line bundle $L$ over a compact Riemann surface
of genus $g\geq 1$ admitting $n+1$ independent holomorphic sections.

Since {\em complex} holomorphic line bundles of negative degree never have holomorphic sections
we need $n\geq 0$ for the estimate to hold, i.e., the bundle in question must have
at least one holomorphic section.  Furthermore, for any complex holomorphic line bundle $L$ 
we always have
\[
d=\deg L\geq \dim H^{0}(L)-1=n\,.
\]
This is due to the fact that $n+1$ independent holomorphic sections produce a holomorphic 
section with a zero of degree at least $n$ by a suitable linear combination. Thus the estimate
\eqref{eq:Willmore_estimate_1} holds in the region 
\[
n\geq 0\qquad\text{and}\qquad n\geq d
\]
and gives the trivial estimate
along $n=d$, which agrees with the complex holomorphic situation.
 
The Riemann-Roch paired
bundle $KL^{-1}$ has degree $\tilde{d}=2g-2-d$ and, by Riemannn-Roch \eqref{eq:Riemann-Roch},
 admits 
$\tilde{n}+1=n-d+g$ many independent holomorphic sections. 
By Definition~\ref{def:Will_energy} the Willmore energy of
$KL^{-1}$ is the same as the one for $L$. Applying \eqref{eq:Willmore_estimate_1} to $KL^{-1}$
yields the further estimate
\begin{equation}\label{eq:Willmore_estimate_2}
W(L)\geq \tfrac{\pi}{g}((n+1)^2-g^2)
\end{equation}
of the Willmore energy for a holomorphic line bundle admitting $n+1$ 
independent holomorphic sections. Translating the restrictions on $n$ and $d$ to $\tilde{n}$ and
$\tilde{d}$ we see that this second estimate holds in the region
\[
n\geq g-1\qquad\text{and}\qquad n\geq d-g+1\,.
\]
Again, along $n=g-1$ and $n=d-g$ there are complex holomorphic line bundles
so that we cannot obtain an estimate for the Willmore energy.

Finally, in the region $n\geq d$ and $n\geq g-1$, where both estimates apply,
\eqref{eq:Willmore_estimate_1} gives the stronger estimate
for $d<g-1$, whereas \eqref{eq:Willmore_estimate_2} gives the stronger estimate
for $d>g-1$. In the case of a spin bundle, where degree $d=g-1$, the two estimates agree.

The diagram below illustrates the possible range
for the number of holomorphic sections of a holomorphic line bundle
in terms of the degree. 
\begin{center}     
\begin{picture}(0,0)%
\epsfig{file=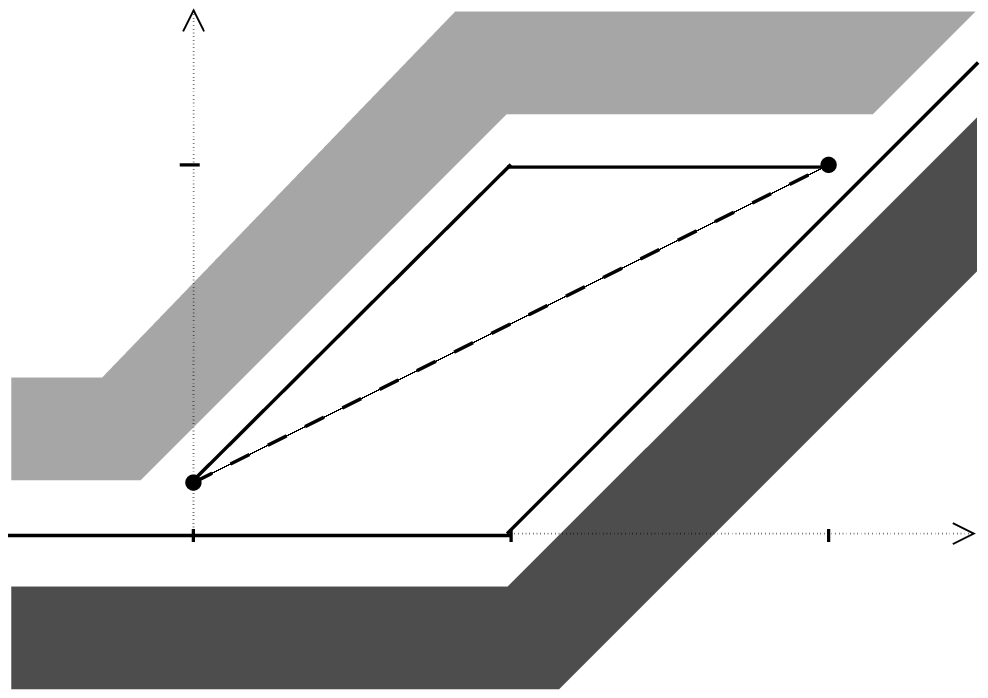}%
\end{picture}%
\setlength{\unitlength}{0.00083300in}%
\begingroup\makeatletter\ifx\SetFigFont\undefined
% extract first six characters in \fmtname
\def\x#1#2#3#4#5#6#7\relax{\def\x{#1#2#3#4#5#6}}%
\expandafter\x\fmtname xxxxxx\relax \def\y{splain}%
\ifx\x\y   % LaTeX or SliTeX?
\gdef\SetFigFont#1#2#3{%
  \ifnum #1<17\tiny\else \ifnum #1<20\small\else
  \ifnum #1<24\normalsize\else \ifnum #1<29\large\else
  \ifnum #1<34\Large\else \ifnum #1<41\LARGE\else
     \huge\fi\fi\fi\fi\fi\fi
  \csname #3\endcsname}%
\else
\gdef\SetFigFont#1#2#3{\begingroup
  \count@#1\relax \ifnum 25<\count@\count@25\fi
  \def\x{\endgroup\@setsize\SetFigFont{#2pt}}%
  \expandafter\x
    \csname \romannumeral\the\count@ pt\expandafter\endcsname
    \csname @\romannumeral\the\count@ pt\endcsname
  \csname #3\endcsname}%
\fi
\fi\endgroup
\begin{picture}(4700,3287)(715,-2736)
\put(4125,-134){\makebox(0,0)[lb]{\smash{\SetFigFont{5}{6.0}{rm}canonical bundle}}}
\put(1740,-1827){\makebox(0,0)[lb]{\smash{\SetFigFont{5}{6.0}{rm}trivial bundle}}}
\put(2915,-2160){\makebox(0,0)[lb]{\smash{\SetFigFont{5}{6.0}{rm}$g-1$}}}
\put(1575,-2166){\makebox(0,0)[lb]{\smash{\SetFigFont{5}{6.0}{rm}$0$}}}
\put(1338,-220){\makebox(0,0)[lb]{\smash{\SetFigFont{5}{6.0}{rm}  $g$}}}
\put(999,413){\makebox(0,0)[lb]{\smash{\SetFigFont{5}{6.0}{rm}$\dim H^0(L)$}}}
\put(5142,-2160){\makebox(0,0)[lb]{\smash{\SetFigFont{5}{6.0}{rm}$\deg L$}}}
\put(4586,-2166){\makebox(0,0)[lb]{\smash{\SetFigFont{5}{6.0}{rm}$2g-2$}}}
\end{picture}

\end{center}
There are no holomorphic line bundles in the dark shaded area. For a holomorphic line
bundle in the light shaded area the Willmore energy is bounded from below by 
the above two estimates. The bounds are quadratic in the number of holomorphic sections. 
By an argument more subtle then Riemann-Roch, namely the Clifford estimate \cite{bible},
one can exclude half of the parallelogram region of
possible {\em complex} holomorphic line bundles above the line  
connecting the trivial bundle with the canonical bundle. Whether there can 
be quaternionic holomorphic line bundles in this region and, if so, which estimates
their Willmore energies satisfy, is at present unknown to us.
We conclude this section by summarizing the above discussion:
\begin{theorem}\label{thm:Willmore_estimate}
Let $L$ be a holomorphic line bundle of degree $d$ over a compact Riemann surface of genus $g$.
Then we have the following estimates on the Willmore energy in terms of the
number $n+1$ of holomorphic sections:
\begin{enumerate}
\item
If $g=0$ we have
\[
W(L)\geq 4\pi(n+1)(n-d)\,,
\]
which is sharp.
\item
If $g\geq 1$ then
\begin{equation}
W(L)\geq
\begin{cases}
\tfrac{\pi}{g}((n+g-d)^{2}-g^{2}) &\text{if}\qquad n\geq 0\,,\, n\geq d\,,\, d\leq g-1\\
\tfrac{\pi}{g}((n+1)^2-g^2) &\text{if}\qquad n\geq d-g+1\,,\, n\geq g-1\,,\, d\geq g-1
\end{cases}
\end{equation}
\end{enumerate}
\end{theorem}
For later applications it will be useful to refine these estimates in the following way:
the parabola $f(n)=\tfrac{1}{4g}(n+g-d)^2$ has as tangents $f_k$, given by
\eqref{eq:f_k}, touching at
$n=d+(2k+1)g$, $k\in\N$. At these special values of $n$ we obtain  
better general estimates for the Willmore energy which we collect for the
following case:
\begin{remark} 
Let $L$ be a holomorphic line bundle of degree $d=0$ over a torus, $g=1$,
admitting $n+1$ independent holomorphic sections.
Then the Willmore energy is bounded from below by
\begin{equation}\label{eq:special_W_estimates}
W(L)\geq
\begin{cases}
\pi(n+1)^2 &\,\,\text{if $n$ is odd, and}\\
\pi((n+1)^2-1)&\,\, \text{if $n$ is even}\,.
\end{cases}
\end{equation}
\end{remark}
This concludes, for the time being, our development of quaternionic holomorphic
geometry. The final two sections of this paper deal with
applications of our theory, especially the general estimates on the Willmore
energy just given, to Dirac eigenvalue estimates and energy estimates
of harmonic tori.

%%% Local Variables:
%%% mode: latex
%%% TeX-master: "willmore"
%%% End: 

\section{Dirac eigenvalue estimates}\label{sec:Dirac}
Our first application of the quaternionic Pl\"ucker formula concerns eigenvalue estimates of
Dirac operators over $2$-dimensional compact surfaces.
To keep the exposition reasonably self contained, and for the 
benefit of readers who are not familiar with
the general theory of spin bundles and Dirac operators \cite{Lawson}, we briefly develop
the necessary terminology from scratch in the language of the present paper.
\subsection{Spin bundles and Dirac operators}
Throughout this section $M$ will denote an oriented, real $2$-dimensional manifold with a
given Riemannian metric $(\,,\,)$.
\begin{definition}\label{def:spin-bundle}
A {\em Riemannian spin bundle} over $M$ is given by the following data:
\begin{enumerate}
\item  
A quaternionic line bundle $L$ over $M$;
\item
A bundle map, the {\em Clifford multiplication}, of the tangent bundle  
\[
\hat{} \, : TM\to \text{End}_{\H}(L)
\]
into the endomorphisms of $L$ satisfying
\[
\hat{X}\hat{Y}+\hat{Y}\hat{X}=-2(X,Y)id_{L}
\]
for any tangent vectors $X, Y\in TM$;
\item
A hermitian inner product $(\,,\,)$ on $L$, i.e., $(\,,\,)$ is non-degenerate, non-negative
on the diagonal, 
\[
\overline{(\psi,\phi)}=(\phi,\psi)\,,\qquad\text{and}\qquad
(\psi\lambda,\phi\mu)=\bar{\lambda}(\psi,\phi)\mu
\]
for $\psi,\phi\in L$ and $\lambda,\mu\in \H$.
\end{enumerate}
Two Riemannian spin bundles are isomorphic 
if there exists a quaternionic linear isomorphism between them intertwining their
respective structures.
\end{definition}
If $M$ is  simply connected, all  Riemannian spin bundles are isomorphic.
In the non-simply connected case, the obstruction for two Riemannian spin bundles to be 
isomorphic is carried by the Clifford multiplication only. If $M$ is compact there
are exactly $2^{2g}$ different Riemannian spin bundles.

Any orthonormal basis $X, Y\in T_p M$ gives rise to a complex structure
\begin{equation}\label{eq:complex-structure}
J_p=\hat{Y}\hat{X}
\end{equation}
on $L_p$, which is independent of the choice of orthonormal basis. In this way we get
-- at least up to sign -- a canonical complex structure $J\in\Gamma(\text{End}(L))$, $J^2=-1$, 
on every Riemannian spin bundle. Notice that the bundle map $\hat{}\colon TM\to \text{End}(L)$
is complex antilinear with respect to the Riemann surface
complex structure on $TM$ -- its metric and orientation make $M$ into a Riemann surface -- 
and with respect to the left complex structure on $\text{End}(L)$ given
by post-composition with $J$. Moreover, if $X\in TM$ is a tangent vector then
$\hat{X}\in \text{End}_{-}(L)$ anticommutes with $J$.
\begin{lemma}\label{lem:basic_properties}
Let $L$ be a Riemannian spin bundle. Then
the complex structure $J$ and $\hat{X}$, for $X\in TM$, are skew-adjoint with respect
to the hermitian product.
\end{lemma}
\begin{proof}
Let 
$T\in\text{End}(L)$ be an endomorphism with $T^2=-\rho\, id_L$ where $\rho\geq 0$.  
Then
$T\psi = \psi\lambda$ with $ \bar{\lambda}=-\lambda$
for non-zero $\psi\in L$, and hence
\[
(T\psi,\psi)
=(\psi\lambda,\psi)=\bar{\lambda}(\psi,\psi)=(\psi,\psi)(-\lambda)\,,
=(\psi,-T\psi)
\]
which means that $T$ is skew-adjoint with respect to the given hermitian product.
In particular, $J$ and $\hat{X}$, for $X\in TM$, are skew-adjoint.
\end{proof}
The final ingredient of a Riemannian spin bundle is its canonical connection:
\begin{lemma}\label{lem:spin_connection}
Let $L$ be a Riemannian spin bundle. Then there is a unique quaternionic
connection $\nabla$ on $L$, the {\em spin connection}, such that
\begin{enumerate}
\item
Clifford multiplication is a parallel bundle map with respect to
the Levi-Civita connection $\nabla$ on $TM$ and with respect to
the connection on $\text{End}(L)$ induced by 
$\nabla$ on $L$, i.e.,
\[
\nabla\hat{X}=\widehat{\nabla X}
\]
for vector fields $X$ on $M$;
\item
the hermitian inner product is parallel, i.e.,
\[
\nabla(\,,\,)=0\,.
\]
\end{enumerate}
\end{lemma}
Note that the complex structure $J$ on $L$ is parallel with respect to the spin connection,
since
\[
\nabla J=\nabla(\hat{Y}\hat{X})=(\nabla\hat{Y})\hat{X}+\hat{Y}(\nabla\hat{X})=0\,.
\]
Here we used that $\nabla X=Y\sigma$ and $\nabla Y=-X\sigma$ for any 
orthonormal local frame
$X,Y$ on $M$, where $\sigma$ is the connection $1$-form of the given frame. 
Thus the spin connection
is a complex connection on the complex rank $2$ bundle $L$. 
\begin{proof}
It is always possible to choose a quaternionic connection $\nabla$ on $L$
which makes $(\,,)$ parallel.
Then any other such connection is of the form $\nabla+\omega$ with 
$\omega$ a skew adjoint $1$-form on $M$ with values in $\text{End}(L)$.
Using this freedom, we choose $\omega$ such that $\nabla+\omega$ becomes
a complex connection, i.e., makes $J$ parallel: the equation
\[ 
0=(\nabla+\omega)J=\nabla J+[\omega,J]
\]
has the $\text{End}_{-}(L)$ valued $1$-form
\[
\omega=\tfrac{1}{2}(\nabla J)J
\]
as a solution. Since any element in $\text{End}_{-}(L)$ squares to a negative multiple of the 
identity, Lemma~\ref{lem:basic_properties} implies that $\omega$ is in fact skew-adjoint.

Notice that we have fixed $\nabla$ up to the addition of $\omega=J\alpha$ for some
real $1$-form $\alpha\in\Omega^1(\R)$.
Now Clifford multiplication is parallel with respect to $\nabla+\omega$ if and only if
\[
[\omega, \hat{X}]=\widehat{\nabla X}-\nabla\hat{X}
\]
or, equivalently
\begin{equation}\label{eq:spin-nabla}
2\alpha J\hat{X}= \widehat{\nabla X}-\nabla\hat{X}
\end{equation}
for all vector fields $X$ on $M$. To see that this equation has a solution $\alpha$ we may,
by $\R$-linearity, assume that $|X|=1$. Let $Y$ be such that $X,Y$ are a local orthonormal
frame of $TM$. Then $\nabla X=Y\sigma$ for some local real $1$-form $\sigma\in \Omega^{1}(\R)$.
Moreover, $\nabla\hat{X}$ anticommutes with $J$ so that
\[
\nabla\hat{X}=\hat{X}\beta+\hat{Y}\gamma
\]
for local real $1$-forms $\beta,\gamma\in\Omega^{1}(\R)$. Multiplying the latter by $\hat{X}$,
and observing that $(\nabla\hat{X}) \hat{X}$ is skew-adjoint, we conclude $\beta=0$.   
Therefore the right hand side of \eqref{eq:spin-nabla} is a multiple of $\hat{Y}$, and
we can solve for $\alpha$ uniquely.
\end{proof}
We are now ready to define the Dirac operator:
\begin{definition}\label{def:Dirac-operator}
Let $L$ be a Riemannian spin bundle with spin connection $\nabla$.
The  first order, linear differential operator
\[
\mathcal{D}\colon \Gamma(L)\to \Gamma(L)
\]
given by 
\[
\mathcal{D}=\hat{X}\nabla_{X}+\hat{Y}\nabla_{Y}\,,
\]
where $X,Y$ is a local orthonormal frame of $TM$, 
is called the {\em Dirac operator} of $L$.
\end{definition} 
It is easy to check that the Dirac operator $\mathcal{D}$ is quaternionic linear, 
self-adjoint with respect to the
hermitian inner product, and anticommutes with $J$. To interpret the
Dirac operator more conveniently as a complex holomorphic structure, 
consider the complex antilinear bundle isomorphism
\begin{equation}\label{eq:hat-map}
\hat{}\, : L\to \bar{K}L\,,\qquad \hat{\psi}=-\tfrac{1}{2}\hat{X}\psi\,.
\end{equation}
The spin connection $\nabla$ is a complex connection on $L$ and thus
$\delbar=\nabla''$ is a complex holomorphic structure on $L$. 
\begin{lemma}\label{lem:Dirac-delbar}
Let $L$ be a Riemannian spin bundle with spin connection $\nabla$.
Then the complex holomorphic structure $\delbar=\nabla''$
factors
\begin{equation}\label{eq:Dirac-delbar}
\delbar=\hat{\mathcal{D}}
\end{equation}
via the Dirac operator and the bundle isomorphism \eqref{eq:hat-map}.
\end{lemma}
\begin{proof}
Let $X,Y$ be an orthonormal local frame of $TM$. Then
\[
\hat{\mathcal{D}}_{X}
=-\tfrac{1}{2}\hat{X}(\hat{X}\nabla_{X}+\hat{Y}\nabla_{Y})=
\tfrac{1}{2}(\nabla+*J\nabla)_{X}
=\delbar_{X}\,.
\]
\end{proof}
Let us summarize what we have done so far: a Riemannian spin bundle $L$ carries a canonical
complex holomorphic structure $\delbar=\nabla''$ given by the 
$\bar{K}$-part of the spin connection which, up to the isomorphism \eqref{eq:hat-map},
is the Dirac operator $\mathcal{D}$. We can use the complex structure $J$ to decompose
\[
L=L_{+}\oplus L_{-}
\]
into the $\pm i$-eigenbundles of
$J$, which are complex line bundles isomorphic to each other. Clifford
multiplication is complex antilinear and thus interchanges $L_{\pm}$, and therefore the
Dirac operator interchanges $\Gamma(L_{\pm})$. On the other hand, the spin connection
$\nabla$ is complex and preserves the eigenbundles $L_{\pm}$ so that $\delbar$
induces isomorphic complex holomorphic structures on the bundles $L_{\pm}$.
The Riemann-Roch paired bundle $KL^{-1}$ carries a unique complex holomorphic structure
$\delbar$ such that
\[
d<\alpha,\psi>=<\delbar\alpha,\psi>-<\alpha\wedge\delbar\psi>
\]
holds, where $\alpha\in\Gamma(KL^{-1})$, $\psi\in\Gamma(L)$ and $<\,,>$ denotes the
Riemann-Roch pairing \eqref{eq:Riemann_Roch_pairing}. 
\begin{lemma}\label{lem:B-isomorphism}
Let $L$ be a Riemannian spin bundle. Then the bundle map
\[
B\colon L\to KL^{-1}\,,\qquad B\psi=(\hat\psi,-)\,,
\]
is a complex holomorphic isomorphism, that is to say,
$L$ is (the double of) a complex holomorphic spin bundle in the sense that $L_{\pm}^2=K$.

In particular, if $M$ is compact and of genus $g$ then
\[
\deg L=g-1\,.
\]
\end{lemma}
\begin{proof}
Since the map \eqref{eq:hat-map} and the  hermitian product in its first slot are
complex antilinear, $B$ is complex linear. Moreover, 
\[
*B\psi=(*\hat\psi,-)=(-J\hat{X}\psi,-)=(\hat{\psi},J(-))=JB\psi\,,
\]
where we used that $\hat{\psi}\in\bar{K}L$, so that $B$ is well-defined. 
Clearly, $B$ is non-zero and thus an isomorphism of complex quaternionic line
bundles. To show that $B$ is holomorphic, we first observe that \eqref{eq:Riemann_Roch_pairing}
implies that
\[
d^{\nabla}=\delbar_{KL^{-1}}\,.
\]
Then, using the definition of $B$, one calculates
\[
(d^{\nabla}B\psi)_{X,Y}=B\delbar_{X}\psi
\]
for $\psi\in\Gamma(L)$ and $X,Y$ a local orthonormal basis. Thus,
\[
\delbar B=B\delbar
\]
and $B$ is a holomorphic bundle isomorphism. 
\end{proof}
\subsection{Eigenvalue estimates} 
With this discussion in mind, consider the eigenvalue problem 
\begin{equation}\label{eq:Dirac_eigenvalue}
\mathcal{D}\psi=\lambda\psi
\end{equation}
for the Dirac operator $\mathcal{D}\colon \Gamma(L)\to\Gamma(L)$ on a Riemannian spin bundle
$L$ over a surface of genus $g$ with fixed metric. Since $\mathcal{D}$ is self-adjoint the eigenvalues $\lambda\in\R$ are real. 
The  {\em multiplicity} $m$ of an eigenvalue $\lambda$  is the {\em quaternionic} dimension
\[
m=\dim_{\H}\ker(\mathcal{D}-\lambda)
\]
of the kernel of $\mathcal{D}-\lambda$. Since  $\mathcal{D}$ anticommutes with $J$, the
eigenvalues come in pairs $\pm\lambda$ of equal multiplicity.
As an example, consider the $2$-sphere with its standard metric of curvature $1$. In this case
the eigenvalues for the Dirac operator are $\lambda=\pm n$, $n\in\N$, 
and occur with multiplicity $n$.

To apply our estimates regarding the Willmore energy to eigenvalue estimates, we slightly
rewrite the  eigenvalue equation \eqref{eq:Dirac_eigenvalue} by  
post-composing with the antilinear isomorphism
\eqref{eq:hat-map} to get
\[
\delbar\psi=\lambda \hat{\psi}\,.
\]
Note that if we put
\[
Q\psi:=-\hat\psi\,,
\]
then $Q\in\Gamma(\bar{K}\text{End}_{-}(L))$ and the eigenvalue equation becomes
\[
(\delbar+\lambda Q)\psi=0\,.
\] 
Now $\delbar+\lambda Q$ is a quaternionic holomorphic structure \eqref{eq:hol_decomp}
on the complex quaternionic line bundle $L$ and the last equation expresses
the holomorphicity of the section $\psi\in\Gamma(L)$. Thus, the multiplicity $m$ of 
the eigenvalue $\lambda$ of the Dirac operator
is precisely 
\[
m=\dim H^{0}(L)
\]
for the holomorphic structure $\delbar+\lambda Q$ on $L$.
To calculate its Willmore energy 
\[
W(L)=2\lambda^2\int <Q\wedge *Q>\,,
\]
note that the $2$-form $2<Q\wedge *Q>\in\Omega^2(\R)$ is just the area form on $M$:
if  $X,Y\in TM$ is an orthonormal basis then
\[
2<Q\wedge *Q>_{X,Y}=-4<(Q_X)^2>=-4<(\tfrac{1}{2}\hat{X})^2>=1\,
\]
and hence
\[
W(L)=\lambda^2\text{area}_M\,.
\]
We now apply Theorem~\ref{thm:Willmore_estimate} to the quaternionic holomorphic line bundle
$L$, which by Lemma~\ref{lem:B-isomorphism} has degree $d=g-1$, and obtain   
the following eigenvalue estimates for the Dirac operator:
\begin{theorem}\label{thm:Dirac_eignavalues}
Let $L$ be a Riemannian spin bundle over a compact surface of genus $g$ with fixed 
metric. If $\lambda$ is an eigenvalue for the Dirac operator $\mathcal{D}$ on $L$
of multiplicity $m$ then
\begin{equation*}
\lambda^2\text{\em area}_M\geq
\begin{cases}
4\pi m^2 &  \text{if}\qquad g=0\,,\\
\tfrac{\pi}{g}(m^2-g^2)&\text{if}\qquad g\geq 1\,.
\end{cases}
\end{equation*}
\end{theorem} 
From the example of the standard sphere, we see that the estimate is sharp in genus zero.
For multiplicity $m=1$ the genus zero estimate has been known \cite{Bar}. 
The estimate for higher multiplicities has been conjectured in \cite{Taimanov},
where it was shown to be true for metrics allowing a $1$-parameter family of symmetries.
In the higher genus case our estimates are new.

%%% Local Variables: 
%%% mode: latex
%%% TeX-master: "willmore"
%%% End: 

\section{Energy estimates for harmonic tori and area estimates
for constant mean curvature tori}\label{sec:harmonic} 
In this final section we apply the quaternionic Pl\"ucker formula to 
obtain estimates for the energy of harmonic 2-tori and the area of constant mean
curvature, CMC, tori. Rather then reinterpreting the usual descriptions of
CMC surfaces in the language of this paper, we 
develop the necessary concepts using quaternionic bundle theory.
This has the advantage of keeping the
exposition reasonably self contained and, at the same time, provides
an application of our quaternionic setup to surface geometry.
\subsection{Willmore connections}\label{subsec:willmore_connection}
The general framework we are working with is as follows:
we have a 
quaternionic vector bundle $V$ of rank $r$ with complex structure $S\in\Gamma(\text{End}(V))$
over a Riemann surface $M$. Given a flat connection
$\nabla$ on $V$, we split
\[
\nabla= \hat{\nabla}+\nabla_{-}
\]
into $S$-commuting and anticommuting parts \eqref{eq:nabla_decomposition}. Here,
$\hat{\nabla}$ is a complex connection on $V$, and $\nabla_{-}\in\Omega^1(\text{End}(V))$
is an endomorphism valued $1$-form. Decomposing further into type, we get 
\[
\hat{\nabla}=\del+\delbar\qquad \text{and}\qquad \nabla_{-}=A+Q\,.
\]
Since $\nabla$ is flat, \eqref{eq:dJ} implies the relations 
\begin{equation}\label{eq:dS}
\nabla S=2*(Q-A)\qquad \text{and}\qquad d^{\nabla}*A=d^{\nabla}*Q\,.
\end{equation}

The $\bar{K}$-part 
\[
\nabla''=\delbar+Q
\]
of the flat connection $\nabla$ defines
a quaternionic holomorphic structure on $V$. Its Willmore energy
\[
W(\nabla'')=2\int<Q\wedge *Q>
\]
defines a functional on the space
of flat connections on $V$ whose critical points we call
{\em Willmore connections}. This terminology is motivated by the fact that
Willmore surfaces give rise to such connections on rank two bundles \cite{bflpp}.

To compute the Euler-Lagrange equation, we vary the flat connection by
\[
\dot{\nabla}=\nabla B
\]
for some (compactly supported) section $B\in\Gamma(\text{End}(V))$. The form of the variation,
an infinitesimal gauge transformation,  guarantees
that $\nabla$ remains flat and that its holonomy does not change. 
Since $\dot{Q}=(\nabla B)''_{-}$ and the subbundles $\text{End}_{\pm}(V)$ are perpendicular 
with respect to
the trace form, we obtain 
\begin{align*}
\tfrac{1}{4}\dot{W}&=\int<\dot{Q}\wedge *Q>=\int<(\nabla B)_{-}'' \wedge *Q>=
\int<\nabla B\wedge *Q>\\
&=-\int< B, d^{\nabla}*Q>\,.
\end{align*}
Besides Stokes' Theorem, we also made use of the fact that 
$(\nabla B)_{-}'\wedge *Q=0$ due to type considerations. 
Together with \eqref{eq:dS} we therefore have
\begin{lemma}\label{lem:harmonic}
A flat connection $\nabla$ on $V$ is Willmore if and only if one of the 
following conditions holds:
\begin{enumerate}
\item $d^{\nabla}*Q=0$,
\item $d^{\nabla}*A=0$, or
\item $d^{\nabla}S*\nabla S=0$.
\end{enumerate}
\end{lemma}
Notice that the last equation expresses the fact that $S$ is a harmonic section of the flat
bundle $\text{End}(V)$ under the constraint that $S^2=-1$. This does not come as
a surprise since
any flat connection $\nabla$ on $V$ allows us to calculate the
{\em Dirichlet energy}
\[
E(S)=\tfrac{1}{2}\int<\nabla S\wedge *\nabla S>
\] 
of the complex structure $S\in\Gamma(\text{End}(V))$. Using the flatness of $\nabla$
and \eqref{eq:dS} we easily compute
\begin{equation}\label{eq:E=2W}
E(S)=2W(\nabla'')+4\pi\deg V\,.
\end{equation}
Any variation of $S$ can be realized by a variation of $\nabla$ via
a gauge transformation. This explains why the complex structure is harmonic
with respect to a Willmore connection and vice versa. 

For later reference we note that $d^{\nabla}*A=0$ can be seen as
a holomorphicity condition on $A$: if $V=W\oplus W$, where $W\subset V$
is the $i$-eigenbundle of $S$, then \eqref{eq:homs} gives
\[
\text{End}_{-}(V)=\text{Hom}_{\C}(\overline{W},W)\,.
\] 
The latter is a complex holomorphic vector bundle with holomorphic structure
induced from the complex connection $\hat{\nabla}$. Now
\[
0=d^{\nabla}*A=d^{\hat{\nabla}}*A+[Q\wedge *A]+[A\wedge *A]=Sd^{\hat{\nabla}}A\,,
\]
where we used $*A=SA$, $Q\wedge A=A\wedge Q=0$ by type considerations, and 
$[A\wedge *A]=0$
due to symmetry. Since $A$ is a section of the bundle $K\text{Hom}_{\C}(\overline{W},W)$
we have 
\begin{equation}\label{eq:A_hol}
\delbar A=d^{\hat{\nabla}}A=0\,,
\end{equation}
which means that $A$ is a complex holomorphic section. Similarly we get 
\[
\del Q=0\,,
\]
so that $Q$ is a complex antiholomorphic section of $\bar{K}\text{Hom}_{\C}(\overline{W},W)$ or,
equivalently, a complex holomorphic section of $K\text{Hom}_{\C}(W,\overline{W})$.
Therefore the compositions $AQ$ and $QA$ are holomorphic 
quadratic differentials with values in $\text{End}_{\C}(W)$.
 
Before continuing, we discuss the 
geometric content of the previous lemma.
The closed $1$-forms $2*A$ and $2*Q$ can be integrated, 
\begin{equation}\label{eq:f&g}
\nabla f=2*A\qquad\text{and}\qquad \nabla g=2*Q\,,
\end{equation}
to yield sections $f$ and $g$ of $\text{End}(\pi^{*}V)$ on the universal covering 
$\pi\colon \tilde{M}\to M$. From \eqref{eq:dS}
we get
\begin{equation}\label{eq:cmc}
d^{\nabla}*\nabla f= d^{\nabla}2S*A=2\nabla S\wedge *A=-\nabla f\wedge\nabla f\,,
\end{equation}
where we again used that $Q\wedge A=0$. Note that \eqref{eq:cmc}
formally looks like the CMC one condition for a surface $f$ in $\R^3$.
%%% considerations explained in section~\ref{subsec:prelim}.%%%
Finally, adjusting constants of integration, we also have   
\[
g=f+\pi^{*}S
\]
\begin{comment}
Note that the first equation,
\begin{equation}\label{eq:cmc}
d^{\nabla}*\nabla f+\nabla f\wedge\nabla f=0\,,
\end{equation}
formally looks exactly like the CMC one condition for a surface $f$ in $\R^3$, and the
second equation says that 
\end{comment}
which, continuing our analogy, says that $g$ is the parallel CMC surface to $f$ in unit normal
distance. 
\subsection{Harmonic maps into the 2-sphere and CMC surfaces}\label{subsec:CMC}
In case the rank
of $V$ is one, this formal analogue becomes precise and
gives the explicit correspondence between CMC surfaces and
harmonic maps into the $2$-sphere. Let $V$ be a quaternionic line bundle
with complex structure $S$, and let $\nabla$ be a Willmore connection. On the universal
cover, $\nabla$ trivializes and we fix a parallel section $\phi$ of $\pi^{*}V$. This allows
us to identify sections of $\pi^{*}V$ and $\text{End}(\pi^{*}V)$ with $\H$ valued maps. 
Therefore the complex structure $S$ is given by the harmonic map 
\[
N\colon \tilde{M}\to S^2\subset \text{Im}\H\,,
\]
where $S\phi=\phi N$. If 
\[
H\colon \pi_1(M)\to \text{\bf Gl}(1,\H)=\H\setminus\{0\}
\]
is the holonomy representation of $\nabla$, then 
\begin{equation}\label{eq:N_holonomy}
\gamma^{*}N=H(\gamma)^{-1}N H(\gamma)\,.
\end{equation}
Conversely, given a harmonic map $N\colon \tilde{M}\to S^2$ satisfying \eqref{eq:N_holonomy}
for some representation $H\colon \pi_1(M)\to \text{\bf Gl}(1,\H)$, we obtain a flat 
quaternionic line bundle $V$. The harmonic map $N$ induces a complex structure
$S$ on $V$ and the flat connection is Willmore. Therefore, up to holonomy, Willmore connections on rank one bundles
are the same as harmonic maps into the 2-sphere.

For our geometric applications it is convenient to work
with hermitian connections on $V$. We may add any
closed, real $1$-form $\omega\in\Omega^{1}(\R)$ to $\nabla$ and 
obtain another Willmore connection $\tilde{\nabla}=\nabla+\omega$:
\[
\tilde{Q}=Q\qquad\text{and}\qquad d^{\tilde{\nabla}}*Q=0\,.
\]
Using this freedom to change $\nabla$, we
ensure that the real line bundle of quaternionic hermitian forms on $V$ 
has trivial holonomy. Thus we can assume that $V$ has a parallel quaternionic hermitian form 
\begin{equation}\label{eq:quat_herm_form}
<\,,\,>\colon V\times V\to \H\,,
\end{equation}
so that the holonomy 
\begin{equation}\label{eq:holonomy_rep}
H\colon \pi_1(M)\to S^3\subset \H
\end{equation}
takes values in the unitary quaternions.

Since $V$ has rank one, $S$ is automatically skew hermitian with respect to $<\,,\,>$. 
Then \eqref{eq:dS} implies that $A$ and $Q$ are skew hermitian, which means
that $<\,,\,>$ is also parallel with respect to the complex connection $\hat{\nabla}$.
If $W\subset V$ is the $i$-eigenbundle
of $S$, then one easily checks that $<\,,\,>$ takes complex values when restricted
to $W$. Therefore $W\subset V$ is a complex hermitian line bundle with hermitian
connection $\hat{\nabla}$ and complex holomorphic structure $\delbar$. As we have seen in
\eqref{eq:A_hol},
the $1$-forms $A$ and $Q$ are holomorphic sections
of $KW^2$ and $KW^{-2}$, where we now identify $W^{-1}=\bar{W}$ via the hermitian form.
Since $\text{End}_{+}(V)=\text{End}_{\C}(W)$ is the trivial bundle, $AQ$ is a
holomorphic quadratic differential on $M$. If $AQ=0$ then either $A$ or $Q$ vanishes identically,
in which case $S$ is, up to holonomy, a holomorphic or antiholomorphic map
into the $2$-sphere. If $AQ\neq 0$ then neither $A$ nor $Q$ vanish identically and
both bundles $KW^2$ and $KW^{-2}$ have non-negative degree, assuming that 
$M$ is compact of genus $g$. Since the degree \eqref{eq:degree} of the quaternionic bundle 
$V$ is $\deg V=\deg W$,
we get  
\begin{equation}\label{eq:deg_L}
|\deg V|\leq g-1\,.
\end{equation}
Put into harmonic maps language, this last relation rephrases the well known result \cite{Eells}
that a harmonic map $N\colon M\to S^2$ of a compact Riemann surface of genus $g$ 
into the $2$-sphere, whose degree
satisfies $|\deg N|\geq g$, is either holomorphic or antiholomorphic.
We will assume from now on that neither $A$ nor $Q$ vanish identically.

To obtain CMC surfaces we have to integrate \eqref{eq:f&g} the closed $1$-forms
$2*A$ and $2*Q$ that, under our identifications, are given by the
$\text{Im}\H=\R^3$ valued forms
\begin{equation}\label{eq:A&Q}
2*A=\tfrac{1}{2}(N*dN-dN)\qquad\text{and}\qquad 2*Q=\tfrac{1}{2}(N*dN+dN)\,.
\end{equation}
If
\[
f,g\colon \tilde{M}\to\R^3
\]
are those integrals then
\begin{equation}\label{eq:f_holonomy}
\gamma^{*}f=H(\gamma)^{-1}fH(\gamma)+T(\gamma)\,,
\end{equation}
where 
\[
T\colon \pi_1(M)\to\R^3
\]
are the translational periods of $f$. Adjusting constants of integration, we also have 
\[
g=f+N\,.
\]
From $*A=SA$ we obtain
\[
*df=Ndf\,,
\]
which says that $f$, away from the zeros of $df$, is a conformal immersion with unit normal
map $N$. But $df=2*A$ so that the zeros of $df$ coincide with the zeros of $A$. We have
seen above that $A$ is a holomorphic section of $KW^2$. Therefore $f$ is an immersion 
if and only if $A$ has no zeros, in which case $KW^2$ is holomorphically trivial, that
is to say $W^{-1}$ is a complex holomorphic spin bundle on $M$. Note that due to 
\eqref{eq:deg_L} this happens if and only if $\deg V=1-g$. Since we are interested in 
surfaces without branch points, we assume from now on 
\begin{equation}\label{eq:deg_eq}
\deg V=1-g\,.
\end{equation}
From \eqref{eq:cmc} we know that $f$ satisfies    
\[ 
d*df+df\wedge df=0\,,
\]
and therefore is an immersion of constant mean curvature one on the universal covering.
To obtain a CMC immersion on $M$, two conditions \eqref{eq:f_holonomy} have to be satisfied:
the flat connection $\nabla$ has to be trivial, i.e., the unit normal map
$N$ has to be defined on $M$, and the translational
periods $T(\gamma)$ have to vanish. In this situation, the line bundle $V$ with the 
quaternionic holomorphic structure $\nabla''$ has the two independent holomorphic sections 
$\phi$ and $\psi=\phi f$ so that $h^{0}(V)\geq 2$. From the explicit expression 
\eqref{eq:A&Q} we calculate the Willmore energy   
\begin{equation}\label{eq:W&area}
W(\nabla'')=2\int Q\wedge *Q=\int H^2-K=\text{area}(f)-4\pi(1-g)
\end{equation}
of the holomorphic structure $\nabla''$, where we used that the mean curvature $H=1$.
Together with \eqref{eq:E=2W} we therefore  get
\begin{equation}\label{eq:E&area}
E(N)=2\,\text{area}(f)+4\pi(g-1)\,.
\end{equation}

Conversely, given a conformal CMC immersion $f\colon M\to\R^3$ with unit normal
map $N\colon M\to\ S^2$, we let $V=M\times \H$ with the trivial connection $\nabla$ and
complex structure $S$ given by $N$. The hermitian product $<\,,\,>$ is the standard one on $\H$.
Since $f$ is conformal we have $*df=Ndf$, and any other $\R^3$ valued $1$-form $\omega$
satisfying $*\omega=N\omega$ is a real multiple of $df$. We now decompose
the second fundamental form
\[
dN=dN'+dN''=\tfrac{1}{2}(dN-N*dN)+\tfrac{1}{2}(dN+N*dN)
\]
into type with respect to the complex structure given by $N$. Since the $\bar{K}$-part
$dN''$ is trace free, we get 
\[
dN'=-Hdf=-df\,.
\]
Comparing with \eqref{eq:A&Q} we thus see that $df=2*A$. The
CMC one condition \eqref{eq:cmc} implies that $d^{\nabla}*A=0$. Thus $S$ is harmonic
or, equivalently, the trivial connection is Willmore by Lemma~\ref{lem:harmonic},
which is the Ruh-Vilms Theorem \cite{Ruh-Vilms}.  
We finally note that the bundle $V$ with the holomorphic structure $\nabla''$ is the 
pullback of the dual tautological bundle over $\H\P^1$ by the conformal immersion
$f\colon M\to \R^3\subset \H\subset \H\P^1$ via the Kodaira correspondence 
discussed in section~\ref{subsec:Kodaira}. The $2$-dimensional linear system in
$H^{0}(V)$ induced by $f$ is spanned by 
the restrictions of the linear coordinate projections
$\psi=\alpha_{|V^{-1}}$ and $\phi=\beta_{|V^{-1}}$
to the pullback $V^{-1}\subset M\times\H^2$ of the tautological bundle over $\H\P^1$,
and $\psi=\phi f$. 
\begin{remark}\label{rem:rank_2}
We have focused on Willmore connections in the rank one case, i.e., harmonic maps
into the 2-sphere and CMC surfaces, since this is the setting for which 
we will give estimates on the Willmore energy. However, it seems worthwhile to
remark briefly on the rank two case, which includes the theory of
Willmore surfaces in 4-space. If $V$ has rank two and $\nabla$ is Willmore,
then the complex structure $S$ is, up to holonomy, a harmonic map into the 
space of oriented 2-spheres in $\H\P^1=S^4$ as explained in Example~\ref{ex:4.2}. 
For such a map to be the mean curvature sphere congruence
of a surface in 4-space is equivalent to $QA=0$. Recall that $A$ and $Q$
are complex holomorphic sections \eqref{eq:A_hol}
and that $QA$ is a holomorphic quadratic differential. In the generic case
neither $A$ nor $Q$ vanish identically. Otherwise $S$ is a holomorphic or
antiholomorphic map and we need $\nabla S\in\Omega^1(\text{End}(L))$ to have rank one.
In both cases there is 
a smooth quaternionic line subbundle $L\subset V$ that satisfies
\[
\text{im} A\subset L\subset \text{ker}Q\,.
\]
Now $\nabla$ is a Willmore connection so 
that $d^{\nabla}*A= d^{\nabla}*Q=0$ by Lemma~\ref{lem:harmonic}. If we interpret $L\subset V$
as a map into $\H\P^1$ with M\"obius holonomy, then these equations express
the fact \cite{bflpp} that $L$ is a Willmore surface with harmonic mean curvature
sphere congruence $S$. In case $A$ or $Q$
vanish identically, the Willmore surface arises via the twistor projection
from a holomorphic curve into $\C\P^3$. 

Conversely, if $f\colon M\to \H\P^1$ is a Willmore surface, we view $f$ as a line subbundle 
$L\subset V$
of the trivial bundle $V=M\times\H^2$. The mean curvature sphere congruence $S$ 
is then a complex structure on $V$, and the trivial connection $\nabla$ is Willmore.
Since $S$ is the mean curvature sphere, Example~\ref{ex:4.2},
Theorem~\ref{thm:mcS} and Corollary~\ref{cor:adapted_nabla_relations} imply
 $Q_{|L}=0$ and  $\text{im}A\subset L$ so that $QA=0$. For a more detailed discussion
of the geometry of the rank two case we refer the reader to \cite{bflpp}.
\end{remark}

\subsection{The spectral curve}
To estimate the Willmore energy of a Willmore connection $\nabla$ on a rank one
vector bundle $V$, we
use the Pl\"ucker formula in Theorem~\ref{thm:Willmore_estimate}. For this to 
succeed, we need to know how many holomorphic sections $h^{0}(V)$ the quaternionic
holomorphic structure $\nabla''$ admits. In case the underlying Riemann
surface $M$ has genus one, integrable system methods allow us to estimate
$h^{0}(V)$ and, as a consequence, also the Willmore energy. 

The basic ingredient in all of this is a certain family of flat connections
associated to a Willmore connection.  
\begin{lemma}\label{lem:lambda_family}
Let $V$ be a quaternionic vector bundle with complex structure $S$ and
flat connection $\nabla$. Consider the family of connections 
\[
\nabla_{\lambda}=\nabla+(\lambda-1)A\,,
\]
where $A=\nabla_{-}'$ and $\lambda=x+yS$ is a complex parameter.

Then $\nabla=\nabla_1$ is Willmore if and only if $\nabla_{\lambda}$ is a flat connection
on $V$ for all unitary $\lambda\in S^1$. If this is the case, $\nabla_{\lambda}$
is also Willmore for all $\lambda\in S^1$. 
\end{lemma}
\begin{proof}
The curvature $R_{\lambda}$ of the connection $\nabla_{\lambda}$ is given by
\begin{align*}
R_{\lambda}=&R+d^{\nabla}(\lambda-1)A+(\lambda-1)A\wedge (\lambda-1)A\\
=& (x-1)\,d^{\nabla}A+y\,d^{\nabla}*A+((x-1)^2+y^2)A\wedge A\,,
\end{align*}
where we used that the connection $\nabla$ is flat. From \eqref{eq:dS}
we have $\nabla S= 2*(Q-A)$ and therefore
\[
d^{\nabla}A=-d^{\nabla}S*A= 2\,A\wedge A-S\,d^{\nabla}*A\,.
\]
Inserting this into the expression for $R_{\lambda}$ we obtain
\[
R_{\lambda}= (|\lambda|^2-1)A\wedge A+(1-\lambda)S\,d^{\nabla}*A\,.
\]
The first term takes values in $\text{End}_{+}(V)$, whereas the second term
takes values in $\text{End}_{-}(V)$. Together with Lemma~\ref{lem:harmonic}, this 
shows that $\nabla$ is Willmore if and only if $\nabla_{\lambda}$ is flat for
unitary $\lambda\in S^1$. 

Finally $Q_{\lambda}=(\nabla_{\lambda}'')_{-}=Q$ so that 
\[
d^{\nabla_{\lambda}}*Q_{\lambda}=d^{\nabla}*Q+(\lambda-1)A\wedge Q + 
Q\wedge (\lambda-1)A=0\,,
\]
where we used that $d^{\nabla}*Q=0$ and that $A\wedge Q=Q\wedge A=0$. 
This shows that also $\nabla_{\lambda}$ is Willmore.
\end{proof}
Recall that the complex structure $S$ on $V$ is harmonic with respect
to a Willmore connection $\nabla$. To interpret $S$ as a harmonic map into $\text{\bf Gl}(r,\H)$
with the constraint $S^2=-1$, we have to trivialize the connection $\nabla$ on the universal covering.
Thus, trivializing the family $\nabla_{\lambda}$ gives an $S^1$ family of harmonic
maps. The correspondence between harmonic maps into symmetric spaces
and $S^1$ families of flat connections is folklore, and the above lemma
is a manifestation of this fact.

Note that the family of Willmore connections $\nabla_{\lambda}$ all induce the same holomorphic
structure $\nabla_{\lambda}''=\nabla''$. Therefore, any parallel section $\psi$ of
$\nabla_{\lambda}$ is automatically a holomorphic section, i.e., $\nabla''\psi=0$. 
Generically though, $\nabla_{\lambda}$ does not admit parallel 
sections for unitary $\lambda\in S^1$. To remedy this, we extend $\nabla_{\lambda}$
to a flat family of complex connections parametrized by $\C\setminus\{0\}$.
Let 
\[
I\colon V\to V\,,\qquad I\psi=\psi i
\]
denote the quaternionic multiplication by $i$. Then $I$ is a complex structure on
$V$, no longer quaternionic linear, which commutes with $S$, $A$, $Q$ and $\nabla$.
We can therefore view $V$ via $I$ as a rank $2r$ complex vector bundle with a flat
complex connection $\nabla$. Consider the family
\[
\nabla_{\lambda}=\nabla+(\lambda-1)A\,,
\]
where, in contrast to Lemma~\ref{lem:lambda_family}, the parameter $\lambda=x+yS$ is
now given by $I$-complex numbers
\[
x=\tfrac{1}{2}(\mu+\mu^{-1})\,,\qquad y=\tfrac{1}{2}(\mu^{-1}-\mu)I\,,\qquad 
\mu=a+bI\in\C\setminus\{0\} 
\]
still satisfying $x^2+y^2=1$. Expressing $\nabla_{\lambda}$ in the complex parameter $\mu$,
we therefore obtain the family
\begin{equation}\label{eq:mu_family}
\nabla_{\mu}=\nabla+\tfrac{1}{2}(\mu+\mu^{-1}-2)A+\tfrac{1}{2}(\mu^{-1}-\mu)I*A
\end{equation}
of complex connections on $V$ parameterized over $\C\setminus\{0\}$. Note that 
for unitary $\mu=a+bI\in S^1$, the connection $\nabla_{\mu}$ restricts to the
connection $\nabla_{\lambda}$ above, where $\lambda=a+bS\in S^1$.
Since $I$ is parallel and commutes with $A$ and $Q$, 
the same formal calculation as in the proof of Lemma~\ref{lem:lambda_family}
shows that $\nabla$ is Willmore if and only if $\nabla_{\mu}$ is flat
for all $\mu \in \C\setminus\{0\}$.

The family of complex connections $\nabla_{\mu}$ possesses a symmetry 
induced by the quaternionic multiplication
\[
J\colon V\to V\,,\qquad J\psi=\psi j\,.
\]
Because  $JI=-IJ$ we have $J\mu=\bar{\mu}J$, and therefore
\begin{align}\label{eq:J_symmetry}
J\nabla_{\mu}&=J\nabla+J\tfrac{1}{2}(\mu+\mu^{-1}-2)A+J\tfrac{1}{2}(\mu^{-1}-\mu)I*A\notag\\
&=\nabla+\tfrac{1}{2}(\bar{\mu}+\bar{\mu}^{-1}-2)AJ-\tfrac{1}{2}(\bar{\mu}^{-1}-\bar{\mu})I*AJ\\
&=\nabla_{\tfrac{1}{\bar{\mu}}}J\notag\,,
\end{align}
where we used that $\nabla$ and $A$ are quaternionic linear, i.e., commute with $J$.
Calculating the $\bar{K}$ part of $\nabla_{\mu}$ with respect to $S$ shows again that 
the holomorphic structure $\nabla_{\mu}''=\nabla''$ is independent
of $\mu$. Thus, parallel sections of $\nabla_{\mu}$ yield holomorphic sections
with respect to $\nabla''$. From the symmetry \eqref{eq:J_symmetry}, we see that 
if $\psi$ is parallel for $\nabla_{\mu}$, then $J\psi$ is parallel for 
$\nabla_{\tfrac{1}{\bar{\mu}}}$. Therefore, the holomorphic sections 
of $V$ obtained from parallel sections of $\nabla_{\mu}$ for all  $\mu \in \C\setminus\{0\}$
span a quaternionic linear subspace of $H^{0}(V)$. 
\begin{lemma}\label{lem:mu_family}
Let $V$ be a quaternionic vector bundle with complex structure $S$ and
flat connection $\nabla$. Then
\begin{enumerate}
\item
$\nabla$ is Willmore if and only if the family $\nabla_{\mu}$ is flat
for every $\mu \in \C\setminus\{0\}$.
\item
The holomorphic structure $\nabla_{\mu}''=\nabla''$ is independent
of $\mu$, so that parallel sections of $\nabla_{\mu}$ yield holomorphic sections
with respect to $\nabla''$.
\item
The linear subspace $U$ spanned by 
\[
\{\psi\in \Gamma(V)\,;\,\nabla_{\mu}\psi=0\,,\,\mu \in \C\setminus\{0\}\}\subset H^{0}(V)
\]
is a quaternionic subspace and thus a linear system in $H^{0}(V)$.
\end{enumerate}
\end{lemma}
In order to estimate the dimension of the linear system $U\subset H^{0}(V)$
in the case when $V$ is a quaternionic line bundle, we
first show that parallel sections of $\nabla_{\mu}$ to distinct $\mu$ are
linearly independent over $\C$. 
\begin{lemma}\label{lem:independent}
Let $V$ be a quaternionic line bundle with complex structure $S$, and let 
$\nabla$ be a Willmore connection such that both $A=\nabla_{-}'$ and
$Q=\nabla_{-}''$ are not identically
zero. 

If $\mu_k\in\C\setminus\{0\}$ are distinct and $\psi_k\in\Gamma(V)$ nontrivial
parallel sections of $\nabla_{\mu_k}$, then $\psi_k$ are complex linearly 
independent as sections of the complex rank two bundle $V$ with complex
structure $I$.

In particular, since $H^0(V)$ is finite dimensional over a
compact Riemann surface, there are at most finitely many
$\mu\in\C\setminus\{0\}$ so that $\nabla_{\mu}$ is trivial. Put differently,
if the family $\nabla_{\mu}$ is trivial, then either $Q=0$ or $A=0$ in which case 
$S$ is an antiholomorphic or holomorphic map into the $2$-sphere. 
\end{lemma}
\begin{proof}
We rephrase the basic fact from linear algebra that eigenvectors 
corresponding to distinct eigenvalues
are linearly independent. Denote by
\[
\nabla^{(1,0)}=\tfrac{1}{2}(\nabla-I*\nabla)
\]
the $K$-part of the connection $\nabla$ with respect to the complex structure $I$ on $V$.
From \eqref{eq:mu_family} we easily calculate that
\[
\nabla^{(1,0)}_{\mu}=\nabla^{(1,0)}+\tfrac{1}{2}(\mu-1)(1-IS)A=:\del+\mu\omega\,,
\]
where the antiholomorphic structure $\del$ collects the $\mu$-independent terms and
$\omega=\tfrac{1}{2}(1-IS)A$ is an $\text{End}(V)$ valued $(1,0)$-form. Note that
parallel sections of $\nabla_{\mu}$ are also in the kernel of $\nabla^{(1,0)}_{\mu}$.
Assume that $\psi_k\in \Gamma(V)$, $k=1,\dots,n$, are linearly independent 
parallel sections of $\nabla_{\mu_k}$ for distinct $\mu_k$, so that 
\[
\del\psi_k+\mu_k\omega\psi_k=0\,.
\]
Let $\psi_0=\sum_{k}c_k \psi_k$, $c_k\in \C$, be
a parallel section of $\nabla_{\mu_0}$ with $\mu_0$ distinct from the $\mu_k$. 
We have to show that all $c_k=0$. Since
\[
0=(\del+\mu_0\omega)\psi_0=\sum_{k=1}^{n}(\mu_0-\mu_k)c_k\omega\psi_k
\]
and $\mu_0-\mu_k\neq 0$, it suffices to show that $\omega\psi_k$ are linearly independent.
Because parallel sections are holomorphic, this will follow if we can show that 
\[
\omega\colon  H^{0}(V)\to \Gamma(KV)
\]
is injective. Away from the isolated zeros \eqref{eq:A_hol} of $A$, 
the bundle map $\omega$ has 
kernel equal to the
$i$-eigenbundle $W\subset V$ of $S$. 
Therefore the kernel of $\omega$ on $H^{0}(V)$
consists of sections $\psi\in \Gamma(W)$ for which 
\[
0=\nabla''\psi=\delbar\psi+Q\psi\,.
\]
But $\delbar$ preserves sections of $W$ whereas $Q$ maps them to sections of $Wj$,
the $-i$-eigenbundle of $S$. Since $Q$ also vanishes at most at isolated points,
we conclude that $\psi=0$. Note that these last arguments used the fact that we are working
on a quaternionic line bundle.
\end{proof}
Recall that in the rank one case we
may assume that $V$ 
carries a parallel quaternionic hermitian
form \eqref{eq:quat_herm_form} with respect to the given Willmore connection $\nabla$. 
Decomposing this form,
\begin{equation}\label{eq:form_decomposition}
<\,,\,>=(\,,\,)+j\det\,,
\end{equation}
we obtain the hermitian form $(\,,\,)$ and the determinant form $\det$ 
on the complex rank two bundle $V$ with complex structure $I$. As we have
observed before, $S$ and therefore $A$ and $SA$ are skew hermitian with respect to 
$<\,,\,>$. This implies that $<\,,\,>$, and hence also $(\,,\,)$ and $\det$, 
are parallel for $\nabla_{\mu}$ as long as
$\mu\in S^1$ is unitary. Moreover, $\det$ is parallel with respect to $\nabla_{\mu}$ for
all $\mu\in\C\setminus\{0\}$ so that $\nabla_{\mu}$ 
is a family of $\text{\bf Sl}(2,\C)$ connections, which is special unitary
along the unit circle $\mu\in S^1$. Fixing a base point on $M$, we thus get a family of holonomy
representations
\[
H_{\mu}\colon \pi_{1}(M)\to\text{\bf Sl}(2,\C)\,,\qquad \mu\in\C\setminus\{0\}\,,
\]
which restrict to special unitary representations
\[
H_{\mu}\colon \pi_{1}(M)\to\text{\bf SU}(2)\,,\qquad \mu\in S^1
\]
along the unit circle. From \eqref{eq:J_symmetry} we obtain the symmetry
\begin{equation}\label{eq:H_symmetry}
H_{\tfrac{1}{\bar{\mu}}}=JH_{\mu}J^{-1}\,,
\end{equation}
expressing the fact that $H_{\mu}$ is quaternionic linear for unitary $\mu\in S^1$.
Since $\nabla_{\mu}$ depends holomorphically on $\mu\in\C\setminus\{0\}$, the holonomy
$H_{\mu}$ varies holomorphically in $\mu$.

We now address how to find values for $\mu$ so that the connection
$\nabla_{\mu}$ admits a parallel section. In other words, we have to find
values for $\mu$ so that the holonomy representation $H_{\mu}$ of $\nabla_{\mu}$ has
a common eigenvector with the eigenvalue one. For the existence of common
eigenvectors one generally needs commuting matrices. This is certainly 
guaranteed if the fundamental group $\pi_1(M)=\Z^2$ is abelian, i.e., if
the underlying compact Riemann surface $M=T^2$ is a 2-torus, which we will
assume from now on.

First we 
exclude the case where the holonomy $H_{\mu}$ is trivial which, by Lemma~\ref{lem:independent},
corresponds to the case of holomorphic respectively antiholomorphic maps 
from the torus $T^2$ into $S^2$. Equivalently \eqref{eq:deg_L}, we assume that 
\begin{equation}\label{eq:deg=0}
\deg V=0\,.
\end{equation}
Since the holonomy is nontrivial and depends holomorphically on
the parameter $\mu$, we have two distinct common eigenlines of $H_{\mu}$ for
generic values of $\mu\in\C\setminus\{0\}$. Off this generic set, $H_{\mu}(\gamma)$
has coinciding eigenvalues for every cycle $\gamma\in\pi_1(M)$, or
equivalently, for two generating cycles: otherwise there would be a
cycle $\gamma$ for which $H_{\mu}(\gamma)$ had distinct eigenvalues, and hence
the representation $H_{\mu}$ would have two distinct common eigenlines.  
Since $H_{\mu}$ is a $\text{\bf Sl}(2,\C)$ representation, having coinciding eigenvalues
is equivalent to 
\begin{equation}\label{eq:H_mu}
H_{\mu}=\pm\begin{pmatrix} 1&*\\ 0&1\end{pmatrix}\,.
\end{equation}
Therefore, on at most a $4$-fold covering of $T^2$, the values $\mu$ for which $H_{\mu}$
have coinciding eigenvalues yield parallel sections of $\nabla_{\mu}$. From
Lemma~\ref{lem:independent} we conclude that there can be at most finitely
many such values. Let 
\[
B\subset \C\setminus\{0\}
\]
be the finite set of $\mu$ satisfying the following conditions:
\begin{enumerate}
\item $H_{\mu}(\gamma)$ has
coinciding eigenvalues for two generating cycles $\gamma\in\pi_1(M)$, and therefore all cycles,
\item $H_{\mu}$ does not have two common distinct eigenlines, and
\item the holomorphic map $\trace H(\gamma)-1$ has odd order zeros at $\mu$ for all
cycles $\gamma\in\pi_1(M)$.
\end{enumerate}
Due to \eqref{eq:H_symmetry} the set $B$ is invariant under
the anti-involution $\mu\mapsto \bar{\mu}^{-1}$. Moreover, no $\mu\in B$ can    
be on the unit circle $S^1$, otherwise $H_{\mu}$ would be the limit
of $\text{\bf SU}(2)$ representations having perpendicular eigenlines and therefore
itself have perpendicular eigenlines. Thus 
\[
B=\{\mu_1,\dots,\mu_g,\bar{\mu_1}^{-1}\dots,\bar{\mu_g}^{-1}\}
\]
is a collection of $g$ pairs of points in $\C\setminus\{0\}$ 
symmetric with respect to the unit circle. 
\begin{definition}\label{def:spectral_curve}
Let $V$ be a quaternionic line bundle with complex structure
$S$ of degree zero over a $2$-torus and let $\nabla$ be a Willmore connection. 
The {\em spectral curve} of $\nabla$ is 
defined to be the hyperelliptic Riemann surface 
\[
\mu\colon \Sigma\to\C\P^1
\]
branched over the values $\mu\in B\cup\{0,\infty\}$, so that the genus of
$\Sigma$ is $g$. We call $g$ the {\em spectral genus}
of the Willmore connection $\nabla$. The same terminology is used for the
harmonic section $S$ and the resulting CMC surface.
\end{definition}
From the definition of the spectral curve, one easily sees that the eigenlines of $H_{\mu}$ 
define a holomorphic line subbundle, the {\em eigenline bundle} $\mathcal{E}$,
of the trivial $\C^2$ bundle 
over the punctured curve $\Sigma\setminus\{\mu^{-1}(0),\mu^{-1}(\infty)\}$.
For $x\in\Sigma\setminus\{\mu^{-1}(0),\mu^{-1}(\infty)\}$ 
the common eigenlines of $H_{\mu(x)}$ are given by
$\mathcal{E}_x$ and $\mathcal{E}_{\sigma(x)}$ where $\sigma$ denotes 
the hyperelliptic involution on $\Sigma$. It can  be shown \cite{Hitchin} that 
the eigenline bundle extends holomorphically into the punctures over $\mu=0$ and
$\mu=\infty$. 

Note that all the constructions so far were done with respect to a fixed base point on $T^2$.
If we change to another point, the holonomy representation $H_{\mu}$ gets conjugated
by parallel transport with respect to $\nabla_{\mu}$, 
so that the spectral curve $\Sigma$ remains unchanged.
What does change, though, is the holomorphic class of the eigenline bundle $\mathcal{E}$ which
can be shown to move linearly with the base point on $T^2$.
For more details and how to reconstruct harmonic $2$-tori from
algebraic geometric data, we refer the reader to Hitchin's paper \cite{Hitchin},
which we have followed in spirit during the above discussion.
\subsection{Energy and area estimates}
Let $V$ be a quaternionic line bundle of degree zero with 
complex structure $S$ over a $2$-torus $T^2$ and let $\nabla$ be a Willmore connection of
spectral genus $g$. Then each branch value $\mu\neq 0,\infty$ of the 
spectral curve $\Sigma$ gives a parallel section for $\nabla_{\mu}$ of $V$ on a 
$4$-fold cover of $T^2$. Moreover, Lemma~\ref{lem:mu_family} and
Lemma~\ref{lem:independent} imply that these sections
are complex linearly independent, holomorphic with respect to $\nabla''$ and
span a linear system  $U\subset H^{0}(V)$. Since the
complex dimension of the linear system is $2g$, its quaternionic dimension is $g$, and
we have 
\[
h^{0}(V)\geq g\,.
\] 
If in addition the Willmore connection $\nabla=\nabla_{1}$ is trivial, i.e., if 
the bundle $V$ is induced by a harmonic map $N\colon T^2\to S^2$, then we obtain one
more independent parallel, and hence holomorphic, section so that 
\[
h^{0}(V)\geq g+1\,.
\]
Finally, if $\nabla$ comes from a CMC torus $f\colon T^2\to\R^3$ then,
as we have seen in Section~\ref{subsec:CMC}, there is an additional
independent holomorphic section yielding
\[ 
h^{0}(V)\geq g+2\,
\]
on a $4$-fold cover of $T^2$. Since $\deg V=0$, the relation \eqref{eq:E=2W} implies
that the energy of the harmonic section $S$ is twice the Willmore energy
\[
E(S)=2W(\nabla'')
\]
of the holomorphic structure $\nabla''$. Moreover, from \eqref{eq:W&area} we get
\[
W(\nabla'')=\text{area}(f)
\]
for the CMC one torus $f$ corresponding to the Willmore connection $\nabla$.

We now apply the Pl\"ucker formula \eqref{eq:special_W_estimates} 
for the special case of a line bundle over a torus $T^2$. Keeping in mind that
we work on a $4$-fold covering of $T^2$ we obtain the following
energy and area estimates in terms of the spectral genus:
\begin{theorem}\label{thm:energy&area_estimates}
\begin{enumerate}
\item
If $N\colon T^2\to S^2$ is a harmonic map of degree zero and spectral genus $g$, then
its energy satisfies the estimate
\[
E(N)\geq 
\begin{cases}
\tfrac{\pi}{2}(g+1)^2 & \,\,\text{if $g$ is odd, and}\\
\tfrac{\pi}{2}((g+1)^2-1)&\,\,\text{if $g$ is even}\,. 
\end{cases}
\]
\item If $f\colon T^2\to \R^3$ is a CMC one torus of spectral genus $g$, then its area
satisfies the estimate
\[
\text{\em area}(f)\geq
\begin{cases}
\tfrac{\pi}{4}(g+2)^2 & \,\,\text{if $g$ is even, and}\\
\tfrac{\pi}{4}((g+2)^2-1)&\,\,\text{if $g$ is odd}\,. 
\end{cases}
\]
\end{enumerate}
\end{theorem}
We conclude by discussing various applications of the above estimates. 
A homomorphism from a $2$-torus $T^2$ into an equator $S^1\subset S^2$ 
is a harmonic map whose energy is proportional to the area of $T^2$. 
This shows that the energy of such harmonic $2$-tori can be arbitrarily small by
choosing a thin fundamental domain. It is therefore natural to 
determine the energy value below which a harmonic $2$-torus is a homomorphism
into $S^1$. The unit normals to the embedded Delauney surfaces 
of constant mean curvature $H=1$ give harmonic maps of $T^2$ into $S^2$. 
Their energies attain 
the minimal value $2\pi^2$ on the round cylinder whose unit normal map is
a homomorphism. This suggests the conjecture that 
a harmonic torus $N\colon T^2\to S^2$ with energy satisfying $E(N)<2\pi^2$ has to be 
a homomorphism into $S^1$. A slight refinement of the estimates in 
Theorem~\ref{thm:energy&area_estimates}
provides some support for this conjecture, even though our result
is off by the scalar factor $\tfrac{\pi}{2}$.
\begin{corollary}
Let $N\colon T^2\to S^2$ be a harmonic torus whose energy satisfies $E(N)<4\pi$.
Then $N$ is a homomorphism into $S^1\subset S^2$.
\end{corollary}
\begin{proof}
We will show that under our assumptions the spectral genus must be zero, which 
implies that $N$ is given in terms of trigonometric functions 
and thus is a homomorphism \cite{Hitchin}. From 
Theorem~\ref{thm:energy&area_estimates} we see that the spectral genus can at most be one. 
In this case the spectral curve $\Sigma$ has at most one pair of branch points,
not counting $0$ and $\infty$. But then we only need to pass to a double cover
\eqref{eq:H_mu} of $T^2$ to obtain global holomorphic sections from the
two branch values. We therefore obtain the improved estimate 
$
E(N)\geq 4\pi
$
for harmonic tori of spectral genus one.
\end{proof}
Our final application of the estimates in Theorem~\ref{thm:energy&area_estimates}
concerns CMC tori in the $3$-sphere $S^3$ and hyperbolic space $H^3$. 
If $V$ is a quaternionic line bundle with complex structure $S$ and Willmore connection
$\nabla$, we have the family $\nabla_{\mu}$ of flat $\text{\bf Sl}(2,\C)$ connections
\eqref{eq:mu_family} which restrict to $\text{\bf SU}(2)$ connections along the
unit circle $\mu\in S^1$. Recall that for $\mu=e^{\theta I}\in S^1$, this family is 
given by the $S^1$-family of flat quaternionic connections
\[
\nabla_{\lambda}=\nabla+(\lambda-1)A\,,\qquad \lambda=e^{\theta S}
\]
leaving a quaternionic hermitian form \eqref{eq:quat_herm_form} on $V$ invariant.

The CMC surfaces in the $3$-sphere $S^3$ arise as the gauge transformations between two
connections in the family $\nabla_{\lambda}$. 
Assume that $\nabla$ and $\nabla_{\lambda}$, for some $\lambda\neq 1$,
are trivial connections and that we have trivialized $V=M\times \H$ by the connection $\nabla$. 
Since flat connections with the same holonomy are gauge equivalent,
there is a smooth map
\[
f\colon M\to S^3\subset \H
\]
satisfying
\[
f^{-1}df=(\lambda-1)A\,.
\]
If $A$ has no zeros, which we have shown \eqref{eq:deg_eq} to be equivalent to 
$\deg V=1-g$, the map 
$f$ is an immersion. From $*A=SA$ we see that $f$ is conformal. Denoting the
Maurer-Cartan form of $f$ by $\alpha=f^{-1}df$, we easily compute that
\[
d*\alpha+\cot\tfrac{\theta}{2}\,\alpha\wedge\alpha=0
\]
by using $d^{\nabla}*A=0$. Therefore $f$ is a CMC immersion with $H=\cot\tfrac{\theta}{2}$.
As in the case of CMC surfaces in $\R^3$, the 
line bundle $V$ with holomorphic structure $\nabla''$ is the pullback by
$f\colon M\to S^3\subset \H\subset \H\P^1$ of the dual tautological bundle over $\H\P^1$.
By Example~\ref{ex:standard_1} the Willmore energy of $\nabla''$ is 
\[
W(\nabla'')=W(f)=\int \tilde{H}^2-K-K^{\perp}\,,
\]
where $\tilde{H}$ denotes the mean curvature vector of $f$ in $4$-space. Since $f$ takes
values in the $3$-sphere, $\tilde{H}^2=H^2+1$ and the normal bundle is flat, i.e., 
$K^{\perp}=0$. For a CMC torus $f\colon T^2\to S^3$ in the $3$-sphere, we therefore obtain
\[
W(f)=W(\nabla'')=\int H^2+1=(H^2+1)\,\text{area}(f)\,.
\]
The CMC surfaces with $H>1$ in hyperbolic $3$-space $H^3$ arise
from trivial connections $\nabla_{\mu}$ where $\mu$ is off the unit circle.
Recall that given a $2$-dimensional complex vector space with determinant, 
hyperbolic $3$-space $H^3$ consists of all positive definite hermitian forms of determinant one.
Assume now that $\nabla$ and $\nabla_{\mu}$ are trivial for some $\mu$ with $|\mu|\neq 1$.
We again trivialize $V=M\times \H$ by the connection $\nabla$. The decomposition
\eqref{eq:form_decomposition} of the quaternionic hermitian form $<\,,\,>$ 
provides us with a determinant form and a complex hermitian form $(\,,\,)$,
both constant in the given trivialization. But viewing $(\,,\,)$ in the trivialization
given by $\nabla_{\mu}$, we obtain a smooth map 
\[
f\colon M\to H^3
\]
into hyperbolic $3$-space. If $F\colon M\to \text{\bf Sl}(2,\C)$ 
gauges $\nabla$ to $\nabla_{\mu}$, then
\[
f=F^{*} F\,.
\]
\begin{comment}
Due to the symmetry \eqref{eq:J_symmetry} also $\nabla_{\tfrac{1}{\bar{\mu}}$ is a 
trivial connection gauged from $\nabla$ by $JFJ^{-1}$. It is easy to see that
$F^{*}=J^{-1}F^{-1}J$ so that $f$ can also be represented by
\[
f= J^{-1}F^{-1}JF\,.
\]
\end{comment}
Similar to the calculation in the case of the $3$-sphere, we can show that $f$ is
a CMC immersion of constant mean curvature $H=\coth|\mu|>1$. Furthermore, by considering
hyperbolic space $H^3\subset \H\P^1$, the Willmore energy of $f$ is again given by
\[
W(\nabla'')=W(f)=\int \tilde{H}^2-K-K^{\perp}\,.
\]
If $M=T^2$ is a $2$-torus, this reduces as above to 
\[
W(\nabla'')=W(f)=\int H^2-1=(H^2-1)\,\text{area}(f)\,,
\]
where we have used that $\tilde{H}^2=H^2-1$ for the case of hyperbolic space.

By considerations similar to the ones discussed for CMC surfaces in $\R^3$, we
can show that every CMC surface in $S^3$ or $H^3$, provided $H>1$ for the latter,
is obtained from the above construction. For more details we refer the reader
to the paper by Bobenko \cite{Sasha}. 

We now can reformulate Theorem~\ref{thm:energy&area_estimates} for 
CMC tori in $3$-dimensional space forms. If $f$ is such a CMC torus of spectral
genus $g$, then the corresponding line bundle $V$ with Willmore connection
$\nabla$ has 
\[
h^{0}(V)\geq g+2
\]
for the holomorphic structure $\nabla''$: in the case
of $S^3$ both connections $\nabla$ and $\nabla_{\lambda}$, where  $\lambda\in S^1$ and
$\lambda\neq 1$, are trivial. This adds  two more independent holomorphic sections
to the already existing $g$ independent holomorphic sections from the branch points
of the spectral curve $\Sigma$. In the case of hyperbolic space $H^3$ the two additional
independent holomorphic sections come from the triviality of $\nabla$ and
$\nabla_{\mu}$, where $\mu$ is off the unit circle and not a branch value for $\Sigma$.
\begin{theorem}
Let $f$ be a CMC torus of spectral genus $g$ into $3$-space $\R^3$, $S^3$ or $H^3$
where in the latter case the mean curvature $H>1$. Then the Willmore energy of $f$
satisfies the estimate  
\[
W(f)\geq
\begin{cases}
\tfrac{\pi}{4}(g+2)^2 & \,\,\text{if $g$ is even, and}\\
\tfrac{\pi}{4}((g+2)^2-1)&\,\,\text{if $g$ is odd}\,. 
\end{cases}
\]
Since $W(f)=(H^2+\epsilon)\text{\em area}(f)$, where $\epsilon=0,\pm 1$ indicates the
curvature of $3$-space, these estimates translate into area estimates. 
In particular, for 
minimal tori in $S^3$ we have $W(f)=\text{area}(f)$.
\end{theorem}
As an immediate consequences of these estimates we see that
minimal tori in $S^3$ of spectral genus $g\geq 4$ have
$W\geq 9\pi>2\pi^2$. In order to verify  the Willmore conjecture on minimal $2$-tori in the 
$3$-sphere it therefore suffices to consider mimimal tori of spectral genus at most $3$.

A similar remark applies to the Lawson conjecture which asserts that the only
embedded minimal torus in the $3$-sphere is the Clifford torus. It has been
shown \cite{Chinese} that a minimal torus in $S^3$ has area less then $16\pi$.
Therefore our estimates imply that it ``suffices" to 
check minimal tori resulting from spectral genera at most five for a verification
of Lawson's conjecture.
   
%%% Local Variables: 
%%% mode: latex
%%% TeX-master: "willmore"
%%% End: 

\end{document}